\def\status{n} 

\newcommand\papertitle{An intersection formula for CM cycles on Lubin-Tate spaces}
\def\href#1{}
\documentclass[12pt]{amsart}

\usepackage{pdfsync}
\usepackage{fullpage}
\usepackage{python}
\usepackage{mathrsfs}
\usepackage{latexsym, amssymb, amsmath, amscd, amsthm, amsxtra}
\usepackage{eucal}
\usepackage[utf8]{inputenc}
\usepackage{verbatim}
\usepackage{hyperref}
\usepackage[english]{babel}
\usepackage{diagbox}
\usepackage{mathtools}
\usepackage{todonotes} 
\usepackage{cite}
\usepackage{graphicx}
\usepackage{amssymb}
\usepackage{epstopdf}
\usepackage{tikz}
\usetikzlibrary{calc}
\usetikzlibrary{graphs}
\usepackage[all]{xy}
\usepackage{color}
\usepackage{fancyhdr}
\usepackage{arydshln}
\usepackage{times}
\usepackage{datetime}
\usepackage{scrtime}
\usepackage{CJKutf8}
\usepackage[refpage]{nomencl}
\usepackage{cancel}
\usepackage{ulem}
\usepackage{bbm}

\long\def\[#1]#2{ \begin{#1}#2\end{#1}}
\renewcommand{\|}{\big|}
\renewcommand{\(}{\left(}
\renewcommand{\)}{\right)}
\renewcommand{\[}{\left[}

\newcommand{\height}{\mathrm{Height}}
\newcommand{\id}{\mathrm{id}}
\newcommand{\matt}[3]{\mathrm{Mat}_{#1\times#2}(#3)}
\newcommand{\ZZ}{\mathbb{Z}}
\newcommand{\QQ}{\mathbb{Q}}

\newcommand{\map}[3]{#1:#2\longrightarrow #3}
\newcommand{\mapp}[2]{#1\longrightarrow #2}

\newcommand{\maps}[5]{\begin{array}{llrll}
#1&:&#2&\longrightarrow &#3\\\\
&&#4&\longmapsto &#5\\
\end{array}}
\newcommand{\mapi}[4]
{\begin{array}{llrll}
&#1&\longrightarrow &#2\\\\
&#3&\longmapsto &#4\\
\end{array}}
\DeclareMathOperator{\Ker}{\mathrm{Ker}}

\newcommand{\CC}{\mathbb{C}}

\newcommand{\dd}{\mathrm{d}}

\newcommand{\Gal}{\mathrm{Gal}}
\newcommand{\CFF}[1]{\overline{\FF_#1}}
\newcommand{\FF}{\mathbb{F}}

\newcommand{\0}{\textbf{0}}
\newcommand{\one}{\mathbbm{1}}
\newcommand{\OO}[1]{\mathcal{O}_{#1}}
\newcommand{\CO}[1]{\mathcal{O}_{\breve{#1}}}

\newcommand{\vv}{\textbf{v}}
\newcommand{\ww}{\textbf{w}}

\newcommand{\ep}{\boldsymbol{\epsilon}}

\DeclareMathOperator{\End}{\mathrm{End}}
\newcommand{\Vol}{\mathrm{Vol}}
\DeclareMathOperator{\Hom}{\mathrm{Hom}}
\DeclareMathOperator{\coker}{\mathrm{coker}}
\DeclareMathOperator{\Spf}{\mathrm{Spf}\text{ }}
\DeclareMathOperator{\Spec}{\mathrm{Spec}}
\newcommand{\Supp}{\mathrm{Supp}}

\newcommand{\Isom}{\mathrm{Isom}}

\newcommand{\Isog}{\mathrm{Isog}}

\DeclareMathOperator{\length}{\mathrm{length}}
\DeclareMathOperator{\Aut}{\mathrm{Aut}}
\DeclareMathOperator{\gl}{\mathrm{M}}

\DeclareMathOperator{\GL}{\mathrm{GL}}

\DeclareMathOperator{\im}{\mathrm{Im}}

\newcommand{\shexkk}[5]{\xymatrix{0\ar[rr]&&{#1}\ar[rrr]^{#4}&&&{#2}\ar[rrr]^{#5}&&&{#3}\ar[rr]&&0}}
\newcommand{\cate}[8]{
\xymatrix{
{#1}\ar[rrr]^{#5}\ar[d]_{#6}&&&{#2}\ar[d]^{#7}\\
{#3}\ar[rrr]^{#8}&&&{#4}\\}}

\newcommand{\triso}[6]{
\xymatrix{
	{#1}\ar[rrrr]^{#4}_{\cong}\ar[rrd]_{#5}&&&&{#2}\ar[lld]^{#6}\\
	&&{#3}\\
}
}

\newcommand{\trii}[6]{
\xymatrix{	
	&&{#3}\ar[rrd]^{#5}\ar[lld]_{#6}\\
	{#1}\ar[rrrr]^{#4}&&&&{#2}\\
}
}

\newcommand{\ezer}[4]{\xymatrix{
{#1}\ar@<-.5ex>[rr]_{#3}\ar@<.5ex>[rr]^{#4}&&{#2}}
}

\makeatletter
\def\t{\@tt|\em\@gobble\@ttt}
\def\em\@tt#1\em{\@tt|\em}
\def\cm\@tt#1\em#2\@ttt{
\begin{tabular}{#1}
\hline
#2
\\\hline
\end{tabular}\@gobble
}
\def\@tt#1\@ttt#2{
\@ifnextchar,{\@tt #1 & #2\\\hline\expandafter\expandafter\expandafter\@gobble\expandafter\@gobble\@gobble\@ttt}{
\@ifnextchar.{\cm\@tt|l #1 & #2\@ttt}{\@tt|l #1 & #2\@ttt}
}
}

\def\m{\@tut\@gobble\@tutt}
\def\@tut#1\@tutt#2{
\@ifnextchar,{\@tut #1 & #2\\\expandafter\expandafter\expandafter\@gobble\expandafter\@gobble\@gobble\@tutt}{
\@ifnextchar.{\begin{pmatrix}#1 &#2\end{pmatrix}\@gobble}{\@tut #1 & #2\@tutt}
}
}

\def\bm{\@tuut\@gobble\@tuutt}
\def\@tuut#1\@tuutt#2{
\@ifnextchar,{\@tuut #1 & #2\cr\expandafter\expandafter\expandafter\@gobble\expandafter\@gobble\@gobble\@tuutt}{
\@ifnextchar.{\bordermatrix{#1 &#2\cr}\@gobble}{\@tuut #1 & #2\@tuutt}
}
}

\makeatother

\newcommand{\tr}{\mathrm{tr}}

\newcommand{\mm}[4]{\left(\begin{array} {cc}
#1 & #2\\
#3 & #4\\
\end{array}
\right)}

\newcommand{\rr}[2]{\left(\begin{array} {cc}
#1 & #2\\
\end{array}
\right)}

\newcommand{\cc}[2]{\begin{pmatrix}
#1\\
#2\\
\end{pmatrix}}

\long\def\[#1]#2{ \begin{#1}#2\end{#1}}

\begin{document}


\let\yourlabel=\label
\def\lb#1{\yourlabel{#1}\comments{{\color{blue}\textbf{Lable}}:\textbf{#1   }  }}
\def\test{t}
\def\comments#1{\ifthenelse{\equal\status\test}
{{\color{red} #1}}
{}
}

\long\def\s#1{\ifthenelse{\equal\status\test}
{{\color{blue} \textbf{{\color{red}]}#1 \color{red}[}}}
{}
}

\long\def\ss#1{\ifthenelse{\equal\status\test}
{{\color{purple} \textbf{{\color{red}\}}#1 \color{red}\{~\\}}}
{}
}

\newcommand{\commenta}[1]{{\color{blue}(\textbf{Comments:} #1)}}
\newcommand{\ada}[1]{#1}
\newcommand{\dla}[1]{}
\newcommand{\ria}[1]{}
\newcommand{\xo}[1]{\input{Content/#1}}
\newcommand{\oo}[1]{\input{Content/#1}}
\newcommand{\ox}[1]{}


\newtheorem{thm}{Theorem}[section]
\newtheorem{pof}{Proof}[section]
\newtheorem{defi}[thm]{Definition}
\newtheorem{prop}[thm]{Proposition}
\newtheorem{cor}[thm]{Corollary}
\newtheorem{summ}{Summary}
\newtheorem{rem}[thm]{Remark}
\newtheorem{ex}[thm]{Problem}
\newtheorem{lem}[thm]{Lemma}
\newtheorem{exa}[thm]{Example}
\newtheorem{proj}{Project}
\newtheorem{conj}{Conjecture}
\numberwithin{equation}{section}
\setcounter{tocdepth}{1}    
\hypersetup{
    colorlinks=true, 
    linktoc=all,     
    linkcolor=red,  
}

\author{Qirui Li}
\title{\papertitle}
\date{}
\maketitle


\newcommand{\haofan}{k}
\newcommand{\papa}{{\pa nm}}
\newcommand{\pia}{\bpb m\ACj{\ACg_0}}
\newcommand{\papapa}{{\pao nm}}
\newcommand{\dpapapa}{{\deg\left(\xymatrix{\DEf{1,m+n}\ar[r]&\DEf{1,n}}\right)}}
\newcommand{\papapapa}{{\pat nm}}
\newcommand{\dpapapapa}{{\deg\left(\xymatrix{\DEf{2,m+n}\ar[r]&\DEf{2,n}}\right)}}
\newcommand{\pbpb}{{\pb nm}}
\newcommand{\dpbpb}{\deg\left(\xymatrix{\DEF{m+n}\ar[r]&\DEF{n}}\right)}
\newcommand{\dgpa}[1]{\deg\left(\xymatrix{\GKH\ar[r]^{\pi^{#1}}&\GKH}\right)}
\newcommand{\dpa}[2]{\deg\left(\xymatrix{\DEf{#1+#2}\ar[r]&\DEf{#1}}\right)}
\newcommand{\dgpb}[1]{\mathrm{deg}(\beta_\infty^{\infty+#1})}
\newcommand{\inter}[2]{\chi({#1}\otimes^{\mathbb{L}}{#2})}
\newcommand{\q}[3]{\newcommand{#1}{#2}
}
\newcommand{\sa}{^{\prime\prime}}
\newcommand{\Deltai}[1]{\Delta_{#1}}
\newcommand{\Deltao}{\Delta_1}
\newcommand{\Deltat}{\Delta_2}
\q{\Nrdd}{\Nrd_{D}}{Nrdd}
\q{\Nrdg}{\Nrd_{\mathbb{G}}}{Nrdg}
\q{\ef}{\ep_F}{ef}
\q{\GG}{\mathcal{G}}{G}
\q{\ek}{\ep_K}{ek}
\newcommand{\equi}[3]{\mathrm{Equi}_{#3}(#1/#2)}
\newcommand{\Fn}{\mathcal F_n}
\newcommand{\Fm}{\mathcal F_m}
\newcommand{\Fmv}{\mathcal{H}_{m+v}}
\newcommand{\Fnvv}{\mathcal{H}_{n+v}}
\newcommand{\pmn}{\pp{m-n}}
\newcommand{\intt}[2]{\chi\left({#1}\litimes{#2}\right)}
\newcommand{\inttt}[2]{\chi\left({#1}\litimes{#2}\right)}
\newcommand{\inttn}[2]{\intt{#1}{#2}_n}
\newcommand{\intti}[2]{\intt{#1}{#2}_\infty}
\newcommand{\jiao}[2]{\mathrm{Int}(#1,#2)}
\newcommand{\ji}[1]{\mathrm{Int}(#1)}
\newcommand{\inta}[1]{\ca{\biso_1}{\blin_1}{#1}\bullet\ca{\biso_2}{\blin_2}{#1}}
\newcommand{\hinta}[1]{\ca{\biso_1}{\blin_1}{#1}\bullet\pp{v}_*\pro^*\ca{\biso_2}{\blin_2}{#1}}
\newcommand{\pcs}[1]{\pi^{#1}\OO F^\times\oplus\mu\OO F^\times}
\newcommand{\litimes}{\otimes^{\mathbb L}}
\newcommand{\HKO}{\mathcal H_{K_1}}
\newcommand{\HKT}{\mathcal H_{K_2}}
\newcommand{\HFF}{\mathcal H_{F}}
\q{\GN}{\mathcal F_\infty^{\gpp}}{GN}
\q{\GO}{\delta_{1,\infty}}{GO}
\q{\GT}{\delta_{2,\infty}}{GT}
\q{\GI}{\delta_\infty}{GI}
\q{\ev}{\N{x\ACj^{-1}y-\pi^v\ACj^{-1}}}{ev}
\q{\pc}{\pi\OO F\oplus\mu\OO F^\times}{pc}
\q{\SPF}{\Spf\CFF q}{SPF}
\q{\fiber}{\GGH{\pi^{v}\gamma^{-1}}gh}{fiber}
\q{\GF}{\GG_F}{GF}
\newcommand{\GFA}{\GG_{F/A}}
\newcommand{\GKA}{\GG_{K/A}}
\q{\G}{\GG_F^{2h}}{G}
\q{\GH}{\GG_F^h}{GH}
\q{\GHO}{\GG_{F,1}^h}{GHO}
\q{\GHT}{\GG_{F,2}^h}{GHT}
\q{\GFK}{\GG_F^{kh}}{GFK}
\q{\GK}{\GG_K}{GK}
\q{\GKB}{\GG_{\overline{K}}}{GKB}
\q{\GKH}{\GG_K^{h}}{GKH}
\q{\superGKH}{\GG_K^{h}\times\overline{\GG_K^h}}{superGKH}
\q{\HG}{\mathrm{Hom}_{\OO F}(\underline{\OO F^{kh\vee}},\GF)}{HG}
\q{\HGK}{\mathrm{Hom}_{\OO K}(\underline{\OO K^{h\vee}},\GK)}{HGK}
\q{\TG}{\GF\otimes_{\OO F}\OO F^{kh}}{TG}
\q{\TGK}{\GK\otimes_{\OO K}\OO K^{h}}{TGK}

\newcommand{\bbpa}[4]{\theta_{#1}^{#1+#2}(#3,#4)}
\newcommand{\pa}[2]{\theta_{#1}^{#1+#2}}
\newcommand{\pao}[2]{\theta_{1,#1}^{#1+#2}}
\newcommand{\pat}[2]{\theta_{2,#1}^{#1+#2}}
\newcommand{\bbpb}[4]{\beta_{#1}^{#1+#2}(#3,#4)}
\newcommand{\pb}[2]{\beta_{#1}^{#1+#2}}
\newcommand{\dpb}[2]{\mathrm{deg}(\alpha_{#1}^{#1+#2})}
\newcommand{\bbpc}[4]{\eta_{#1}^{#1+#2}(#3,#4)}
\newcommand{\bccc}[4]{\eta_{#1}^{#2}(#3,#4)}
\newcommand{\ppc}[2]{\eta_{#1}^{#1+#2}}

\newcommand{\bpa}[3]{\bbpa n{#1}#2#3}
\newcommand{\ppa}{\pa nm}
\newcommand{\gpa}{\pa {\infty}m}
\newcommand{\Gpa}{[\pi^M]_{\G}}

\newcommand{\bpb}[3]{\bbpb n{#1}{#2}{#3}}
\newcommand{\ppb}{\pb nm}
\newcommand{\bppb}[2]{\bbpb{#1}{#2}\ACj\ACg}
\newcommand{\pppb}{\beta^{\bullet\sim}_{\bullet\sim}(\ACj,\ACg)}
\newcommand{\ppbb}[2]{\beta^{#2}_{#1}(\ACj,\ACg)}
\newcommand{\bpbb}[1]{\bpb{#1}\ACj\ACg}

\newcommand{\bpc}[3]{\bbpc n{#1}{#2}{#3}}
\newcommand{\gpc}[3]{\bbpc \infty{#1}{#2}{#3}}
\newcommand{\bppc}[2]{\bbpc{#1}{#2}\biso\blin}
\newcommand{\bco}[2]{\bbpc{#1}{#2}{\biso_1}{\blin_1}}
\newcommand{\bct}[2]{\bbpc{#1}{#2}{\biso_2}{\blin_2}}
\newcommand{\bppco}[2]{\bbpc{#1}{#2_1}{\biso_1}{\blin_1}}
\newcommand{\bppct}[2]{\bbpc{#1}{#2_2}{\biso_2}{\blin_2}}
\newcommand{\bppch}[2]{\eta_{#1-#2_2}^{#1+#2_1}(\biso_3,\blin_3)}

\newcommand{\pppc}{\eta^{\bullet\sim}_{\bullet\sim}(\biso,\blin)}
\newcommand{\pppcc}{\eta^{\bullet\sim}_{\bullet\sim}(\biso,\id)}
\newcommand{\ppcc}[2]{\eta^{#1+#2}_{#1}(\biso,\blin)}
\newcommand{\pccc}[4]{\eta^{#1+#2}_{#1}(\biso,\blin)}
\newcommand{\bpcc}[1]{\bpc{#1}\biso\blin}

\newcommand{\bppd}[2]{\bbpc{#1}{#2}{\ACj\cdot\biso}{\ACg\cdot\blin}}
\newcommand{\bpdd}[1]{\bbpc n{#1}{\ACj\cdot\biso}{\ACg\cdot\blin}}

\q{\bpp}{\bbpc nm\biso\blin}{bpp}
\q{\bppo}{\bppco nm}{bppo}
\q{\bpph}{\bppch nm}{bpph}
\q{\bppt}{\bppct nm}{bppt}

\newcommand{\ssc}[2]{i_{#2+#1}^{#2}(\biso,\blin)}
\newcommand{\sgc}[2]{i_{\infty}^{#2}(\biso,\blin)}
\newcommand{\sbpb}[1]{i_{m+#1}^m(\ACj,\ACg)}
\newcommand{\sbpp}{\sbpb{n}}
\newcommand{\sbpc}[1]{i_{m+#1}^m(\biso,\blin)}
\newcommand{\sbppc}{\sbpc{n}}

\newcommand{\sss}{{i_{m+n}^m}}
\newcommand{\sssss}{{i_{n}^M}}
\newcommand{\ssssss}{{i_{\infty}^M}}
\newcommand{\ssss}{{i_\infty^m}}

\q{\bp}{(\biso,\blin)}{bp}

\q{\bpn}{(\biso_0,\blin_0)}{bpn}
\newcommand{\gppp}[3]{\eta_\infty^{\infty+#3}(#1,#2)}
\newcommand{\gppb}[3]{\beta_\infty^{\infty+#3}(#1,#2)}
\newcommand{\gppbj}[2]{\gppb{\ACj_{#1}}{\ACg_{#1}}{#2}}
\newcommand{\gpb}{\gppbj{}m}

\newcommand{\sgpb}{i_\infty^m(\ACj,\ACg)}

\newcommand{\gppi}[1]{\gppp{\biso_{#1}}{\blin_{#1}}{}}
\newcommand{\gppj}[2]{\gppp{\biso_{#1}}{\blin_{#1}}{#2}}
\q{\gp}{\gppp{\biso}{\blin}{}}{gp}
\q{\gpo}{\gppp{\biso_1}{\blin_1}{m_1}}{gpo}
\q{\gpt}{\gppp{\biso_2}{\blin_2}{m_2}}{gpt}
\q{\gpp}{{\gppj{}m}}{gpp}
\q{\sgpp}{i_\infty^m(\biso,\blin)}{sgpp}
\q{\gppo}{{\gppj1m}}{gppo}
\q{\gppt}{{\gppj2m}}{gppt}

\newcommand{\gppdj}[2]{\gppp{\ACj\cdot\biso_{#1}}{\ACg\cdot\blin_{#1}}{#2}}
\newcommand{\gc}{Z}
\newcommand{\ca}[3]{Z_{#3}(#1,#2)}
\newcommand{\cao}[1]{\ca{\biso_1}{\blin_1}{#1}}
\newcommand{\cat}[1]{\ca{\biso_2}{\blin_2}{#1}}
\newcommand{\caon}{\cao{m+n}}
\newcommand{\catn}{\cat{m+n}}
\newcommand{\cycy}[1]{\ca{\biso}{\blin}{#1}}
\newcommand{\cycw}[2]{\ca{\biso}{\pi^{#2}\blin}{#1}}
\newcommand{\cycn}[2]{\ca{\biso_{#1}}{\blin_{#1}}{#2}}
\newcommand{\cycwn}[3]{\ca{\biso_{#1}}{\pi^{#2}\blin_{#1}}{#3}}
\q{\cyio}{\ca{\biso_1}{\blin_1}{\infty}}{cyio}
\q{\cyit}{\ca{\biso_2}{\blin_2}{\infty}}{cyit}
\q{\cyc}{\cycy{n}}{cyc}
\q{\cycj}{\gc^{(j)}_{n}(\biso,\blin)}{cycj}
\q{\cyco}{\cycn1n}{cyco}
\q{\cyct}{\cycn2n}{cyct}
\q{\cywo}{\ca{\biso_1}{\pi^w\blin_1}{\infty}}{cywo}
\q{\cywt}{\ca{\biso_2}{\pi^w\blin_2}{\infty}}{cywt}
\q{\cy}{\cycy{\infty}}{cy}

\q{\ACj}{\gamma}{ACj}
\q{\ACx}{x}{ACx}
\q{\GE}{\GG_\eta}{GE}
\q{\GjE}{\GG_{\ACj_c\eta}}{GjE}
\q{\GEN}{{\GG_0^\prime}}{GEN}
\q{\QP}{Quotient Projection}{QP}
\q{\RP}{Restriction Projection}{RP}
\q{\LatF}{\OO F^{kh}}{LatF}
\q{\LatK}{\OO K^h}{LatK}
\q{\ACg}{g}{ACg}
\q{\ACpair}{\AC{\ACj}{\ACg}}{ACpair}
\newcommand{\pp}[1]{\pbpb}
\newcommand{\pg}[1]{(\pi^{#1})}
\q{\pro}{\bpb m\ACj{\ACg_0}}{pro}
\q{\proi}{(\ACj,\ACg)_{\infty}}{proi}
\q{\proh}{({\ACj,\pi^vh\ACg})}{proh}
\q{\pron}{(\ACj_0,\pi^v\ACg_0)}{pron}
\q{\proc}{(\ACg\ACg_0,\ACj_0\ACj)}{proc}
\q{\prng}{(\ACj_0,\ACg)}{prng}
\q{\ACpairnot}{\AC{\gamma_0}{g_0}}{ACpairnot}
\q{\EndMF}{\OO {{D_F}}^{\times}\times \GL_h(\OO F)}{EndMF}
\q{\EndD}{\OO D}{EndD}
\q{\EndM}{\MM_{kh}(\OO F)}{EndM}
\q{\GLMF}{\GL_h(F)\cap\EndM}{GLMF}
\q{\dm}{\DEFf m{}h}{dm}
\q{\dwm}{\DEFf {m+w}{\eta_0}h}{dwm}
\q{\Dm}{\DEFF mF{kh}}{Dm}
\q{\djm}{\DEFf m{\ACj_c\eta}h}{djm}
\q{\dqm}{{\DEFf m{}h}^{(w)}}{dqm}
\newcommand{\Qua}[2]{[#1,#2]}
\newcommand{\infQua}[2]{(#1,#2)_{\infty}}
\newcommand{\Quaa}[1]{[\biso_{#1},\blin_{#1}]}
\newcommand{\Quasii}[1]{\Quasi_{#1}}
\q{\Quasi}{{(\pi^{m}\biso,\blin)}}{Quasi}
\q{\Quasio}{[\biso_1,\pi^{w_1}\blin_1]}{Quasio}
\q{\Quasit}{[\biso_2,\pi^{w_2}\blin_2]}{Quasio}
\q{\QQuasi}{\overline{\Quasi}}{QQuasi}
\q{\QQQuasi}{\overline{\Quasi}}{QQQuasi}
\q{\Quasin}{{[\biso_0,\blin_0]}}{Quasin}
\q{\superQuasi}{\infQua{\eta}{\pi^w(\tau\times\overline{\tau})}}{superquasi}

\newcommand{\TIMES}[1]{\!\!\!\!\underset{#1}{\times}\!\!\!\!}

\q{\HBN}{\pgDEFf{m+v+u+w}\eta hw}{HBN}
\q{\FBN}{haha}{FBN}
\q{\HSN}{\pgDEFf{m+u+w}\eta hw}{HSN}
\q{\FSN}{\GEN^h[\pi^{m+w}]}{FSN}
\q{\HBM}{\DEFF{m+u+v}F{kh}}{HBM}
\q{\FBM}{\G[{\pi^{m+v}\ACj^\prime}^{-1}\ACj^{-1}\ACg^\prime\ACg]}{FBM}
\q{\HSM}{\DEFF{m+u}F{2h}}{HSM}
\q{\FSM}{\G[{\pi^m\ACj^\prime}^{-1}\ACg^\prime]}{FSM}
\q{\RFSM}{\G[\pi^m\ACj\ACg^{-1}]}{RFSM}
\q{\BM}{\DEFF uF{2h}}{BM}
\q{\BN}{\DEFf{u}{\eta_0} h}{BN}

\q{\Mat}{M_{\eta,\tau}}{Mat}
\q{\Diso}{\Delta(\biso,\blin)}{Disc}
\newcommand{\Disoo}{\Delta(\biso_1,\blin_1)}
\newcommand{\Disot}{\Delta(\biso_2,\blin_2)}
\q{\Discn}{\Delta_{\eta_0,\tau_0}}{Discn}
\q{\Up}{\Gamma^{-1}}{Up}
\q{\Upn}{\Gamma^`}{Upn}
\q{\embed}{(\eta,\tau)}{embed}
\q{\MG}{\mathfrak{m}_{\G}}{MG}
\q{\MGEN}{\mathfrak{m}_{\GEN}}{MGEN}
\q{\RG}{\CFF q[[X_1,\cdots,X_{2h}]]}{RG}
\q{\RGEN}{\CFF q[[T_1,\cdots,T_{h}]]}{RGEN}

\q{\INTER}{\mathcal{F}_{m+v}\otimes{\pron^*\mathcal{F}_{m}}}{INTER}
\q{\SINTER}{s^*\mathcal{F}_{m+v}\otimes s^*{\pron^*\mathcal{F}_{m}}}{SINTER}
\q{\INTERNUM}{\inter{\mathcal{F}_{m+v}}{\pron^*\mathcal{F}_{m}}}{INTERNUM}
\q{\Tor}{\mathcal{T}or(\mathcal{F}_{m+v},\pron^*\mathcal{F}_{m})}{Tor}
\q{\INFINT}{\infint{\gamma_0}{\tau_0}{\pi^v\gamma_0^{-1}}{g_0}}{INFINT}
\q{\HM}{\mathcal{H}_m}{HM}
\q{\FM}{\mathcal{F}_{m+u}}{FM}
\q{\HV}{\mathcal{H}_{m+v}}{HV}
\q{\FV}{\mathcal{F}_{m+v+u}}{FV}

\newcommand\len{\length_{\CO F}}
\newcommand\mor{\mathrm{Mor}_{\CO F}}

\newcommand\gro[1]{$\QQ$-vector space that is freely generated by coherent sheaves of #1. }
\newcommand\groo[1]{\widetilde{K}\left(#1\right)\otimes_{\ZZ}\QQ}
\newcommand\graa[1]{\widetilde{K}\left(#1\right)}

\newcommand{\nmb}{\varrho}
\newcommand{\ff}{\mathrm{V}}
\newcommand{\fff}{\alpha}
\newcommand{\ffff}{\beta}

\newcommand{\Heegner}{CM }
\newcommand{\yaba}{|\mathrm{Disc}_{K/F}|_F}
\newcommand{\yabi}{|\mathrm{Disc}_{K/F}|_F^{-h^2}}
\newcommand{\yapi}{|\mathrm{Disc}_{K_1\otimes_F L/L}|_L^{-h^2}}
\newcommand{\yabii}{|\mathrm{Disc}_{K_1/F}|_F^{-h^2}}

\newcommand{\J}{\GKH[\pi^m\biso\blin]}
\newcommand{\Jw}{\GKH[\pi^{m+w}\biso\blin]}
\newcommand{\Jo}{\GKH[\pi^{m_1}\biso_1\blin_1]}
\newcommand{\Jt}{\GFK[\pi^{m_2}\biso_2\blin_2]}
\newcommand{\Jh}{\GKH[\pi^{m_3}\biso_3\blin_3]}
\newcommand{\JJ}{\GFK[\pi^m\ACj\ACg]}
\newcommand{\JJo}{\GFK[\pi^{m_1}\ACj_1\ACg_1]}
\newcommand{\JJt}{\GFK[\pi^{m_2}\ACj_2\ACg_2]}
\newcommand{\JJJ}{\GFK[\pi^m]}
\newcommand{\JJJJ}{\GF^{2h}[\pi^m]}
\newcommand{\JJJJJ}{\GF^{2h}[\pi^M]}
\newcommand{\CCC}{\mathcal{C}}
\newcommand{\CCCC}{\mathcal{C}\otimes\CFF q}
\def\m{\mathrm{cond}}
\newcommand{\standard}{standard}
\newcommand{\qz}[1]{{\mathcal D^{#1}}}
\newcommand{\bz}[1]{{\mathcal D_{#1}}}
\newcommand{\N}[1]{|\mathrm{Nrd}(#1)|_F^{-1}}
\tikzstyle{every picture}+=[remember picture]
\tikzstyle{na} = [shape=rectangle,inner sep=0pt,text depth=0pt]

\newcommand{\bisoo}{\biso^{-1}}
\newcommand{\blinn}{\blin^{\vee}}
\newcommand{\bbbbiso}{\Isog_{\OO E}(\GG_{E},\GF)}
\newcommand{\bbbblin}{\Isom_E(E^{lkh},F^{kh})}

\newcommand{\bbbiso}{\Isom_{\OO F}(\GK,\GF)}
\newcommand{\bbblin}{\Isom_{\OO F}(\OO K^{h},\OO F^{2h})}
\newcommand{\bilin}{\Isom_{\OO F}(\OO K^{h},\OO F^{kh})}
\newcommand{\bbiso}{\Isog_{\OO F}(\GK,\GF)}
\newcommand{\bblin}{\Isom_F(K^{h},F^{2h})}
\newcommand{\obbiso}{\Isog_{\OO F}(\GG_{K_1},\GF)}
\newcommand{\obblin}{\Isom_F(K_1^h,F^{2h})}
\newcommand{\tbbiso}{\Isog_{\OO F}(\GG_{K_2},\GF)}
\newcommand{\tbblin}{\Isom_F(K_2^h,F^{2h})}
\newcommand{\balin}{\Isom_F(K^h,F^{kh})}
\newcommand{\bibiso}{\Isog_{\OO F}(\GF,\GF)}
\newcommand{\biblin}{\Isom_F(F^{2h},F^{2h})}
\newcommand{\bublin}{\Isom_F(F^{kh},F^{kh})}
\newcommand{\boblin}{\Isom_F(F^2,K)}
\newcommand{\eqh}{\mathrm{Equi}_h(F,K)}
\newcommand{\eqth}{\mathrm{Equi}_{2h}(F,F)}
\newcommand{\heeg}{\mathrm{Heeg}_{h}(F,K)}
\newcommand{\Ltimes}{\otimes^{\mathbb L}}
\newcommand{\ms}[1]{\mm{#1_{11}}{#1_{12}}{#1_{21}}{#1_{22}}}
\newcommand{\mv}{\mm1{}{}{\pi^v}}
\newcommand{\mf}{\mm111{-1}}
\newcommand{\es}{\ep_S}
\newcommand{\er}{\ep_R}
\newcommand{\qqm}[1]{\{#1\}}
\newcommand{\qqsigma}{\sigma}
\newcommand{\qq}[2]{\qqm{#1}+\qqm{#2}\qqsigma}
\newcommand{\qspace}[1]{V_{#1}}
\newcommand{\U}[1]{U_{#1}}
\newcommand{\PD}[3]{\prod_{\vv\in{#3}}(#1-#2(\vv))}
\newcommand{\bs}{\overline{s}}

\newcommand{\biso}{\varphi}
\newcommand{\blin}{\tau}
\newcommand{\infint}[4]{c_{#1,#2}({#3}{#4})}
\newcommand{\CA}[2]{Q_{{#1},{#2}}}
\newcommand{\GGGG}[1]{\mathfrak{G}_{#1}}
\newcommand{\AC}[2]{({#1},{#2})}
\newcommand{\FO}[2]{R_{#1,#2}}
\newcommand{\GGT}[1]{\mathcal{G}[\pi^{#1}]}
\newcommand{\GGS}[1]{\mathcal{G}_0^{#1}}
\newcommand{\GGH}[3]{\GG_0^{#3}[\AC{#1}{#2}]}
\newcommand{\GGSS}[2]{\GG_0[\pi^{#1}]^{#2}}
\newcommand{\AAA}[1]{\mathcal{A}^{(#1)}}
\newcommand{\DEF}[1]{{\mathcal{M}_{#1}}}
\newcommand{\DEf}[1]{{\mathcal{N}_{#1}}}

\newcommand{\OD}[2]{D_{#1}}
\def\delta{Z}

\newcommand{\gth}{\mathrm{GL}_{2h}(F)}
\newcommand{\Nrd}{\mathrm{Nrd}}
\newcommand{\nrd}{\mathrm{nrd}}

\newcommand{\suga}{\Gamma_U}
\newcommand{\uni}{\omega}\newcommand{\rc}[1]{W_{#1}}
\newcommand{\la}{\lambda}
\newcommand{\lo}{a}
\newcommand{\loo}[2]{\lo_{#1}(#2)}
\newcommand{\losj}{\loo{s}{\ACj}}
\newcommand{\losk}{\loo{s}{k}}
\newcommand{\loskj}{\loo{s}{k\ACj}}
\newcommand{\losoj}{\loo{s_1}{\ACj}}
\newcommand{\losjo}{\loo{s}{\ACj_1}}
\newcommand{\losjt}{\loo{s}{\ACj_2}}
\newcommand{\laa}[2]{\la_{#1}(#2)}
\newcommand{\lasj}{\laa{s}{\ACj}}
\newcommand{\lasp}{\laa{s}{\ACj^\prime}}
\newcommand{\lasjo}{\laa{s}{\ACj_1}}
\newcommand{\lasjt}{\laa{s}{\ACj_2}}
\newcommand{\laskj}{\laa{s}{k\ACj}}
\newcommand{\lasoj}{\laa{s_1}{\ACj}}
\newcommand{\lastj}{\laa{s_2}{\ACj}}
\newcommand{\lask}{\laa{s}k}
\newcommand{\lastk}{\laa{s_2}{k}}
\newcommand{\lasok}{\laa{s_1}{k}}
\newcommand{\univ}{\mathrm{univ}}

\newcommand{\tar}[2]{\End(\GG_{#1}\otimes_{\rc {#1}}\rc {#1}/\uni^{#2})}

\newcommand{\AF}{{A_F}}
\newcommand{\AK}{{A_K}}
\newcommand{\AD}{{A_D}}
\renewcommand{\AAA}{{A_?}}
\newcommand{\GD}{{G_D}}
\newcommand{\GGG}{{G_?}}
\newcommand{\HF}{{H_F}}
\newcommand{\HK}{{H_K}}
\newcommand{\HD}{{H_D}}
\newcommand{\HHH}{{H_?}}
\newcommand{\TF}{{T_F}}
\newcommand{\TK}{{T_K}}
\newcommand{\TD}{{T_D}}
\newcommand{\TTT}{{T_?}}
\newcommand{\OOF}{{\mathrm{Orb}_F}}
\newcommand{\OOK}{{\mathrm{Orb}_K}}
\newcommand{\ooF}{{\mathrm{orb}_F}}
\newcommand{\ooK}{{\mathrm{orb}_K}}
\newcommand{\OOD}{{\mathrm{Orb}_D}}
\newcommand{\OOO}{{\mathrm{Orb}_?}}
\newcommand{\OBF}{{\mathbb{O}^{reg}_F}}
\newcommand{\OBK}{{\mathbb{O}^{reg}_K}}
\newcommand{\OBB}{{\mathbb{O}^{reg}_?}}
\newcommand{\Res}{\mathrm{Res}}
\newcommand{\orb}{\mathrm{orb}}

\newcommand\overlina[1]{#1^\sigma}

\newcommand\falg{F^{\mathrm{alg}}}

\newcommand\iii{i_{\biso,\blin}}
\newcommand\iiio{i_{\biso_1,\blin_1}}
\newcommand\iiit{i_{\biso_2,\blin_2}}
\newcommand\ppp{p_{\biso,\blin}}
\newcommand\pppo{p_{\biso_1,\blin_1}}
\newcommand\pppt{p_{\biso_2,\blin_2}}
\newcommand\gaga{(\Gamma_{\blin})_{\GK}}
\newcommand\gagao{(\Gamma_{\blin_1})_{\GG_{K_1}}}
\newcommand\gagat{(\Gamma_{\blin_1})_{\GG_{K_2}}}

\newcommand\ggga{\left(\Gamma_{\blin}^{-1}\right)_{\GK}}
\newcommand\gggao{\left(\Gamma_{\blin_1}^{-1}\right)_{\GG_{K_1}}}
\newcommand\gggat{\left(\Gamma_{\blin_2}^{-1}\right)_{\GG_{K_2}}}

\newcommand\pipa{(P_{\blin})_{\GK}}
\newcommand\pipao{(P_{\blin_1})_{\GG_{K_1}}}
\newcommand\pipat{(P_{\blin_2})_{\GG_{K_2}}}

\newcommand\pppa{\left(P_{\blin}^{-1}\right)_{\GK}}
\newcommand\pppao{\left(P_{\blin_1}^{-1}\right)_{\GG_{K_1}}}
\newcommand\pppat{\left(P_{\blin_2}^{-1}\right)_{\GG_{K_2}}}

\newcommand\qaqa{(Q_{\blin})_{\GK}}
\newcommand\qaqao{(Q_{\blin_1})_{\GG_{K_1}}}
\newcommand\qaqat{(Q_{\blin_2})_{\GG_{K_2}}}

\newcommand\qqqa{\left(Q_{\blin}^{-1}\right)_{\GK}}
\newcommand\qqqao{\left(Q_{\blin_1}^{-1}\right)_{\GG_{K_1}}}
\newcommand\qqqat{\left(Q_{\blin_2}^{-1}\right)_{\GG_{K_2}}}

\newcommand\blan{M_{\blin}}
\newcommand\blano{M_{\blin_1}}
\newcommand\blant{M_{\blin_2}}

\newcommand\ttor{\mathcal{T}or}

\begin{abstract}
We give an explicit formula for the arithmetic intersection number of \Heegner cycles on Lubin-Tate spaces for all levels. We prove our formula by formulating the intersection number on the infinite level. Our \Heegner cycles are constructed by choosing two separable quadratic extensions $K_1,K_2$ over a non-Archimedean local field $F$. Our formula works for all cases: $K_1$ and $K_2$ can be either the same or different, ramified or unramified over $F$. This formula translates the linear Arithmetic Fundamental Lemma (linear AFL) into a comparison of integrals. As an example, we prove the linear AFL for $\GL_2(F)$ in this article.

\end{abstract}
\tableofcontents

\section{Introduction}

\subsection{Global Motivation}
In this article, we study an intersection problem of special cycles on Lubin--Tate towers motivated by the Gross--Zagier formula and its generalizations to certain higher-dimensional Shimura varieties. 

The Gross--Zagier formula \cite{gross1986heegner}\cite{yuan2013gross} relates the Neron–Tate height of Heegner points on Shimura curves to the first central derivative of certain L-functions. There have been conjectural generalizations of the  Gross--Zagier formula to higher dimensional Shimura varieties, notably the Gan--Gross--Prasad conjectures \cite{GGP} and the Kudla--Rapoport conjecture \cite{KRGlobal}. Recently,  Zhang has proposed another generalization  \cite{zhang2017twisted}  from which  the intersection problem in this paper arises.

We briefly recall the geometric  construction in Zhang's proposal.  Let $F_0=\QQ$ (in general, a totally real field), and let $F$ (resp. $E_0$) be an imaginary (resp. a real) quadratic field. Let  $E=E_0\otimes_{F_0} F$, a bi-quadratic field extension of $F_0$.  Let $V$ be an $n$-dimensional vector space over $E$ with an $E/E_0$-Hermitian form $\langle -,-\rangle_V$. Consider the unitary group $H$ (an algebraic group over $E_0$) associated to $V$ and  $\langle -,-\rangle_V$, and denote by  $\mathrm{Res}_{E_0/\QQ} H$ (an algebraic group over $\QQ$) the restriction of scalars of $H$.  Let $V'=V$, viewed as an $F$ vector space (of dimension $2n$), and let 
 $\langle -,-\rangle_{V'}=\mathrm{tr}_{E/F}\langle -,-\rangle_{V}$ be the induced $F/F_0$-Hermitian form.  Let $G$ be the unitary group  associated to $V'$ and  $\langle -,-\rangle_{V'}$. We have a natural embedding of algebraic groups over $\QQ$
 $$
\mathrm{Res}_{E_0/\QQ} H\longrightarrow G.
$$
 For a cuspidal automorphic representation $\pi$ of $G$ and $\phi\in\pi$, we can consider the $H$-period integral of $\phi$. A conjecture in \cite{zhang2020period} states that this period integral is related to the central value of a certain L-function  associated to $\pi$, analogous to the Waldspurger formula and the global Gan--Gross--Prasad conjecture \cite{GGP}. 

We assume further that the signatures of $V$ at the the two archimedean places of $E_0$ are $((n-1,1),(n,0))$. Then the signature of $V'$ at the archimedean place of $F_0=\QQ$ is $(2n-1,1)$.  The Shimura varieties   $\mathrm{Sh}_G$ and $\mathrm{Sh}_H$ associated to $G$ and $\mathrm{Res}_{E_0/\QQ} H$ have dimension $2n-1$ and $n-1$ respectively (see section 27 in \cite{GGP}). The algebraic group inclusion $\mathrm{Res}_{E_0/\QQ} H\longrightarrow G$ induces a cycle $Z$ on $\mathrm{Sh}_G$ of codimension $n$.  Then the height pairing of the $\pi$-isotypic part of the cycle $Z$ is conjecturally related to the first central derivative of a certain L-function associated to $\pi$, analogous to the Gross--Zagier formula and the arithmetic Gan--Gross--Prasad conjecture \cite{GGP}.

In \cite{zhang2012arithmetic}, Zhang proposed a relative trace formula (RTF) approach to the arithmetic Gan–-Gross-–Prasad conjecture. The approach leads to local conjectures, notably the arithmetic fundamental lemma (AFL) conjecture formulated by Zhang in \cite{zhang2012arithmetic} (see \cite{zhang2019arithmetic} for recent progress), and the arithmetic transfer (AT) conjecture formulated by Rapoport--Smithling--Zhang {\cite{rapoport2017arithmetic,rapoport2018regular}}. The AFL conjecture relates the special value of the derivative of a relative orbital integral to an arithmetic intersection number on a Rapoport--Zink (RZ) formal moduli space of p-divisible groups attached to a unitary group. In this paper, we study an analog of Zhang's AFL conjecture in this new setting. More precisely,  let  $v$ be a place of $F_0=\QQ$ that is split in $F$ and inert in $E_0$. Then the local unitary group $G_{F_{0,v}}$ (resp. $H_{E_{0,v}}$) is isomorphic to the general linear group $\GL_{2n, F_{0,v}}$ (resp. $\GL_{n,E_{0,w}}$ for the unique place $w$ of $E_0$ above $v$). Then the corresponding RZ spaces are the Lubin--Tate formal moduli spaces, and we will  call the resulting cycle the {\Heegner cycle} (relative to the quadratic extension $E_{0,w}/F_{0,v}$). We then consider the intersection number of the \Heegner cycle with their translation in the ambient Lubin--Tate space.  Zhang conjectured in his unpublished notes \cite{zhang2017notes} that this intersection number is related to the first derivative of a relative orbital integral. He called this conjecture  {the linear AFL}, to be distinguished from the AFL conjecture \cite{zhang2012arithmetic} in the context of the arithmetic Gan--Gross--Prasad conjecture.  We will recall the precise statement in \S\ref{ss lAFL}.

In this article, we focus on the geometric side of the linear AFL and we will establish an explicit formula for the arithmetic intersection numbers. In fact, we will work in a much more general setting. First of all, our Theorem \ref{SUPER} calculates the intersection number on the Lubin--Tate space with arbitrary Drinfeld level structure, and we allow both ramified and unramified quadratic extensions in characteristic $0$ or an odd prime. It also works in characteristic 2 for separable quadratic extensions. Note that  in the function field case, Yun and Zhang have  discovered a ``higher Gross--Zagier formula" \cite{yun2015shtukas}, relating higher derivatives of L-functions to intersection numbers of special cycles on the moduli space of Drinfeld Shtukas of rank two. Our theorem in the positive characteristic case is related to a generalization of their theorem to the moduli space of Shtukas of higher rank {See Remark 4.1 of \cite{zhang2017periods}}.  Secondly, we obtain results for the intersection of  \Heegner cycles relative to  two different quadratic extensions (see Theorem \ref{BIQ}).  In the function field case, this is related to the recent result of Howard--Shnidman on the Gross--Kohnen--Zagier formula \cite{howard2019gross}.

Our formula can be expressed as an integral over $\GL_{2h}(F)$, where $F$ is some non-Archimedean local field. By using this formula in the case of $h=1$, the author gives a new proof of Gross and Keating’s result on lifting endomorphisms of formal modules in \cite{li2017quasi}. Further, our formula reduces the linear AFL to an analytic identity between two integrals. In \cite{li2017arithmetic}, we prove the identity by direct calculation when $h=2$, and hence verify the linear AFL in these cases. Moreover, we could use our main theorem to verify various analogs of the arithmetic transfer conjecture in this new setting. We will pursue these directions in the future. In the ongoing work of the author and Howard, we conjectured that the linear AFL holds for more general settings, where we allow \Heegner cycles associated with different quadratic extensions  \cite{howard2020bi}.

\subsection{Main Results --- a baby case}\lb{bababa}
We abandon all previous notations and start new ones. Let $h$ be a non-negative integer. Fix a non-Archimedean local field $F$ (of any characteristic) with uniformizer $\pi$, ring of integers $\mathcal O_F$, residue field $\mathbb F_q$, and algebraic closure $F^{\mathrm{alg}}$. Let $K$ be a separable quadratic extension of $F$ with ring of integers $\OO K$ and residue field $\mathbb F_{q_K}$. Let $\zeta_K(s)=(1-q_K^{-s})^{-1}$, $\zeta_F(s)=(1-q^{-s})^{-1}$ be their zeta functions, and let $|\bullet|_F$, $|\bullet|_K$ be their corresponding absolute values, which are normalized so that $|x|_F^{-1}=\#\left(\OO F/x\right)$ and $|y|_K^{-1}=\#\left(\OO K/y\right)$ for any $0\neq x\in\OO F$ and $0\neq y\in\OO K$.  Let $\CFF q$ be the algebraic closure of $\FF_q$. We fix a field embedding $\mathbb F_{q_K}\rightarrow\CFF q$.

Let $\breve K$ (resp. $\breve F$) be the completion of the maximal unramified extension of $K$ (resp. $F$), with ring of integers $\CO K$(resp. $\CO F$). Let $\sigma\in\Gal(K/F)$ be the non-trivial element, and denote by $x^\sigma$ the Galois conjugate of $x\in K$. Let $|\mathrm{Disc}_{K/F}|_F$ be the norm of the discriminant of $K/F$, which equals $|(x-x^\sigma)^2|_F$ for any $x\in K$ such that $\OO K=\OO F[x]$.  Let $\GF$ be a formal $\mathcal O_F$-module of height $2h$ over $\CFF q$. It is well-known that 
$$
D_F:=\End(\GF)\otimes_{\mathcal O_F}F
$$
is a central division algebra of invariant $\frac1{2h}$, which implies $D_F\otimes_FF^{\mathrm{alg}}\cong\matt {2h}{2h}{F^{\mathrm{alg}}}$. We fix an $\OO F$-algebra map 
\[equation]{\lb{suchaway}\OO K\rightarrow \End(\GF)}
in such a way that the induced action of $\OO K$ on the $\CFF q$-vector space $\mathrm{Lie}(\GF)$ is through the reduction map $\OO K\rightarrow \FF_{q_K}\rightarrow\CFF q$, which yields a formal $\OO K$-module $\GK$ of height $h$. 

Denote the Lubin--Tate deformation space of $\GF$ (resp. $\GK$) by $\DEF0$ (resp. $\DEf0$), which admits a natural action of $\Aut(\GF)$ (resp. $\Aut(\GK)$). It is well known that $\DEF0$ is the formal specturm of a power series ring in $2h-1$ variables over $\CO F$. 
Deforming $\GF$ with the extra $\OO K$-module structure through $\OO K\rightarrow \End(\GF)$ gives rise to a morphism of formal schemes $\eta_0^0:\DEf0\rightarrow\DEF0$, which yields a cycle $Z_0:=(\eta_0^0)_*\OO{{\DEf0}}$ with $h$-dimensional support on $\DEF0$. Translating $Z_0$ by an element $\ACj\in\Aut(\GF)$ yields another cycle $Z_0(\ACj)$. Consider the intersection number
\[equation*]{
	\chi(Z_0\litimes_{\DEF0}Z_0(\ACj))=\sum_{i=0}^{\infty}(-1)^i\length_{\CO F}\left(\mathrm{Tor}^i_{\DEF0}(Z_0,Z_0(\ACj))\right),
	}
which is in fact well-defined and only has one non-zero term $\length_{\CO F}(Z_0\otimes_{\DEF0}Z_0(\ACj))$ for almost all elements $\ACj\in D_F^\times$. The main result of this article is a formula identifying $\chi(Z_0\litimes_{\DEF0}Z_0(\ACj))$ with an integral involving invariant polynomials of $\ACj$, which can be described as follows.

Let $K\subset D_F$ be the embedding corresponding to \eqref{suchaway}. Denote by $D_K\subset D_F$ the centralizer of $K$. Note that $D_K\otimes_F\falg\cong \matt hh\falg\times\matt hh\falg$, and that $D_F\otimes_F\falg\cong\matt{2h}{2h}\falg$. Up to conjugation, the embedding 
$D_K\otimes_FF^{\mathrm{alg}}\rightarrow D_F\otimes_FF^{\mathrm{alg}}$
is identified with the blockwise diagonal embedding 
\begin{equation}\lb{ebdd}\matt hh{F^{\mathrm{alg}}}\times\matt hh{F^{\mathrm{alg}}}\rightarrow \matt{2h}{2h}{F^{\mathrm{alg}}},\end{equation}
which associates a degree $h$ polynomial $P_g$ to any element $g\in\GL_{2h}(\falg)$ as described by the following definition.
\begin{defi}[Invariant Polynomial]\label{invpoly}\comments{invpoly}
 For any $A,B,C,D\in\matt hh\falg$ be such that
$$g:=\mm ABCD$$
is an element in $\GL_{2h}(\falg)$, let $X',X''\in\matt hh\falg$ be such that
	\begin{equation}\label{gprime}\comments{gprime}\mm{X'}{}{}{X''}=\[pmatrix]{A&\\&D}\[pmatrix]{A&B\\C&D}^{-1}\[pmatrix]{A&\\&D}\[pmatrix]{A&-B\\-C&D}^{-1}.\end{equation}
		Then $X'$ and $X''$ have the same characteristic polynomial, which we denote by $P_g$. Call $P_g$ the invariant polynomial of $g$ with respect to the embedding $\GL_h(\falg)\times\GL_h(\falg)\subset \GL_{2h}(\falg)$.
\end{defi}
It is straightforward to see that $P_{g'}=P_g$ for any 
$$g'\in \left(\GL_h(\falg)\times\GL_h(\falg)\right)\cdot g\cdot \left(\GL_h(\falg)\times\GL_h(\falg)\right).$$
Therefore, $P_g$ is defined up to the double quotient 
$$\GL_h(\falg)\times\GL_h(\falg)\backslash \GL_{2h}(\falg)/\GL_h(\falg)\times\GL_h(\falg).$$

For any $\ACj\in D_F^\times$, we define the invariant polynomial of $\ACj$ with respect to $D_K^\times\subset D_F^\times$ by regarding it as an element in $D_F\otimes \falg$. The invariant polynomial of $\ACg\in\GL_{2h}(F)$ with respect to the embedding $\GL_h(K)\subset \GL_{2h}(F)$ is defined similarly. In fact, the coefficients of these polynomials are in $F$ (see Proposition \ref{coeffa}).
The main result of this article is the following formula and its generalizations.
\begin{thm}\lb{babys}Let $P_\ACj$ be the invariant polynomial of $\ACj$ with respect to $D_K^\times\subset D_F^\times$. If $P_\ACj$ is irreducible and $P_\ACj(0)P_\ACj(1)\neq 0$, then $\chi(Z_0\litimes_{\DEF 0}Z_0(\ACj))$ is finite and we have the following formula 
$$
	\chi(Z_0\litimes_{\DEF 0}Z_0(\ACj))=\frac{\left(\prod_{i=1}^h\zeta_K(i)\right)^2}{\prod_{i=1}^{2h}\zeta_F(i)}\cdot\yabi\cdot\int_{\GL_{2h}(\mathcal O_F)}\left|\mathrm{Res}(P_\ACj,P_\ACg)\right|_F^{-1}\dd g,
$$
	where $P_\ACg$ is the invariant polynomial of $\ACg\in\GL_{2h}(F)$ with respect to an arbitrary embedding $\GL_{h}(K)\subset\GL_{2h}(F)$ which restricts to $\GL_{h}(\OO K)\subset\GL_{2h}(\OO F)$, and $\mathrm{Res}(P_\ACj,P_\ACg)$ is the resultant of $P_\ACj$ and $P_\ACg$. Here $\dd g$ is the normalized Haar measure of $\GL_{2h}(\OO F)$.
\end{thm}

\subsection{Main results --- generalizations}
Following Drinfeld \cite{drinfel1974elliptic}, there is a projective system of formal schemes
$$
\DEF\bullet = \cdots\rightarrow \DEF n\rightarrow \DEF {n-1}\rightarrow \cdots\rightarrow\DEF 0
$$
such that each member $\DEF n$ parametrizes deformations $\GG$ of $\GF$ with additional data of $\pi^n$-level structures defined by an $\OO F$-linear homomorphism
$$\map{\alpha_n}{\OO F^{2h\vee}/\pi^n\OO F^{2h\vee}}{\GG[\pi^n]}$$
with some non-degeneracy conditions (See Definition \ref{levelstru}). Here by $\OO F^{2h\vee}$ we mean the $\OO F$-module of $1\times 2h$-matrices over $\OO F$.
The action of the group $\GL_{2h}(\OO F)$ on $\pi^n$-level structures induces a surjective homomorphism
$$
\GL_{2h}(\OO F)\longrightarrow \Aut_{\DEF 0}(\DEF n)\cong\GL_{2h}(\OO F/\pi^n)
$$
with the kernel denoted by $R_n$ for any $n\geq 0$ . 
Clearly,
$$
R_0=\GL_{2h}(\OO F);\qquad R_n=I_{2h}+\pi^n\matt {2h}{2h}{\mathcal O_F}
$$
for $n\geq 1$.
A pair
\begin{equation}\label{pairs}\comments{pairs}
\begin{split}
\map{\blin}{K^{h}}{F^{2h}},\\
\map{\biso}{\GK}{\GF}
\end{split}
\end{equation}
consisting of an $\OO F$-quasi-isogeny and an $F$-linear isomorphism induces $F$-algebra embeddings
$$
\varphi\cdot K\cdot \varphi^{-1}\subset \varphi\cdot D_K\cdot \varphi^{-1}\subset D_F;
$$
$$
K\subset \matt hhK\subset\matt{2h}{2h}F.
$$
Let $\mathcal O:=K\cap \matt{2h}{2h}{\OO F}$. Loosely speaking, up to some multiplicities, we define the CM cycle $Z_n(\biso,\blin)$ on $\DEF n$ by deforming $\GF$ with an extra $\mathcal O$-module structure through $\varphi\cdot\mathcal O\cdot\varphi^{-1}\subset \End(\GF)$, and with a $\pi^n$-level structure $\alpha_n$ chosen in such a way that $\mathcal O$ acts on $\alpha_n$ through $\mathcal O\subset K\subset \matt{2h}{2h}F$. Precomposing \eqref{pairs} with two elements $\ACj\in D_F^\times$ and $\ACg\in\GL_{2h}(F)$ yields another cycle $Z_n(\ACj\biso,\ACg\blin)$. Based on the discussion in \S\ref{removal}, the cycle $Z_n(\varphi,\tau)$ (which is also denoted by $Z_n^{(0)}(\varphi,\tau)$ in this article) is either empty or isomorphic to a cycle $Z_n(\varphi',\tau')$ with the pair $(\varphi',\tau')$ normalized in the sense that $\mathrm{Height}(\biso')=\mathrm{Height}(\blin')$. Therefore, without loss of generality, for the rest of the section, we may fix $\biso$ and assume $\mathrm{Height}(\biso)=\mathrm{Height}(\blin)$ (See Definition \ref{heto}). Denote by $P_{\ACj}$ the invariant polynomial of $\ACj$ with respect to $\varphi\cdot D_K\cdot \varphi^{-1}\subset D_F$ for any $\ACj\in D_F^\times$.
Denote by $\nrd(\ACj)$ the reduced norm of $\ACj$. We have the following generalization of our theorem.
Note that the restriction $|\nrd(\ACj)|_F^{-1}=|\det(\ACg)|_F$ is necessary for $(\ACj\biso,\ACg\blin)$ to be an equi-height pair. 
\begin{thm}\lb{girls}
Suppose  $|\nrd(\ACj)|_F^{-1}=|\det(\ACg)|_F$. If $P_\ACj$ is irreducible and $P_\ACj(0)P_\ACj(1)\neq 0$, or if $\mathrm{Res}(P_\ACj,P_{x\ACg})\neq 0$ for any $x\in R_n$, then the intersection number of $Z_n(\biso,\blin)$ and $Z_n(\ACj\cdot\biso,\ACg\cdot\blin)$ is finite, and we have the following formula for the intersection numbers for $n\geq 1$,
$$
	\chi(Z_n(\biso,\blin)\litimes_{\DEF n}Z_n(\ACj\cdot\biso,\ACg\cdot\blin))=\yabi\cdot\int_{R_n}\left|\mathrm{Res}(P_\ACj,P_{x\ACg})\right|_F^{-1}\dd x,
$$
	where $P_{x\ACg}$ is the invariant polynomial of $x\ACg\in\GL_{2h}(F)$ with respect to the embedding $\GL_{h}(K)\subset\GL_{2h}(F)$ induced by $\blin$. Here $\dd x$ is the normalized Haar measure of $R_n$.
\end{thm}
This theorem generalizes the Theorem \ref{babys} by allowing $\ACg\notin\GL_{2h}(\OO F)$ as well as admitting level structures. By allowing $g\in\GL_{2h}(F)$ one can construct more \Heegner cycles. For example, when $n=0$, for any open $\OO F$-subalgebra $\mathcal O \subset \OO K$, by an appropriate choice of $(\ACj,\ACg)\in D_F^\times\times\GL_{2h}(F)$, one can construct the cycle $Z_0(\ACj\cdot\biso,\ACg\cdot\blin)$ corresponding to the deformation of $\GF$ as a formal $\mathcal O$-module such that $\mathcal O$ acts on its Tate module through a preselected embedding $\mathcal O\rightarrow\matt{2h}{2h}{\OO F}$. Such a deformation is also known as a quasi-canonical lifting of the formal $\OO F$-module $\GF$ in the literature (see \cite{gross1986canonical}).


Now we generalize our formula for correspondences. 
For any non-negative integer $n$, let $$\mathscr H(R_n\backslash\GL_{2h}(F)/R_n):=\{f\in C_c^\infty(\GL_{2h}(F)):f(xg)=f(gx)=f(g) \text{ for any }x\in R_n\}$$ be the set of double $R_n$-invariant compactly supported real-valued functions on $\GL_{2h}(F)$. 
Suppose $g\in\GL_{2h}(F)\cap\matt{2h}{2h}{\OO F}$. Let $m$ be an integer such that $g(\OO F^{2h})\supset \pi^m\OO F^{2h}$, and let $\ACj\in D_F^\times$ be an element such that $|\nrd(\ACj)|_F^{-1}=|\det(\ACg)|_F$. The pair $(\ACj,\ACg)$ induces a natural morphism
$$
\map{\beta_n^{m+n}(\ACj,\ACg)}{\DEF{m+n}}{\DEF n}
$$
for any $n\geq 0$ (see \eqref{haww}).
Denote by $\beta^{m+n}_n$ the transition map from $\DEF{m+n}$ to $\DEF n$. Let
\[equation]{\lb{jiaolvjiaolv}
	f_{R_n\cdot g\cdot R_n,m}(x):=\frac{\deg(\xymatrix{\DEF{m+n}\ar[rr]^{\beta^{m+n}_n}&&\DEF n})}{\Vol(R_n\cdot g\cdot R_n)}\one_{R_n\cdot g\cdot R_n}(x)\in \mathscr H(R_n\backslash\GL_{2h}(F)/R_n),
	}
where $\one_{R_n\cdot g\cdot R_n}(x)$ is the characteristic function of $R_n\cdot g\cdot R_n$. Then the correspondence associated to the pair $(\ACj,f)$ is defined by
$$
\xymatrix{\DEF n&&\DEF{m+n}\ar[rr]^{\beta^{m+n}_n(\ACj,\ACg)}\ar[ll]_{\beta^{m+n}_n}&&\DEF n}.
$$

To introduce a general formula, we fix the following constants
$$
C_0(K_1,K_2):=\frac{\left(\prod_{i=1}^h\zeta_{K_1}(i)\right)\left(\prod_{i=1}^h\zeta_{K_2}(i)\right)}{\prod_{i=1}^{2h}\zeta_F(i)};\qquad C_n(K_1,K_2):=1
$$
for any $n>0$, and any two quadratic extensions $K_1,K_2$ over $F$.
\begin{thm}\lb{SUPER}
	Let $f=f_{R_n\cdot g_0\cdot R_n,m}$ and $C_n=C_n(K,K)$. Suppose $|\nrd(\ACj)|_F^{-1}=|\det(\ACg)|_F$ and $|\nrd(\ACj_0)|_F^{-1}=|\det(\ACg_0)|_F$. If $P_{\ACj_0^{-1}\ACj}$ is irreducible and $P_{\ACj_0^{-1}\ACj}(0)P_{\ACj_0^{-1}\ACj}(1)\neq 0$, or if $$\Res\left(P_{\ACj_0^{-1}\ACj},P_{x^{-1}g}\right)\neq 0$$ for any $x$ such that $f(x)\neq 0$, then the following intersection number is finite and given by the formula 
	$${\[split]{&
	\chi(Z_n(\biso,\blin)\litimes_{\DEF n} \beta^{m+n}_{n*}\beta^{m+n}_n(\ACj_0,\ACg_0)^* Z_n(\ACj\cdot\biso,\ACg\cdot\blin))\\=\;&C_n\yabi\int_{\GL_{2h}(F)}f(x)\left|\Res\left(P_{\ACj_0^{-1}\ACj},P_{x^{-1}g}\right)\right|_F^{-1}\dd x,
	}}$$
	where $P_{x^{-1}\ACg}$ is the invariant polynomial of $x^{-1}\ACg\in\GL_{2h}(F)$ with respect to the embedding $\GL_{h}(K)\subset\GL_{2h}(F)$ induced by $\blin$. Here the Haar measure $\dd x$ is normalized by $\GL_{2h}(\OO F)$.
\end{thm}
Note that when $\ACj_0=\id$, $\ACg_0=\id$ and $m=0$, this formula specializes to Theorem \ref{babys} and Theorem \ref{girls}. The formula can also be generalized to the case of $K_1\not\cong K_2$. For $i=1,2$, we put extra subscripts and denote by $Z_n(\biso_i,\blin_i)$ the cycle constructed by the data $\blin_i: K_i^h\rightarrow F^{2h}$ and $\biso_i:\GG_{K_i}\rightarrow \GF$, and we associate a matrix $\Delta(\biso_i,\blin_i)\in\matt{2h}{2h}{D_F}$ to $(\biso_i,\blin_i)$ (See Definitions \ref{quaiso} and \ref{matu}). Again they are all equi-height pairs. 
The following theorem gives a formula in the most general case without using invariant polynomials.
\begin{thm}\lb{BIQ}
	Let $f=f_{R_n\cdot g_0\cdot R_n,m}$ and $C_n=C_n(K_1,K_2)$. Suppose $|\nrd(\ACj)|_F^{-1}=|\det(\ACg)|_F$ and $|\nrd(\ACj_0)|_F^{-1}=|\det(\ACg_0)|_F$. If the integrand is bounded, then the following intersection number is finite and given by the formula 
	$${\[split]{
		&\chi(Z_n(\biso_1,\blin_1)\litimes_{\DEF n} \beta^{m+n}_{n*}\beta^{m+n}_n(\ACj_0,\ACg_0)^* Z_n(\ACj\cdot\biso_2,\ACg\cdot\blin_2))
	\\=\;&C_n\yabii\int_{\GL_{2h}(F)}f(x)\left|\mathbf R(x)\right|_F^{-1}\dd x,
}}$$
	where $\mathbf R(x)$ is the reduced norm of a matrix $R(x)\in\matt hh{D_F}$, and
$$
	R(x):=\[pmatrix]{0_h&I_h}\Delta(\biso_1,\blin_1)^{-1}\cdot \ACj_0^{-1}\ACj\cdot x^{-1}\ACg \cdot \Delta(\biso_2,\blin_2)\[pmatrix]{I_h\\0_h}.
$$
Here the Haar measure $\dd x$ is normalized by $\GL_{2h}(\OO F)$.
\end{thm}
\subsection{The linear AFL}\lb{ss lAFL}
When $K/F$ is an unramified extension with odd residue characteristic, Zhang has a conjecture called the linear AFL, which gives another formula for the intersection numbers of \Heegner cycles. Recall $Z_0:=(\eta_0^0)_*\OO{{\DEf0}}$ as defined in Section \ref{bababa}. Let $\ACj\in D_F^\times$ and $g\in \GL_{2h}(F)$ be elements such that $|\nrd(\ACj)|_F^{-1}=|\det(g)|_F$. Let $f=f_{R_0gR_0,m}$ as defined in \eqref{jiaolvjiaolv}. The linear AFL predicts that the intersection number
$$
\chi(Z_0\litimes_{\DEF0}\beta^m_{0*}\beta^{m}_0(\ACj,\ACg)^* Z_0)
$$
is given by the central derivative of a certain orbital integral $\OOF(g(\ACj),f,s)$, which we describe now.

Let $G'=\GL_{2h}(F)$ and $H^\prime=\GL_h(F)\times\GL_h(F)$. Consider $H^\prime$ as a subgroup of $G^\prime$ by the blockwise diagonal embedding.
Let $g(\ACj)\in G^\prime$ be an element having the same invariant polynomial $P_\ACj$ with $\ACj$. Let $\eta$ be the non-trivial quadratic character associated to $K/F$. We regard $\eta$ and $\left|\bullet\right|_F$ as characters on $H^\prime$ by precomposing them with $(g_1,g_2)\mapsto \det(g_1^{-1}g_2)$ (note the inverse on $g_1$). Consider the following orbital integral
\begin{equation}\label{afl}\comments{afl}
\OOF\left(g(\ACj),f,s\right)=\int_{\frac{H^\prime\times H^\prime}{I(g(\ACj))}}f\left(h_1^{-1}g(\ACj)h_2\right)\eta(h_2)|h_1h_2|_F^s\dd h_1\dd h_2,
\end{equation}
where for any $g\in G'$,
$$
I\left(g\right)=\{(h_1,h_2)|h_1g=gh_2\}.
$$
\begin{conj}[linear AFL]\label{po}\comments{po}
	Let $\ACj,g,f$ and $g(\ACj)$ be as defined above. If $P_\ACj$ is irreducible and $P_\ACj(0)P_\ACj(1)\neq 0$, we have
\begin{equation}
\pm(2\ln q)^{-1}\left.\frac{\dd}{\dd s}\right|_{s=0}\OOF(g(\ACj),f,s)=
	\chi(Z_0\litimes_{\DEF0} \beta^{m}_{0*}\beta^{m}_0(\ACj,\ACg)^* Z_0),
\end{equation}
where the sign is chosen so that the left hand side is non-negative.
\end{conj}
In Section \ref{EVI}, we verify this conjecture for $h=1$. By calculating both sides of this identity, the author also proves the linear AFL in the case of $h=2$ for the characteristic function of $\GL_4(\OO F)$ in \cite{li2017arithmetic}. 

\subsection{Outline of contents}The main idea is to raise the problem to the infinite level. We review some history. In Theorem 6.4.1 of the paper \cite{scholze2012moduli} of Scholze-Weinstein, and also in the paper \cite{weinstein2013semistable} of Weinstein, they show that the projective limit of the generic fiber of the Lubin--Tate tower for $\GF$ is a perfectoid space $M_{\infty}$, which can be embedded into the universal cover of $\GG^{2h}$, where $\GG$ is a certain deformation of $\GF$. The main idea is to describe the infinite level Lubin--Tate spaces by $\GG^{2h}$.

In this article, after establishing definitions of Lubin--Tate towers and \Heegner cycles in \S \ref{chenmi}, we realize their idea on the integral model of the Lubin--Tate tower on finite levels by a new construction. In \S\ref{mo}, we prove that the preimage of the closed point under the transition map $\DEF n\rightarrow \DEF 1$ is canonically isomorphic to $\GF^{2h}[\pi^{n-1}]$. In other words, the following diagram is Cartesian (See Proposition \ref{haolemma})
$$
\cate{\GF^{2h}[\pi^{n-1}]}{\Spec\CFF q}{\DEF n}{\DEF 1.}{}{}{}{}
$$
These heuristic examples let us regard $\GF^{2h}$ as an approximation of the special fiber of $\DEF n$ over $\DEF1$ when $n\rightarrow\infty$. Therefore, it is natural to construct \Heegner cycles $\ca\biso\blin\infty$ on $\GF^{2h}$ and formulate a similar intersection problem. 
In \S\ref{mi}, by Proposition \ref{haolemma}, we show that the corresponding intersection number on $\GF^{2h}$ is the same as the intersection number on the Lubin--Tate tower if the level is sufficiently high.
We calculated in \S\ref{zhaojile} Proposition \ref{xiaomubiao} that the intersection number of $\ca\biso\blin\infty$ and $\ca{\ACj\biso}{\ACg\blin}\infty$ is related to $\|\Res(P_\ACj,P_\ACg)\|_F^{-1}$. 
In \S \ref{FORMULAA}, we prove our main Theorem \ref{SUPER} by using the projection formula. The essential property for the method in \S \ref{FORMULAA} to work is that the transition maps of the Lubin--Tate tower are generically \'etale. In \S\ref{EVI}, we prove the linear AFL for the case of $h=1$ as an application.


\subsection{Acknowledgement}
The author would like to thank Professor Wei Zhang, for suggesting this topic and for his comments on this article. The author also thanks Professor Johan de Jong for answering his several questions. The author would like to express his sincere gratitude to the three anonymous referees for their reports which help the author make significant improvements on both the mathematical contents and the exposition. Last but not least, the author thanks Professor Liang Xiao for his encouragement, Ziquan Yang and Stanley Yao Xiao for their suggestions on writing.


\section{\Heegner cycles on Lubin--Tate towers}\label{chenmi}\comments{chenmi}

We fix a non-Archimedean local field $F$ with ring of integers $\OO F$ and residue field $\FF_q$. Let $K$ be a finite separable extension of $F$ of degree $k$ (of any characteristic), with standard discrete valuation $\vv_K$, ring of integers $\OO K$, and residue field $\FF_{q_K}$. We fix a uniformizer $\pi$ of $F$, which is not necessarily a uniformizer of $K$. Let $\breve K$ (resp. $\breve F$) be the completion of the maximal unramified extension of $K$ (resp. $F$) with ring of integers $\CO K$ (resp. $\CO F$). Denote by $\CCC$ (resp. $\CCC_K$) the category of complete Noetherian local $\CO F$-(resp. $\CO K$-)-algebras with residue field $\CFF q$. For $R=K,F,\OO K$ and $\OO F$, denote by $R^n$ the $R$-module of $n\times 1$ matrices over $R$, and by $R^{n\vee}$ the $R$-module of $1\times n$ matrices over $R$. The pairing 
\[equation]{\lb{papa}\langle-,-\rangle:R^{n\vee}\otimes_R R^{n}\rightarrow R}
is given by the matrix multiplication.

We first construct the Lubin--Tate towers, which admit an action by $D_F^\times\times\mathrm{GL}_{kh}(F)$ where $D_F$ is the $F$-central division algebra of invariant $1/kh$. We then construct \Heegner cycles corresponding to algebraic group embeddings $D_K^\times\subset D_F^\times$ and $\mathrm{GL}_{h}(K)\subset \mathrm{GL}_{kh}(F)$. Here $D_K$ is the $K$-central division algebra of invariant $\frac1h$. Denote by $\CFF q$ the algebraic closure of $\FF_q$.

\subsection{The Lubin--Tate tower}In this subsection, we give a precise definition of the Lubin--Tate tower associated to a formal $\OO K$-module $\GK$ of height $h$ over $\CFF{{q_K}}$. 
\subsubsection{Formal modules}
Let $A$ be a complete Noetherian local algebra over a local ring $B$ with structure map $s:B\longrightarrow A$ and residue field $\CFF q$. A formal $B$-module over $A$ is a one dimensional formal group law $\mathcal G$ endowed with a $B$-action $i:B\longrightarrow \End(\GG)$, where we insist that the induced action on the Lie algebra 
$$B\longrightarrow \End(\mathrm{Lie \GG})\cong A$$
is identical to the structure map $s$. A homomorphism of two formal $B$-modules $\beta:\mathcal \GG_1\longrightarrow \mathcal \GG_2$ over $A$ is a homomorphism of formal group laws that commutes with actions of $B$. When $\beta$ is surjective and $pA=0$ , the height of $\beta$ is defined by
\[equation]{\lb{heso}
\mathrm{Height}(\beta):=\log_{q_B}(\deg \beta)
}
where $q_B$ is the cardinality of the residue field of $B$. Furthermore, if $B$ is a discrete valuation ring with uniformizer $\pi_B$, then we define the height of $\GG$ as the height of the endomorphism induced by $\pi_B\in B\subset \mathrm{End}(\GG)$. Furthermore, we denote by $[b]_\GG:\GG\longrightarrow \GG$ the action of $b\in B$, $[+]_\GG:\GG\times \GG\longrightarrow \GG$ (resp.$[-]_\GG$) the addition (resp. subtraction) of $\GG$, and $\overline\GG$ the base change of $\GG$ to the residue field of $A$. For any homomorphism of formal modules $\map\psi{\GG_1}{\GG_2}$, we denote its reduction by $\map{\overline\psi}{\overline{\GG_1}}{\overline{\GG_2}}$.


\[defi]{
A deformation of $\GK$ over $A\in \CCC_K$ is a pair $(\GG,\iota)$ such that
\[itemize]{
\item $\GG$ is a formal $\OO K$-module over $A$,
\item $\map{\iota}{\GK}{\overline{\GG}}$ is an $\OO K$-quasi-isogeny.
	}
}
We define level structures following Drinfeld \cite{drinfel1974elliptic}.
\[defi]{\label{levelstru}\comments{levelstru}
For any $x\in\OO K$, an $x$-level structure of a formal $\OO K$-module $\GG$ over $A\in\CCC$ is a homomorphism of left $\OO K$-modules
$\map{\alpha}{\OO K^{h\vee}}{\GG(A)}$ such that 
\[itemize]{
\item $\alpha$ factors through $\OO K^{h\vee}\longrightarrow\OO K^{h\vee}/x\OO K^{h\vee}\longrightarrow \GG(A)$,
\item $\alpha$ satisfies the non-degeneracy condition that the power series $[x]_\GG(X)$ is divisible by $$\prod_{\vv\in\OO K^{h\vee}/x\OO K^{h\vee}}\left(X[-]_\GG\alpha(\vv)\right).$$
	}
}
\[defi]{For any $x\in\OO K$, an $x$-level deformation is a triple $(\GG,\iota,\alpha)$ obtained by attaching an $x$-level structure $\alpha$ to a deformation $(\GG,\iota)$ of $\GK$ over $A\in \CCC_K$.}

\[defi]{\lb{ahi}Two $x$-level deformations $(\GG_1,\iota_1,\alpha_1)$ and $(\GG_2,\iota_2,\alpha_2)$ are equivalent if there is an isomorphism $\map{\zeta}{\GG_1}{\GG_2}$ of formal $\OO K$-modules such that the following two diagrams commute
$$
\xymatrix{
	\GG_K\ar[rrr]^{\mathrm{id}}\ar[d]_{\iota_1}&&&\GG_K\ar[d]^{\iota_2}\\
	\overline{\GG_1}\ar[rrr]^{\overline\zeta}&&&\overline{\GG_2}
	}
\qquad
\xymatrix{
	\mathcal O_K^{h\vee}\ar[rrr]^{\mathrm{id}}\ar[d]_{\alpha_1}&&&\mathcal O_K^{h\vee}\ar[d]_{\alpha_2}\\
	\GG_1\ar[rrr]^{\zeta}&&&\GG_2.
	}
$$
We denote by $[\GG,\iota,\alpha]_n$ the equivalence class of $(\GG,\iota,\alpha)$, where the subscript $n$ is indicating that $\alpha$ is an $x$-level structure with $\vv_K(x)=n\cdot \vv_K(\pi)$.
}
Note that for any $x,x'\in\OO K$ with $\vv_K(x)=\vv_K(x')$, an $x$-level structure is an $x'$-level structure. Therefore the subscript $n$ is sufficient to describe the level of $\alpha$ in the triple $[\GG,\iota,\alpha]_n$. For simplicity, we will call the $x$-level structure as the $\pi^n$-level structure. The use of $\pi^n$ is symbolic and it may not be a well-defined element in $\OO K$ when $n$ is not an integer.

\[thm]{[Drinfeld, Lubin--Tate \cite{drinfel1974elliptic}\cite{lubin1966formal}] \label{LTlevel}\comments{LTlevel}The functor sending $A\in\CCC_K$ to the set of equivalence classes of $\pi^n$-level deformations of $\GK$ over $A$ is representable by a formal scheme $\DEf n^{\sim}$ of dimension $h$ over $\Spf\CO K$. This formal scheme decomposes as
$$
\DEf n^\sim = \coprod_{j\in \ZZ}\DEf n^{(j)}
$$
where $\DEf n^{(j)}$ is the open and closed formal subscheme on which $\iota$ has height $j$, and there is a regular local $\CO K$-algebra $A_n$ such that $\DEf n^{(j)}=\Spf(A_n)$ for every $j$.

}

\begin{defi}\label{LTtower}\comments{LTtower}
	The Lubin--Tate tower $\DEf\bullet^\sim$ associated to $\GK$ is a projective system consisting of $\DEf n^\sim$ for $n\in\frac1{\vv_K(\pi)}\cdot\ZZ_{\geq 0}$ with transition maps, functorial in $A\in\CCC_K$, given by
$$
\mapi{\DEf{n+m}^\sim(A)}{\DEf n^\sim(A)}{[\GG,\iota,\alpha]_{n+m}}{\left[\GG,\iota,[\pi^m]_{\GG}\circ\alpha\right]_{n}}
$$
	for any integer $m\geq0$. We denote this transition map by $\map{\theta_n^{m+n}}{\DEf{m+n}^\sim}{\DEf n^\sim}$.
\end{defi}
These transition maps are finite flat, and they keep the superscript of components --- they induce maps between $\DEf {n+1}^{(j)}$ and $\DEf{n}^{(j)}$ for every $n$ and $j$. We denote by $\DEf\bullet^{(j)}$ the subtower of $\DEf\bullet^\sim$ with superscript $(j)$.
\begin{rem}
The subscript $n$ of $\DEf n^{(j)}$ refers to the $\pi^n$-level structure. Since $\pi$ is not necessarily a uniformizer of $\OO K$, a fractional subscript like $\DEf {\frac{n}k}^{(j)}$ can make sense if $K/F$ is ramified. But we do not need these spaces with fractional-subscripts in our article. From now on, all subscripts $n$ will be a non-negative integer.
\end{rem}

\subsection{Maps between Lubin--Tate towers}

Let $\mathcal G_F$ be the formal $\mathcal O_F$-module of height $kh$ obtained by only remembering the $\mathcal O_F$-action of $\mathcal G_K$, and denote its associated Lubin--Tate tower by $\DEF\bullet^\sim$, which is a regular formal scheme of dimension $kh$ over $\Spf\CO F$. Denote the degree of the corresponding residue field extension of $K/F$ by
$$
r:=[\FF_{q_K}:\FF_q]=\log_qq_K.
$$
From now on, we work over the base $\Spf\CO F$ and we regard each formal scheme as a functor from $\CCC$ to sets. We can extend the domain of the functor $\DEf n^\sim$ from $\CCC_K$ to $\CCC$ as follows. For any $A\in\CCC$, we define $\DEf n^\sim(A)$ to be the set of $\pi^n$-level deformations $[\GG,\iota,\alpha]_n$ over $A$ with extra endomorphisms by $\OO K$. The action of $\OO K$ on $\mathrm{Lie}(\GG)\cong A$ gives $A$ an $\OO K$-algebra structure, therefore $\DEf n^\sim(A)\neq\emptyset$ implies that $A\in \CCC_K$.

A map from $\DEf\bullet^\sim$ to $\DEF\bullet^\sim$ over $\Spf\CO F$ can be constructed by choosing a pair of homomorphisms $(\biso,\blin)$
\begin{equation}\label{pairsa}\comments{pairsa}
\begin{split}
\map{\biso}{\GK}{\GF};\\
\map{\blin}{K^{h}}{F^{kh}},
\end{split}
\end{equation}
where $\blin$ is an $F$-linear isomorphism and $\biso$ is a quasi-isogeny of formal $\OO F$-modules. In fact, this $\blin$ is purely symbolic. 
 
Our later definitions and proofs will work instead with the dual map $$\map{\blinn}{F^{kh\vee}}{K^{h\vee}},$$ which is defined by $l\mapsto \blinn(l)$ where $\blinn(l)\in K^{h\vee}$ is the unique element such that $\langle l,\blin(v)\rangle=\tr_{K/F}\langle \blinn(l),v\rangle$ for any $v\in K^h$. Here the paring $\langle-,-\rangle$ is defined in \eqref{papa}.
The uniqueness of $\blinn(l)$ requires the non-degeneracy of the form $\tr_{K/F}\langle -, -  \rangle$, which is valid because the extension $K/F$ is separable.

Roughly speaking, the map $\DEf\bullet^\sim\rightarrow \DEF\bullet^\sim$ is obtained by identifying special fibers of deformations of $\GK$ with the ones of $\GF$ through $\biso$, and identifying their level structures through $\blin^\vee$. 
We prefer working with $\blinn$ since the action of $\GL_h(\OO K)$ on the $\pi^n$-level structures is a right action, while the action of $\GL_h(\OO K)$ on $\blin$ is a left action.

Note that the data $\tau$ and $\biso$ induce algebraic group embeddings $\mathrm{GL}_h(K)\subset \mathrm{GL}_{hk}(F)$ and $D_K^\times\subset D_F^\times$, which will be the essential data parametrizing our \Heegner cycles. Therefore, for our convenience, among those $\tau$ which induce the same group embedding, we can select one such that
\[equation]{\label{import}\comments{import}
\blinn(\OO F^{kh\vee})\subset \OO K^{h\vee}.
}

Denote the set of morphisms $\map fXY$ of two formal schemes over $\Spf\CO F$ by $\mathrm{Mor}_{\CO F}(X,Y)$.

\begin{rem}\label{towermor}\comments{towermor}
	In our article, a map for two towers $$\map{\pppc}{\DEf\bullet^\sim}{\DEF\bullet^\sim}$$ means an element in 
$$\varprojlim_j\varinjlim_i\mor(\DEf i^\sim,\DEF j^\sim).$$
	This kind of element is uniquely determined by a system of compatible morphisms $\DEf{m+n}^\sim\rightarrow \DEF n^\sim$ over $\Spf\CO F$ for all $n\geq 0$ with some fixed $m\geq 0$. In such a system we denote each member by
	\[equation]{\label{hawa}\comments{hawa}
	\map{\pccc{n}{m}{k}{j}}{\DEf{m+n}^{(j)}}{\DEF n^{(k)}}
	}
	where the superscript is intentionally chosen to specify the domain, and the subscript is chosen for the codomain. 
\end{rem}


Next we define maps of Lubin--Tate towers induced by the pair $(\biso,\blin)$. We start with the simplest case that $\tau^\vee(\OO K^{h\vee})=\OO F^{kh\vee}$, in which the corresponding algebraic group embedding $\mathrm{GL}_h(K)\subset \mathrm{GL}_{hk}(F)$ restricts to $\mathrm{GL}_h(\mathcal O_K)\subset \mathrm{GL}_{hk}(\mathcal O_F)$. Moreover, if $\alpha$ is a $\pi^m$-level structure for $\GG$, then so is $\alpha\circ\blin^\vee$. 
\begin{defi}\label{bisoindu}\comments{bisoindu}
	If $\tau^\vee(\OO K^{h\vee})=\OO F^{kh\vee}$, we define the map  $\map{\pppc}{\DEf\bullet^\sim}{\DEF\bullet^\sim}$ by  the following compatible morphisms, functorial in $A\in\CCC$,  
\begin{equation}\label{generalk}\comments{generalk}
	\maps{\eta^n_n(\biso,\blin)}{\DEf n^{(j)}(A)}{\DEF n^{(rj+\height(\biso^{-1}))}(A)}{[\GG,\iota,\alpha]_n}{[\GG,\iota\circ\biso^{-1},\alpha\circ \blin^\vee]_n}
\end{equation}
for each $n\geq 0$.
In other words, $\eta^n_n(\biso,\blin)$ maps $[\GG,\iota,\alpha]_n$ to $[\GG,\iota',\alpha']_n$ such that the following two diagrams
	\[equation]{\lb{dagaidagai}
	\cate {\GK}{\GF}{\overline\GG}{\overline\GG}{\biso}{\iota}{\iota'}{\id_{\overline{\GG}}}\qquad 
\xymatrix{
	\OO F^{kh\vee}\ar[rr]^{{\blin^\vee}}\ar[d]_{\alpha'}&&\OO K^{h\vee}\ar[d]^{\alpha}\\
	{\GG}&&{\GG}\ar[ll]_{\id_{\GG}}\\
	}
	}
commute.
\end{defi}
This definition is obviously well-defined as this map carries equivalent triples to equivalent triples.
In general, if $\tau^\vee(\OO K^{h\vee})\neq\OO F^{kh\vee}$, the map $\alpha\circ\blin^\vee$ may not be a well-defined level structure. We fix this issue and generalize our definition in the remaining part of this subsection.  
\begin{defi}\label{defi:condu}\comments{defi:condu}
For any $F$-linear map $\map{\blin}{K^{h}}{F^{kh}}$, the conductor $\m(\blin)$ of $\blin$ is the minimal integer $m$ such that
$$
\blin^\vee\left(\OO F^{kh\vee}  \right) \supset\pi^m\OO K^{h\vee}.
$$
\end{defi}
The idea to generalize Definition \ref{bisoindu} is allowing the identity maps in the diagram \eqref{dagaidagai} to be replaced by certain isogenies, which we describe now. For any pairs $(\biso,\blin)$, any integer $m\geq\m(\blin)$, and any equivalence class of $\pi^{m+n}$-level deformation $[\GG,\iota,\alpha]_{m+n}$ over $A$, we construct a formal $\OO K$-module $\GG'$ and an isogeny $\map\psi{\GG}{\GG'}$ over $A$ as follows. Note that 
\begin{equation}\label{largeen}\comments{largeen}\blin^\vee(\OO F^{kh\vee})\supset \pi^m\OO K^{h\vee}.\end{equation}
Put
\begin{equation}\label{lattice}\comments{lattice}V=\blinn(\pi^n\OO F^{kh\vee})/\pi^{n+m}\OO K^{h\vee}.\end{equation}
	Since $\map\alpha{\OO F^{kh\vee}}{\GG(A)}$ is a Drinfeld $\pi^{n+m}$-level structure, by definition, it induces a morphism
		$$\map{\widetilde\alpha}{\OO F^{kh\vee}/\pi^{n+m}\OO F^{kh\vee}}{\GG(A)}.$$
Consider a power series $\psi$ defined by 
\begin{equation}\label{baga}\comments{baga}
\psi(X)=\prod_{\vv\in V}(X[-]_\GG\widetilde{\alpha}(\vv)).\end{equation}
By Serre's construction, there exists a formal $\OO F$-module $\GG'$ over $A$ such that $\map\psi\GG{\GG'}$ is an isogeny. Note that this construction is functorial in $A\in\CCC$. 
\[defi]{We call the isogeny $\map\psi\GG{\GG'}$ in \eqref{baga} Serre's isogeny associated to $m$,$(\biso,\blin)$ and $\GG$.}

\begin{lem}\label{serre}\comments{serre}
	Let $\map{\psi}{\GG}{\GG'}$ be Serre's isogeny associated to $m,(\biso,\blin)$ and $\GG$ for some $m\geq\m(\blin)$.
Suppose $\alpha$ is a $\pi^{m+n}$-level structure for $\GG$, then $\psi\circ\alpha\circ\blinn$ is a $\pi^n$-level structure for $\GG'$.
\end{lem} 
\begin{proof} 
It suffices to show that $[\pi^n]_{\GG'}(X)$ is divisible by
$$\prod_{\vv\in\OO F^{kh\vee}/\pi^n\OO F^{kh\vee}}\Big(X[-]_{\GG'}\psi\circ\alpha\circ\blinn(\vv)\Big).$$
Since the construction of Serre's isogeny is functorial in $A\in \CCC$, without loss of generality, we assume that $\GG$ is the universal formal $\mathcal O_K$-module over $A$ where $A$ is chosen so that $\DEf{m+n}^{(0)}=\Spf A$. As $A$ is a regular local ring and, in particular, a unique factorization domain, 
it suffices to prove that $\psi\circ\alpha\circ\blinn(\vv)$ are distinct solutions of $[\pi^n]_{\GG'}(X)=0$ for $\vv\in\OO F^{kh\vee}/\pi^n\OO F^{kh\vee}$. First, for any $\vv,\ww\in\OO F^{kh\vee}/\pi^n\OO F^{kh\vee}$ such that $\vv\neq\ww$, we have $\psi\circ\alpha\circ\blinn(\vv)\neq\psi\circ\alpha\circ\blinn(\ww)$ because 
$$\blinn(\vv)-\blinn(\ww)\notin \ker\big(\psi\circ\alpha\big)=\blinn(\pi^n\OO F^{kh\vee}).$$ 
Second, $X = \psi\circ\alpha\circ\blinn(\vv)$ is a solution for  $[\pi^n]_{\GG'}(X)=0$ for any $\vv\in\OO F^{kh\vee}/\pi^n\OO F^{kh\vee}$. Indeed, 
$$
{\begin{split}
[\pi^n]_{\GG'}(\psi\circ\alpha\circ\blinn(\vv))=\;&\psi([\pi^n]_{\GG}\circ\alpha\circ\blinn(\vv))\\
=\;&\psi\circ\alpha\circ\blinn(\pi^n\vv)\\
=\;&0
\end{split}
}
$$ as desired. This completes the proof.\end{proof}

\begin{defi}\label{blinindu}\label{tower-define}\comments{tower-define}\comments{blinindu}
For any pair $(\biso,\blin)$ as described in \eqref{pairsa}, and any interger $m\geq\m(\blin)$, we define the map  $\map{\pppc}{\DEf\bullet^\sim}{\DEF\bullet^\sim}$ by the following compatible morphisms over $\Spf\CO F$ , functorial in $A\in\CCC$, 
	\begin{equation}\label{DEFI:CANOA}\comments{DEFI:CANOA}\maps{\bppc nm}{\DEf{n+m}^{(j)}(A)}{\DEF n^{\left(rj+\height(\overline{\psi}\circ\pi^{-m})+\height(\bisoo)\right)}(A)}{[\GG,\iota,\alpha]_{n+m}}{[\GG',\overline{\psi}\circ\iota\circ\pi^{-m}\biso^{-1},\psi\circ\alpha\circ\blinn]_n}\end{equation}
for any $n\geq 0$, where $\map\psi{\GG}{\GG'}$ is Serre's isogeny associated to $(\biso,\blin)$, $m$, and $\GG$. In other words, $\bppc nm$ maps $[\GG,\iota,\alpha]_{n+m}$ to $[\GG',\iota',\alpha']_{n}$ such that the following diagrams commute
\[equation]{\label{comaa}\comments{comaa}
\cate {\GK}{\GF}{\overline\GG}{\overline\GG',}{\pi^m\biso}{\iota}{\iota'}{\overline\psi}\qquad 
\xymatrix{
	\OO F^{kh\vee}\ar[rr]^{{\blin^\vee}}\ar[d]_{\alpha'}&&\OO K^{h\vee}\ar[d]^{\alpha}\\
	{\GG'}&&{\GG.}\ar[ll]_{\psi}\\
	}
}
\end{defi}
This definition is well-defined by previous lemmas and propositions. It is clear that the induced map  $\pppc$ does not depend on the choice of $m$ since one sees that $\eta^{n+m+1}_n(\biso,\blin)=\eta^{n+m}_n(\biso,\blin)\circ\theta^{n+m+1}_{n+m}$ for any $m\geq\m(\blin)$, where we recall that $\map{\theta^{n+m+1}_{n+m}}{\DEf{n+m+1}}{\DEf{n+m}}$ is the transition map. Moreover, one can show that those morphisms are finite by Proposition \ref{haolemma}.

From now on, by saying Serre's isogeny attached to $\bppc nm$ and $\GG$, we mean Serre's isogeny associated to $m$, $(\biso,\blin)$ and $\GG$. When no confusion can arise, we will simply say Serre's isogeny attached to $\bppc nm$.


\begin{cor}\lb{rmindu}
	If $\bppc nm[\GG,\iota,\alpha]_{n+m} = [\GG',\iota',\alpha']_{n}$, then there exists an isogeny $\map{\psi}{\GG}{\GG'}$ such that the diagram \eqref{comaa} commutes.
\end{cor}
\begin{proof}
Let $\GG''$ be the formal module obtained from Serre's isogeny $\map{\psi'}{\GG}{\GG''}$ attached to $\bppc nm$ and $\GG$, and write $$\bppc nm[\GG,\iota,\alpha]_{n+m} = [\GG'',\iota'',\alpha'']_{n}.$$
	Note that $\GG''$ may not equal $\GG'$. Instead, there is an isomorphism $\map{\zeta}{\GG''}{\GG'}$ translating $[\GG'',\iota'',\alpha'']_{n}$ to $[\GG',\iota',\alpha']_{n}$. Letting $\psi=\zeta\circ\psi'$, we have the following diagrams 
$$
\xymatrix{
\GK\ar[r]^{\pi^m\biso}\ar[d]_{\iota}&\GF\ar[r]^{\id_{\GF}}\ar[d]^{\iota''}&\GF\ar[d]^{\iota'}\\
	\overline\GG\ar@/_1pc/@{-->}[rr]_{\overline\psi}\ar[r]^{\overline{\psi'}}&\overline\GG''\ar[r]^{\overline\zeta}&\overline\GG'}
\qquad 
\xymatrix{
\OO K^{h\vee}\ar[d]^{\alpha}
	&\OO F^{kh\vee}\ar[l]_{\blinn}\ar[d]^{\alpha''}
	&\OO F^{kh\vee}\ar[l]_{\id}\ar[d]_{\alpha'}
	\\
	\GG\ar[r]^{\psi'}\ar@/_1pc/@{-->}[rr]_{\psi}
	&\ar[r]^{{\zeta}}\GG''
	&\GG'&\\
	&&&&
}
$$
where for each diagram, the left square commutes because of Definition \ref{blinindu}, and the right square commutes because of Definition \ref{ahi}. The corollary follows from the commutativity of outer squares.
\end{proof}

\subsection{Heights}
Note that the superscript of the codomain of the map in \eqref{DEFI:CANOA} is shifted by 
$$\mathrm{Height}(\biso^{-1}) + \mathrm{Height}(\overline\psi\circ\pi^{-m}),$$ 
where the shifting contributed by $\blin$ is $\height(\overline{\psi}\circ\pi^{-m})$. For ease of exposition, we call this number the height of $\blin$. Note that
$$
\height(\pi^{-m}\overline{\psi})=\height(\pi^{-m})+\log_q \#V=\log_q \Vol(\tau^\vee(\OO F^{kh\vee})).
$$
\begin{defi}\lb{heto}
Recall that $q$ is the cardinality of the residue field of $\OO F$. For an $F$-linear map $\map{\blin}{K^{h}}{F^{kh}}$, define the height of $\blin$ by
	\begin{equation}\label{bagong}\comments{bagong}\mathrm{Height}(\blin)=\log_{q} \Vol(\tau^\vee(\OO F^{kh\vee})),\end{equation}
		where the volume function $\Vol(\bullet)$ is normalized so that $\Vol(\OO K^{h\vee})=1$.
\end{defi}

\subsection{Group actions on Lubin--Tate towers}From now on, we fix an isomorphism $\End(\GF)\otimes_{\OO F}F\cong D_F$.
The Lubin--Tate tower $\DEF\bullet^\sim$ admits a natural action of the group $D_F^\times\times\GL_{kh}(F)$, which can be described as follows. An element $(\ACj,\ACg)\in D_F^\times\times\GL_{kh}(F)$ induces the following homomorphisms
\begin{equation}\label{pairsb}\comments{pairsb}
\begin{split}
\map{\ACj}{\GF}{\GF};\\
\map{\ACg}{F^{kh}}{F^{kh}},
\end{split}
\end{equation}
where $\ACg$ is an $F$-linear isomorphism and $\ACj$ is a quasi-isogeny of formal $\OO F$-modules. We may assume further that
\[equation]{\lb{referee}
\ACg^\vee(\OO F^{kh\vee})\subset\OO F^{kh\vee}.
}
Therefore, the pair $(\ACj,\ACg)$ induces a map functorial in $A\in\CCC$
\begin{equation}\lb{haww}
	\maps{\ppbb{n}{m+n}}{\DEF{n+m}^{(j)}(A)}{\DEF n^{\left(l\right)}(A)}{[\GG,\iota,\alpha]_{n+m}}{[\GG',\overline{\psi}\circ\iota\circ\pi^{-m}\ACj^{-1},\psi\circ\alpha\circ\ACg^\vee]_n}
\end{equation}
in a similar way as in Definition \ref{tower-define}, where $l=j-\height(\ACj)+\height(\ACg)$, and $\map\psi\GG{\GG'}$ is Serre's isogeny attached to $\ppbb{n}{m+n}$. Here $\mathrm{Height}(g)$ is defined similarly as in \eqref{bagong}. Note that the construction of $\psi$ and $\ppbb{n}{m+n}$ are essentially the same as in Definition \ref{tower-define} if we replace $K/F$ by the trivial extension $F/F$. Therefore the same arguments prove that the map $\ppbb{n}{m+n}$ is well-defined, and hence we obtain a map $\map{\beta_{\bullet\sim}^{\bullet\sim}(\ACj,\ACg)}{\DEF\bullet^\sim}{\DEF\bullet^\sim}$ for the Lubin--Tate tower.

If both $\ACj$ and $\ACg$ are identity, we abbreviate $\bbpb nm{\id}{\id}$ to 
$$\map{{\ppb}}{\DEF{n+m}^\sim}{\DEF n^\sim},$$
which coincides with the transition map between $\DEF{n+m}^\sim$ and $\DEF n^\sim$ in the Lubin--Tate tower.

\subsection{\Heegner cycles on the Lubin--Tate tower}\lb{kgroup}
In this subsection, we introduce \Heegner cycles on the Lubin--Tate tower $\DEF\bullet^\sim$ induced by
\begin{equation}\label{inflev}\comments{inflev}
	\map{\pppc}{\DEf\bullet^\sim}{\DEF\bullet^\sim}.
\end{equation}
In this article, we construct cycles as coherent sheaves. To define a system of compatible cycles, we need to uniformize their multiplicities. Therefore we need to consider things like $\mathcal F^{\oplus\frac12}$ for some coherent sheaf $\mathcal F$  on a formal scheme $X$. This motivates us to define the following concept. Denote by $\graa X$ the free abelian group generated by isomorphism classes of coherent sheaves of $\OO X$-modules, modulo the relations by $[\mathcal F^{\oplus n}]=n[\mathcal F]$ for any coherent sheaf  $\mathcal F$ of $\OO X$-modules on $X$, where we use $[\mathcal F]$ to denote the corresponding element in $\graa X$ associated to $\mathcal F$.
\begin{defi}\label{tcycle}\comments{tcycle}
	Let $\pppc$ be the map between Lubin--Tate towers induced by $(\biso,\blin)$ as in Definition \ref{tower-define}. The corresponding \Heegner cycle $Z_\bullet^\sim(\biso,\blin) $ is a collection of cycles $$\cycj\in\groo{\DEF n^{(j)}}$$ for each $n$ and $j$, where each cycle $\cycj$ is defined as follows: Suppose the map $\pppc$ from $\DEf{n+m}^\sim$ to $\DEF n^\sim$ is given by 
	$$\map{\pccc{n}{m}{j}{l}}{\DEf{n+m}^{\left(l\right)}}{\DEF n^{(j)}}$$ 
where $l=\frac{j}r-\height(\blin)+\height(\biso)$. If $l$ is an integer, we define 
	\[equation]{\lb{zhaogebishengtianhaodelaopo}
	\cycj:=\frac1{\deg\left(\DEf{n+m}^{\left(l\right)}\rightarrow \DEf{n}^{\left(l\right)}\right)} 
\left[\pccc nmjl_*
\OO{\DEf{n+m}^{\left(l\right)}}\right]
	}
where the map $\DEf{n+m}^{\left(l\right)}\rightarrow \DEf{n}^{\left(l\right)}$ is the transition map. Otherwise, we define $\cycj:=0$.
\end{defi}
\begin{rem}
	The Definition \ref{tcycle} does not depend on $m$ because each transition map $$\map{\theta_{m}^{m+1}}{\DEf{m+1}^{(l)}}{\DEf m^{(l)}}$$ is a finite flat map of formal spectra of regular local rings, and therefore $\theta_{m*}^{m+1}\OO{\DEf{m+1}^{(l)}}\cong \OO{\DEf m^{(l)}}^d$ for $d:=\deg\left(\DEf{m+1}^{(l)}\rightarrow\DEf m^{(l)}\right)$. Further, since each $(\gamma_k,g_k)\in D_K^\times\times \GL_h(K)$ induces an isomorphism of the tower $\DEf\bullet^\sim$, one finds that $Z_\bullet^\sim(\biso\circ \gamma_k,\blin\circ g_k) = Z_\bullet^\sim(\biso,\blin)$ by an analogue of Lemma \ref{wanquanbuzhi}. Therefore, the cycle only depends on the algebraic group embeddings $D_K^\times\rightarrow D_F^\times$ and $\GL_h(K)\rightarrow \GL_{2h}(F)$ induced by $\biso$ and $\blinn$. 
\end{rem}

\subsection{Classical Lubin--Tate spaces}Recall that the Lubin--Tate tower is a disjoint union of components indexed by integers
$$
\DEF\bullet^\sim=\coprod_{j\in\ZZ} \DEF\bullet^{(j)}.
$$
Note that the action of an element $\omega\in\OD F{kh}^\times$ with valuation $j$ maps $\DEF n^{(j)}$ isomorphically onto $\DEF n^{(0)}$, which reduces all problems to a single space $\DEF n^{(0)}$. From now on, without loss of generality, we only consider Lubin--Tate spaces with superscript $(0)$ for ease of exposition, whereas we come back to the general situation later in \S\ref{removal}.
\begin{defi}\lb{clacla}
We call a space in the Lubin--Tate tower with superscript (0) ( for example $\DEf n^{(0)}$ or $\DEF n^{(0)}$) a classical Lubin--Tate space. For simplicity, we omit their superscript and denote them by $\DEf n$ or $\DEF n$.
\end{defi}
When considering \Heegner cycles on $\DEF n$, we may also omit the superscript $(0)$ of a cycle and simply denote it by 
$$
\cycy n:=Z_n^{(0)}(\biso,\blin). \quad
$$
\subsubsection{The height restrictions for \Heegner cycles on classical Lubin--Tate spaces}
To induce a morphism from $\DEf \bullet$ to $\DEF \bullet$, a pair $(\biso,\blin)$ must be restricted so that the superscript $(0)$ remains unchanged after the map $\pppc$ (see Definition \ref{tower-define}). This is equivalent to requiring $\height(\blin)=\height(\biso)$.
\begin{defi}\label{equi}\comments{equi}
We call a pair of morphisms in \eqref{pairsa} an equi-height pair if 
$$\height(\blin)=\height(\biso).$$
Denote by $\equi KFh$ the set of equi-height pairs.
\end{defi}
\begin{rem}\label{remequi}\comments{remequi}
	Similarly, applying Definition \ref{equi} to the trivial extension $F/F$, one obtains that $\equi FF{kh}$ is the set of elements $(\ACj,\ACg)\in \OD F{kh}^\times\times\GL_{kh}(F)$ such that $\mathrm{Height}(\ACj)=\mathrm{Height}(\ACg)$, or equivalently,  
$$
	|\det(\ACg)|_F=|\nrd(\ACj)|_F^{-1},
$$
	where $\nrd(\bullet)$ is the reduced norm of $\OD F{kh}$.
\end{rem}




\subsection{Composition of maps of Lubin--Tate spaces}\label{yange}\comments{yange}  
In this subsection, we introduce several properties of compositions of maps between Lubin--Tate towers for our applications later.
\begin{prop}\label{Composition}\comments{Composition}
Letting $(\biso,\blin)\in\equi KFh$, $\ACpair\in\equi FF{kh}$, $m_1\geq \m(\blin)$ and $m_2\geq\m(\ACg)$, we have
$$
\bpbb{m_2}\circ \bppc{n+m_2}{m_1}=\bpdd{m_1+m_2}
.$$ In other words, the following diagram commutes
$$
	\trii {\DEF {n+m_2}}{\DEF n.}{\DEf{n+m_1+m_2}}{\bpbb{m_2}}{\bpdd{m_1+m_2}}{\bppc{n+m_2}{m_1}}
$$
\end{prop}

\begin{proof}
	For any $A\in\CCC$, and any $[\GG,\iota,\alpha]_{n+m_1+m_2}\in \DEf{n+m_1+m_2}(A)$, let
\begin{equation}
\begin{split}
	[\GG^\prime,\iota^\prime,\alpha^\prime]_{n+m_2}:=\;&\ppcc{n+m_2}{m_1}[\GG,\iota,\alpha]_{n+m_1+m_2};\\
	[\GG^{\prime\prime},\iota^{\prime\prime},\alpha^{\prime\prime}]_n:=\;&\bpbb{m_2}[\GG^\prime,\iota^\prime,\alpha^\prime]_{n+m_2},
\end{split}
\end{equation}
and denote Serre's isogenies attached to $\ppcc{n+m_2}{m_1}$ and $\bpbb{m_2}$ by $$\map{\psi_1}\GG{\GG'}, $$$$
	\map{\psi_2}{\GG'}{\GG''}.$$ 
	By Definition \ref{blinindu} and using \eqref{comaa}, we have the following commutative diagrams
	\[equation]{\lb{azi}
	\xymatrix{
		\GK\ar[r]^{\pi^{m_1}\biso}\ar[d]^{\iota}&\GF\ar[r]^{\pi^{m_2}\ACj}\ar[d]^{\iota'}&\GF\ar[d]^{\iota''}\\
		\overline\GG\ar[r]^{\overline{\psi_1}}&\overline\GG'\ar[r]^{\overline{\psi_2}}&\overline\GG''
		}
		\qquad
	\xymatrix{
		\OO K^{h\vee}\ar[d]^{\alpha}&\OO F^{kh\vee}\ar[l]_{\blin^\vee}\ar[d]^{\alpha'}&\OO F^{kh\vee}\ar[d]^{\alpha''}\ar[l]_{\ACg^\vee}\\
		\GG\ar[r]^{{\psi_1}}&\GG'\ar[r]^{{\psi_2}}&\GG''.
		}
	}	
Letting $\map{\psi_3}\GG{\GG''}$ be Serre's isogeny attached to $\bpdd{m_1+m_2}$, it suffices to prove	
$$\psi_3=\psi_2\circ\psi_1.$$  
By definition of $\map{\psi_1}{\GG}{\GG^\prime}$ and $\map{\psi_2}{\GG^\prime}{\GG^{\prime\prime}}$, we have
	\[equation]{\lb{fafafanimei}
\psi_2(\psi_1(X)) = \prod_{\vv\in U(g)}\left(\psi_1(X)[-]_{\GG^\prime}\alpha^\prime(\vv)\right)
	}
	where $U(g)=g^\vee(\pi^n\OO F^{kh\vee})/\pi^{n+m_2}\OO F^{kh\vee}$. Using the diagram \eqref{azi}, we have $\alpha^\prime = \psi_1\circ\alpha\circ\blin^\vee$, which allows us to write each factor of \eqref{fafafanimei} as
	$$\psi_1(X)[-]_{\GG^\prime}\alpha^\prime(\vv) = \psi_1(X)[-]_{\GG^\prime}\psi_1(\alpha(\blinn(\vv))).
$$
Denote $\ww:=\blin^\vee(\vv)$. Since $\psi_1$ is an isogeny from $\GG$ to $\GG^\prime$, the factor further simplifies to 
$$
	\psi_1(X)[-]_{\GG^\prime}\psi_1(\alpha\(\ww\)) = \psi_1\left(X[-]_{\GG}\alpha(\ww)\right) = \prod_{\vv\in U(\blin)}\left(X[-]_{\GG}\alpha(\ww+\vv)\right),
$$
where $U(\blin)=\pi^{n+m_2}\blinn(\OO F^{kh\vee})/\pi^{n+m_1+m_2}\OO K^{h\vee}$, and the last equality is obtained from the definition of $\psi_1$. 
Putting these simplified factors back to the original product \eqref{fafafanimei}, we obtain
$$
\psi_2\circ\psi_1(X) = \prod_{\ww\in U(g\blin,\blin)}\prod_{\vv\in U(\blin)}\left(X[-]_{\GG}\alpha(\ww+\vv)\right),
$$
where $U(g\blin,\blin)= \blinn\circ\ACg^\vee(\pi^n\OO F^{kh\vee})/\pi^{n+m_2}\blinn(\OO F^{kh\vee})$. 
 Therefore, by further simplifying this expression we obtain
$$
\psi_2\circ\psi_1(X) = \prod_{\vv\in   U(g\blin)}(X[-]_{\GG}\alpha(\vv)) = \psi_3(X),
$$
where $U(g\blin)=(\ACg\blin)^\vee(\pi^n\OO F^{kh\vee})/ \pi^{n+m_1+m_2}\OO K^{h\vee}$. This implies $\psi_3=\psi_2\circ\psi_1$, which completes the proof of the lemma.\end{proof}

The next lemma shows that the action of $\ACpair$ translates the cycle $Z_n(\biso,\blin)$ to $Z_n(\ACj\cdot\biso,\ACg\cdot\blin)$.
\[lem]{\label{wanquanbuzhi}\comments{wanquanbuzhi}
For any $n\geq 0$, any $(\biso,\blin)\in\equi KFh$, any $(\ACj,\ACg)\in\equi FF{kh}$, and any $m\geq\m(\ACg)$, we have
\begin{equation}\label{HAYA}\comments{HAYA}
	\frac1{\dpa nm}\bppb nm_{*}\ca{\biso}{\blin}{n+m}=
	\ca{\ACj\cdot\biso}{\ACg\cdot\blin}n.
\end{equation}
}

\begin{proof}
Let $M=\m(\blin)$. 
	By Definition \ref{tcycle}, it suffices to prove that
\[equation*]{\[split]{
	&\frac{1}{\dpa {n+m}M\dpa nm}\bppb nm_*\bppc {n+m}M_*\left[\OO{\DEf{n+M+m}}\right]\\=&\frac{1}{\dpa n{m+M}}\bpdd{m+M}_*\left[\OO{\DEf{n+M+m}}\right],
}}
which follows directly from the observation that
$$\dpa n{M+m}=\dpa{n+m}M\dpa nm,$$ 
and that
$$\bppb nm\circ\bppc {n+m}M=\bpdd{m+M}$$ by Proposition \ref{Composition}. \end{proof}

\section{An approximation for \Heegner cycles on the infinite level}\label{mo}\comments{mo}

  
In this section, we keep our general assumption that $[K:F]=k$, and use the same notation as in the previous section. As we described in Definition \ref{blinindu}, an element $(\biso,\blin)\in\equi KFh$ induces a system of compatible morphisms 
$$
\map{\bpp}{\DEf {m+n}}{\DEF n}
$$
for all $n$, and all $m\geq \m(\blin)$. Let $\GKH$ (resp.$\GFK$) be the $h$-(resp.$kh$-)-fold self-direct product of $\GK$ (resp.$\GF$) over $\Spec\CFF q$. In this section, we will define a similar map
$$
\map{\gpp}{\GKH}{\GFK}
$$
 in Definition \ref{aho}, and use it to construct \Heegner cycles $\cycy\infty$ on $\GFK$. We intentionally use $\infty$ as our subscript since $\cycy\infty$ should be thought of as an approximation of $\cycy n$ when $n\rightarrow\infty$. Our main result in this section is Proposition \ref{juda}, which shows a comparison identity
$$
\sss^*\cycy{n+m}=\ssss^*\cycy\infty\in \groo{\JJJJ}
$$
for the following closed embeddings for any $n>0$, 
$$
\xymatrix{
	\DEF {m+n}&&\JJJJ\ar[rr]^\ssss\ar[ll]_\sss&&\GFK,
}
$$ 
which will be defined in \eqref{fifinato} and \eqref{ififnato}. 

We need to introduce some constructions before proving our main results. Our strategy is to interpret $\GFK$ and $\GKH$ as moduli spaces of level structures of $\GF$ and $\GK$ respectively by writing 
$$\GKH
\cong\HGK$$
and
$$\GFK
\cong\HG,$$
where 
by $\HGK$ we mean the functor defined by assigning each $A\in\CCC$ to the $\OO K$-module
$$
\Hom_{\OO K}\left(\OO K^{h\vee},\GK(A)\right),
$$
and the meaning for 
$\HG$ is similar.

Since all formal $\OO K$-modules of height $h$ are isomorphic over $\CFF q$, throughout this subsection, without loss of generality, we may assume that $$[\pi]_{\GK}(X)=X^{q^{kh}}.$$ 
Moreover, by $\CCCC$ we mean the full subcategory of $\CCC$ consisting of those $A\in\CCC$ such that $\pi=0$ in $A$. In this section, all pairs $(\biso,\blin)$ and $(\ACj,\ACg)$ are equi-height pairs, and the integer $m$ is chosen so that $m\geq \m(\blin)$. Therefore 
$$
\mathrm{Height}(\biso^{-1})=-\mathrm{Height}(\blin)=-\log_q\Vol(\blinn(\OO F^{kh\vee}))\geq 0,
$$
and
$$
\mathrm{Height}(\pi^m\biso) =-\mathrm{Height}(\pi^m(\blin)^{-1})= \log_q\frac{\Vol(\blinn(\OO F^{kh\vee}))}{\Vol(\pi^m\OO K^{h\vee})}\geq 0,
$$
which implies that $\map{\biso^{-1}}\GK\GF$ and $\map{\pi^m\biso}\GK\GF$ are actual isogenies.

\subsection{Maps from $\GKH$ to $\GFK$}
\begin{defi}\label{aho}\comments{aho} For any $(\biso,\blin)\in\equi KFh$ and $m\geq\m(\blin)$, let $\map\gpp\GKH\GFK$ be the quasi-homomorphism of formal $\OO F$-modules defined 
functorially in $A\in\CCC$
by
\begin{equation}
\maps{\gpp}{\Hom_{\OO K}(\OO K^{h\vee},\GK(A))}{\Hom_{\OO F}(\OO F^{kh\vee},\GF(A))}{f}{\pi^m\biso\circ f\circ\blinn.}
\end{equation}
Moreover, let $\J$ be the group scheme sitting inside the following exact sequence
$$
	\xymatrix{
		0\ar[rr]&&\J\ar[rr]&& \GKH\ar[rr]^\gpp&&\GFK.
	}$$
 \end{defi}
\begin{rem}
	For any $(\ACj,\ACg)\in\equi FF{kh}$ and $m\geq\m(\ACg)$, functorially in $A\in\CCC$, we define the map $\gpb$ by 
\begin{equation}
\maps{\gpb}{\Hom_{\OO F}(\OO F^{kh\vee},\GF(A))}{\Hom_{\OO F}(\OO F^{kh\vee},\GF(A))}{f}{\pi^m\ACj\circ f\circ\ACg^\vee.}
\end{equation}
Moreover, $\JJ$ is the group scheme sitting inside the following exact sequence
$$
	\xymatrix{
		0\ar[rr]&&\JJ\ar[rr]&& \GFK\ar[rr]^\gpb&&\GFK.
	}$$
\end{rem}
\begin{prop}[Analogue to Proposition \ref{Composition}]\label{haha}\comments{haha}
For any $(\biso,\blin)\in\equi KFh$, $(\ACj,\ACg)\in\equi FF{kh}$, $m_1\geq\m(\ACg)$ and $m_2\geq\m(\blin)$, we have 
$$
	\gppbj{}{m_1}\circ\gppj{}{m_2}=\gppdj{}{m_1+m_2}.
$$
In other words, the following diagram commutes

$$
	\trii {\GFK}{\GFK.}{\GKH}{\gppbj{}{{m_1}}}{\gppdj{}{{m_1+m_2}}}{\gppj{}{m_2}}
$$

\end{prop}
\begin{proof} For any $A\in\CCC$, and any $f\in\GKH(A)$, by Definition \ref{aho}, we have
$$
	\gppbj{}{m_1}\circ\gppj{}{m_2}(f) = \gppbj{}{m_1} (\pi^{m_2}\biso\circ f\circ \blinn) = \pi^{m_1+m_2} \ACj\circ\biso\circ f\circ\blinn\circ \ACg^\vee,
$$
which equals
$$
	\pi^{m_1+m_2}(\ACj\cdot\biso)\circ f\circ (\ACg\cdot\blin)^\vee = \gppdj{}{m_1+m_2}(f)
$$
as desired.	
\end{proof}

\subsection{Thickening comparison}
This part is the technical core of the article. Loosely speaking, thinking of $\GKH$ as the deformation space on the infinite level, we temporarily denote $\DEf\infty:=\GKH$ and $\DEF\infty:=\GFK$. Then the exact sequence in Definition \ref{aho} is
$$
	\xymatrix{
		0\ar[rr]&&\J\ar[rr]&& \DEf\infty\ar[rr]^\gpp&&\DEF\infty.
}$$
Since intuitively $\DEF\infty$ is an approximation of $\DEF n$ when $n$ is sufficiently large, one may carelessly expect an exact sequence
$$
	\xymatrix{
		0\ar[rr]&&\J\ar[rr]&& \DEf{n+m}\ar[rr]^\bpp&&\DEF n
}$$
for sufficiently large $n$. Going back to the rigorous math, the exact sequence here does not make any sense since $\DEF n$ is not a group scheme. Nevertheless, we realize this idea in Proposition \ref{haolemma} by identifying $\J$ with the preimage of the closed point $\Spec(\CFF q)$ under the map $\bpp: {\DEf{n+m}\longrightarrow\DEF n}$, which is the main result of this subsection.

 
Before the construction of the map $\J\longrightarrow\DEf{m+n}$, we need some lemmas for level structures of $\GK$.
\begin{lem}\label{LEMM:REDUC}\comments{LEMM:REDUC}
All elements in
$$\Hom\big({\OO K^{h\vee}},\GK[\pi^n](A)\big)$$
are $\pi^{n+1}$-level structures of $\GK$ over $A$ for any $A\in\CCCC$.
\end{lem}
\begin{proof}
Suppose $A\in\CCCC$. We have $\pi=0\in A$.
Let $f\in\Hom\big({\OO K^{h\vee}},\GK[\pi^n](A)\big)$. 
It suffices to verify the non-degeneracy condition for $f$,
that is, to show that 
\begin{equation}\label{rea}\comments{rea}\prod_{\ww\in\OO K^{h\vee}/\pi^{n+1}\OO K^{h\vee}}(X-f(\ww))=[\pi^{n+1}]_{\GK}(X).\end{equation}
We prove it by induction on $n$. If $n=0$, the expression (\ref{rea}) is clearly true since $f=0$. 
Assuming the following induction hypothesis 
\begin{equation}\lb{hshs}\prod_{\ww\in\OO K^{h\vee}/\pi^{n}\OO K^{h\vee}}(X-f(\ww))=[\pi^{n}]_{\GK}(X).\end{equation}	
	Abbreviating $[\pi^n]_{\GK}$ to $[\pi^n]$, the right hand side of \eqref{rea} simplifies to 
$$
[\pi^{n+1}](X)=[\pi^n]([\pi](X)),
$$
	which, by the induction hypothesis \eqref{hshs}, equals
$$
\prod_{\ww\in\OO K^{h\vee}/\pi^{n}\OO K^{h\vee}}\left([\pi]X-[\pi]f\left(\ww\right)\right).
$$
Since $[\pi](X)=X^{q^{kh}}$, each factor of the preceding product simplifies to  
$$
[\pi]X-[\pi]f\left(\ww\right)= X^{q^{kh}}-f(\ww)^{q^{kh}} = (X-f(\ww))^{q^{kh}}.
$$
Putting the simplified factors back into the product, we obtain
$$\prod_{\ww\in\OO K^{h\vee}/\pi^{n}\OO K^{h\vee}}(X-f(\ww))^{q^{kh}}=\prod_{\ww\in\OO K^{h\vee}/\pi^{n+1}\OO K^{h\vee}}(X-f(\ww)),$$
	which is the left hand side of \eqref{rea}. Here the last equality is obtained from the observation that $f(\ww+\vv)=f(\ww)+f(\vv)=f(\ww)$ for any $\vv\in\pi^n\OO K^{h\vee}$. 
This completes the proof.\end{proof}

\begin{lem}
	All elements $f\in \J(A)$ are $\pi^{m+n}$-level structures of $\GK$ over $A$ for any integer $n>\m(\blin)$, any integer $m\geq\m(\blin)$ any $(\biso,\blin)\in\equi KFh$, and any $A\in\CCCC$.
\end{lem}
\begin{proof}
	Since $\mathrm{Height}(\biso^{-1})\geq 0$, the map
	\begin{equation}\label{backone}
	\map{\biso^{-1}}\GK\GF
	\end{equation}
	is an actual isogeny. Since $n\geq\m(\blin)+1$, we have $\blinn(\OO F^{kh\vee})\supset \pi^{n-1}\OO K^{h\vee}$, which implies that 
	$$\pi^{n-1}(\tau^{\vee})^{-1}\left(\OO K^{h\vee}\right)\subset\OO F^{kh\vee}.$$
	Therefore $\pi^{n-1}(\blin^{\vee})^{-1}$ restricts to a map
	\begin{equation}\label{backtwo}
		\map{\pi^{n-1}(\blin^{\vee})^{-1}}{\OO K^{h\vee}}{\OO F^{kh\vee}}.
	\end{equation}
	By Definition \ref{aho}, $f\in \J(A)$ implies that $\pi^m\biso\circ f\circ\blinn=0$, which can be put into the following commutative diagram
$$
	\xymatrix{
		\mathcal O_F^{kh\vee}\ar[rr]^{\blinn}\ar[d]_0&&\mathcal O_K^{h\vee}\ar[d]^f\\
		\GF&&\GK.\ar[ll]^{\pi^m\biso}
		}
$$
Extending this diagram using the maps in \eqref{backone} and \eqref{backtwo}, we have the following commutative diagram
$$
	\xymatrix{
		\mathcal O_K^{h\vee}\ar@{-->}@/^1pc/[rrrr]^{\pi^{n-1}}\ar[rr]_{\pi^{n-1}(\blin^{\vee})^{-1}}\ar[d]_0&&\mathcal O_F^{kh\vee}\ar[rr]_{\blinn}\ar[d]_0&&\mathcal O_K^{h\vee}\ar[d]^f\\
		\GK&&\GF\ar[ll]_{\biso^{-1}}&&\GK.\ar[ll]_{\pi^m\biso}\ar@{-->}@/^1pc/[llll]^{\pi^m}
		}
$$
Note that this diagram implies that $[\pi^{m+n-1}]_{\GK}\circ f=0$, and so that $f$ is a $\pi^{m+n}$-level structure by Lemma \ref{LEMM:REDUC} , which completes the proof.\end{proof}
The formal scheme $\GK^h$ can be naturally thought of as a functor from $\CCC$ to the category of sets by extending its definition from $\CCCC$ to $\CCC$, where we define $\GK^h(A)=\emptyset$ for any $A\notin\CCCC$.

\begin{defi}\lb{rafara}
	For any $A\in\CCCC$, by $\GG_{K/A}$ (resp. $\GG_{F/A}$) we mean the formal $\OO K$ (resp.$\OO F$) module obtained by the base change of $\GK$ (resp. $\GF$) to $A$.
\end{defi}
Note that $\GK$ is defined by a lot of power series describing its addition and its scalar multiplication by $\OO K$. Those power series have coefficients in $\CFF q$. The formal module $\GG_{K/A}$ has exactly the same power series with $\GK$, but we view those coefficients as coefficients in $A$ through $\CFF q\subset A$. Therefore, any endomorphism $\map{\ACj}{\GK}\GK$ can be directly viewed as an endomorphism of $\GG_{K/A}$, and by abuse of notation, we use the same symbol $\map{\ACj}{\GG_{K/A}}{\GG_{K/A}}$ to represent this endomorphism.
\begin{defi}\label{Thes}\comments{Thes}
We define the map $\J\longrightarrow\DEf{m+n}$ functorially in $A\in\CCC$ by 
	\begin{equation}\label{FIBE:DEFIN}\comments{FIBE:DEFIN}\mapi {\J(A)}{\DEf{m+n}(A)}{f}{[\GG_{K/A},\id,f]_{m+n}}\end{equation}
for any $n> \m(\blin)$ and any $(\biso,\blin)\in\equi KFh$.
\end{defi}

The following lemma will be used in the proof of Proposition \ref{haolemma}.
\[lem]{\lb{stupid}For any $A\in\CCCC$ and $[\GG_{K/A},\ACj,f]_n\in\DEf n(A)$, we have 
$$[\GG_{K/A},\ACj,f]_n=[\GG_{K/A},\id,\ACj^{-1}\circ f]_n.$$}
\[proof]{Note that $\DEf n$ is an abbreviation of $\DEf n^{(0)}$ (for notation see Theorem \ref{LTlevel}). Since $[\GG_{K/A},\ACj,f]_n\in\DEf n^{(0)}(A)$, the morphism $\map\ACj\GK{\overline{\GG_{K/A}}}$ is of height $0$, and thus an isomorphism. Since $\overline{\GG_{K/A}}=\GK$, the isomorphism $\ACj$ is an automorphism of $\GK$. By base change to $A$, the isomorphism $\ACj:\GG_{K/A}\rightarrow\GG_{K/A}$ identifies the triple $(\GG_{K/A},\id,\ACj^{-1}\circ f)$ to the triple $(\GG_{K/A},\ACj,f)$. Therefore these two triples are equivalent.}


\begin{prop}\label{haolemma}\comments{haolemma} In the following diagrams, all unlabeled vertical maps are defined by the same way as in Definition \ref{Thes}. We have:
\begin{enumerate}
	\item[(1)]For any $(\biso,\blin)\in\equi KFh$, any $m\geq\m(\blin)$, and any $n>\m(\blin)$, the following diagram
\begin{equation}\label{shaya}\comments{shaya}
\cate {\J}{\Spec\CFF q}{\DEf{m+n}}{\DEF n}{}{}{}{\bpp}
\end{equation}
		is Cartesian. In particular, the vertical maps are closed embeddings.
\item[(2)]  For any
$(\biso_1,\blin_1)\in\equi KFh$, any $(\ACj,\ACg)\in\equi FF{kh}$, any $m_1>\m(\blin_1)$, and any $m_2>\m(\ACg)$, we denote $(\biso_3,\blin_3):=(\ACj\biso_1,\ACg\blin_1)\in\equi KFh$, and $m_3:=m_1+m_2$. Then the following diagram
\begin{equation}\label{fanren}\comments{fanren}
	\cate {{\Jh}}{{\mathcal G^{kh}_F[\pi^{m_2}\ACj\ACg]}}{\DEf{n+m_1}}{\DEF{n}}{\gpo}{\dla{s(\pi^{m_3}\biso_3,\blin_3)_{n-m_2}}}{\dla{s(\pi^{m_2}\biso_2,\blin_2)_{n-m_2}}}{\bppo}
\end{equation}
		is Cartesian for any $n>m_2+\m(\blin_3)$.
\end{enumerate}
\end{prop}
\begin{proof}[Proof of statement (1)] It suffices to prove that the diagram (\ref{shaya}) is Cartesian. Then the morphism $\J\longrightarrow\DEf{m+n}$ is automatically a closed embedding since it is obtained by base change of the closed embedding $\Spec\CFF q\rightarrow \DEF n$. In other words, it suffices to show that for any $A\in\CCC$ and any $[\GG,\iota,\alpha]_{m+n}\in \DEf{m+n}(A)$, 
\[equation]{\label{iffa}{\bpp}[\GG,\iota,\alpha]_{m+n}=[\GFA,\id,0]_n}
	if and only if
	\[equation]{\label{iffb}[\GG,\iota,\alpha]_{m+n}=[\GKA,id,f]_{m+n}\;\text{ for some }f\text{ such that }\;\pi^m\biso\circ f\circ \blinn=0.}
	We prove that \eqref{iffb} implies \eqref{iffa} as follows. Let $n'=\max\{n,m\}$. Let $\map{\psi}{\GKA}{\GG^\prime}$ be Serre's isogeny attached to $\eta^{m+n}_{n'}(\biso,\blin)$ and $\GKA$. By definition,
$$\psi(X)=\prod_{\ww\in V}\Big(X[-]_{\GKA} f(\ww)\Big),$$
	where $V=\blinn(\pi^{n'}\OO F^{kh\vee})/\pi^{n+m}\OO K^{h\vee}$.
	Since $\pi^m\biso\circ f\circ \blinn=0$ and $\mathrm{Height}(\bisoo)\geq 0$, we have 
$$
	f\circ\blin^\vee(\pi^{n'}\vv)=\pi^{n'-m}\cdot\biso^{-1}(\pi^m\biso\circ f\circ\blinn(\vv))=\pi^{n'-m}\cdot\biso^{-1}(0)=0
$$
for any $\vv\in\OO F^{kh\vee}$,
which implies that $f(\ww)=0$ for any $\ww\in V$. Therefore
$$\psi(X)=\prod_{\ww\in V}\Big(X[-]_{\GKA} f(\ww)\Big)=X^{\# V}.$$
Since
$$
	\psi\circ[\pi]_{\GKA}(X)=X^{q^{kh+\# V}}=[\pi]_{\GFA}\circ\psi(X),
$$
	we have $\GG'=\GFA$. Thus we can write ${\eta^{m+n}_{n'}(\biso,\blin)}[\GKA,\id,f]_{n+m}=[\GFA,\iota,\alpha']_{n'}$ for some $\pi^{n'}$-level structure $\alpha':\OO F^{kh\vee}\rightarrow\GFA(A)$ and some automorphism $\iota:\GF\rightarrow\GF$. 
	Then $$\eta^{n+m}_n(\biso,\blin)[\GKA,\id,f]_{n+m}=\beta^{n'}_n\circ\eta^{n+m}_{n'}(\biso,\blin)[\GKA,\id,f]_{n+m} = \beta^{n'}_n[\GFA,\iota,\alpha']_{n'}=[\GFA,\iota,\alpha]_n,$$  
	where $\alpha = \alpha'\circ \pi^{n'-n}$. By Corollary \ref{rmindu}, we obtain an isogeny $\map{\psi'}\GKA\GFA$ such that the following diagrams commute
	\[equation]{\lb{buzhidaozenmeban}
	\xymatrix{
		\GK\ar[rr]^{\pi^m\biso}\ar[d]_{\id}&&\GF\ar[d]^{\iota}\\
		\GK\ar[rr]^{\overline{\psi'}}&&\GF,
		}
\qquad
	\xymatrix{
		\OO K^h\ar[d]_{f}&&\OO F^{kh}\ar[ll]_{\blin^\vee}\ar[d]^{\alpha}\\
		\GKA\ar[rr]^{{\psi'}}&&\GFA.
		}
	}
Note that the lifting of endomorphisms of $\GF$ is unique, so that $\overline{\psi'}=\psi'$.
From the left diagram of \eqref{buzhidaozenmeban}, one sees that
$$
	\overline{\psi'}=\iota\circ \pi^m\biso.
$$
Since $\pi^m\biso\circ f\circ \blin^\vee = 0$, from the right diagram of \eqref{buzhidaozenmeban}, we obtain that
$$
\alpha = \psi'\circ f\circ \blin^\vee = \iota\circ \pi^m\biso\circ f\circ \blin^\vee=\iota\circ 0=0,
$$
which implies that
$${\bpp}[\GKA,\id,f]_{n+m}=[\GFA,\iota,0]_n.$$
	Furthermore, note that $\mathrm{Height}(\blin)=\mathrm{Height}(\biso)$. Therefore $\map{\pi^m\biso}\GFA\GKA$ and $\map{\psi'}\GK\GF$ are elements of the same height in $\End(\GF)$, which implies that $\iota=\psi'\circ(\pi^m\biso)^{-1}$ is an isomorphism of $\GF$. Therefore, by Lemma \ref{stupid}, $[\GFA,\iota,0]_n=[\GFA,\id,0]_n$, which proves \eqref{iffa}.

	Conversely, let $$\map{\biso_0}{\GK}{\GF}, \quad \map{\blinn_0}{\OO F^{kh\vee}}{\OO K^{h\vee}}$$ be natural forgetful isomorphisms defined by only remembering $\OO F$-actions. Then there exists $\ACj\in D_F$ and $\ACg\in\GL_{kh}(F)$ such that the following two diagrams 
$$
	\triso\GK\GF{\GF,}{\biso_0}{\biso}{\ACj}\qquad\qquad 
	\xymatrix{
		{K^{h\vee}}&&&&\ar[llll]_{\blin_0^\vee}^{\cong}{F^{kh\vee}}\\
		&&{F^{kh\vee}}\ar[llu]^{\blinn}\ar[rru]_{\ACg^\vee}
	}
$$
commute.
Since $\mathrm{Height}(\blin_0)=\mathrm{Height}(\biso_0)=0$ and $\mathrm{Height}(\blin)=\mathrm{Height}(\biso)$, the pair $(\ACj,\ACg)$ is an equi-height pair. 
	Letting $u=\m(\blin)+1$, we have $\OO F^{kh\vee}\supset g^\vee(\OO F^{kh\vee})\supset \pi^{u-1}\OO F^{kh\vee}$, which implies that 
	$$
	\pi^{u-1}\OO F^{kh\vee}\subset \pi^{u-1}\ACg^{-1\vee}\left(\OO F^{kh\vee}\right)\subset \OO F^{h\vee}.
	$$
Therefore $\pi^{u-1}\ACg^{-1}\in\matt{kh}{kh}{\OO F}$, and 
$$\m(\pi^{u-1}\ACg^{-1})\leq u-1.$$
	Since $\mathrm{Height}(\pi^{1-u}\ACj^{-1})=\mathrm{Height}(\pi^{u-1}\ACg^{-1})$, the pair $(\pi^{1-u}\ACj^{-1},\pi^{u-1}\ACg^{-1})$ is an equi-height pair.  
Now we prove that \eqref{iffa} implies \eqref{iffb}. Recall that \eqref{iffa} is ${\bpp}[\GG,\iota,\alpha]_{n+m}=[\GFA,\id,0]_n$. Composing this identity with $\beta_{n-u+1}^n(\pi^{1-u}\ACj^{-1},\pi^{u-1}\ACg^{-1})$, we obtain
\begin{equation}
\begin{split}
&\beta_{n-u+1}^n(\pi^{1-u}\ACj^{-1},\pi^{u-1}\ACg^{-1})\circ{\bpp}[\GG,\iota,\alpha]_{n+m}\\=\;&\beta_{n-u+1}^n(\pi^{1-u}\ACj^{-1},\pi^{u-1}\ACg^{-1})[\GFA,\id,0]_{n}\\
	=\;&
	[\GFA,\overline\psi\circ\id\circ\pi^{u-1}(\pi^{1-u}\ACj^{-1})^{-1},\psi\circ0\circ (\pi^{1-u}\ACg^{-1})^\vee]_{n-u+1}\\
	=\;&
[\GFA,\id,0]_{n-u+1},
\end{split}
\end{equation}
where $\map\psi\GFA\GFA$ is the Serre's isogeny attached to $(\pi^{1-u}\ACj^{-1},\pi^{1-u}\ACg^{-1})$ and the last equality follows from Lemma \ref{stupid}. Moreover, this expression can be calculated in another way as follows,
\begin{equation}
\begin{split}
	&\beta_{n-u+1}^n(\pi^{1-u}\ACj^{-1},\pi^{u-1}\ACg^{-1})\circ {\bpp}[\GG,\iota,\alpha]_{n+m}\\
	=\;&\bccc{n-u+1}{n+m}{\pi^{1-u}\ACj^{-1}\biso}{\pi^{u-1}g^{-1}\blin} [\GG,\iota,\alpha]_{n+m} \\
	=\;&\bccc{n-u+1}{n+m}{\pi^{1-u}\biso_0}{\pi^{u-1}\blin_0}[\GG,\iota,\alpha]_{n+m}\\
	=\;&
[\GG,\iota\circ\bisoo_0,\alpha\circ\pi^{m+u-1}\tau_0^\vee]_{n-u+1}.
\end{split}
\end{equation}
Therefore, we have $[\GG,\iota\circ\bisoo_0,\alpha\circ\pi^{m+u-1}\tau_0^\vee]_{n-u+1}=[\GFA,\id,0]_{n-u+1}$, which implies that two triples $(\GG,\iota\circ\bisoo_0,\alpha\circ\pi^{m+u-1}\tau_0^\vee)$ and $(\GFA,\id,0)$ are equivalent as formal $\OO F$-modules. Recall that $\GG$ is deformed from $\GK$, and that $\GF$ is obtained by only remembering $\OO F$-actions of $\GK$. Therefore, in fact, 
$$
[\GG,\iota,\alpha\circ\pi^{m+u-1}\tau_0^\vee]_{n-u+1}=[\GKA,\id,0]_{n-u+1},
$$
which implies that we have an isomorphism $\mapp\GG\GKA$ of formal $\OO K$-modules such that the following diagrams commute
$$
\xymatrix{
	&\GK\ar[ld]_{\iota}\ar[rd]^{\id}\\
	\GG\ar[rr]&&\GK,
	}
\qquad
\xymatrix{
	&\OO K^{h\vee}\ar[ld]_{\alpha\circ \pi^{m+u-1}\blin_0^\vee}\ar[rd]^{0}\\
	\GG\ar[rr]&&\GKA.
	}
$$
Letting $f=\iota^{-1}\circ\alpha$, we have
\[equation]{\lb{eqeqa}
	[\GG,\iota,\alpha]_{n+m}=[\GKA,\id,f]_{n+m}.
}
Therefore by the assumption ${\bpp}[\GG,\iota,\alpha]_{n+m}=[\GFA,\id,0]_n$, we obtain
$$
{\bpp}[\GKA,\id,f]_{n+m}=[\GFA,\id,0]_n,
$$
which implies that there is an isogeny $\map{\psi}\GKA\GFA$ such that the following diagrams commute
$$
\xymatrix{
	\GK\ar[rr]^{\pi^m\biso}\ar[d]_\id&&\GF\ar[d]^\id\\
\GK\ar[rr]^{\psi}&&\GF,
	}
\qquad\qquad
\xymatrix{
	\OO K^{h\vee}\ar[d]_f&&\OO F^{kh\vee}\ar[ll]_{\blinn}\ar[d]^0\\\
\GKA\ar[rr]^{\psi}&&\GFA.
	}
$$
Therefore, we have $\psi=\pi^m\biso$ and $\psi\circ f\circ\blinn=0$, which implies that $\pi^m\biso\circ f\circ \blinn=0$. Combining this with \eqref{eqeqa}, we complete the proof of \eqref{iffb}.

\end{proof}
\begin{proof}[Proof of statement (2):]
	Before starting the proof, we extend the diagram \eqref{fanren} into the following diagram
$$
	\xymatrix{
		\Jh\ar[rr]^\gpo\ar[d]&&\mathcal G_F^{kh}[\pi^{m_2}\ACj\ACg]\ar[rr]\ar[d]&&\Spec \CFF q\ar[d]\\
		\DEf{n+m_1}\ar[rr]^\bppo&&\DEF n\ar[rr]^{\beta^{n}_{n-m_2}(\ACj,\ACg)}&&\mathcal \DEF_{n-m_2}.
		}
$$
	Here the extended right square is constructed in the same way as \eqref{shaya}, which is a Cartesian diagram. 
If the left square commutes, then the whole diagram commutes, and therefore the left square is automatically a Cartesian diagram because the right and the outer squares are Cartesian by statement (1).
Therefore, it suffices to show that the left square commutes, that is, for any $A\in\mathcal C\otimes\CFF q$ and $f\in\Jh(A)$, one has
$$\bppo[\GKA,\id,f]_{n+m_1}=[\GFA,\id,\pi^{m_1}\biso_1\circ f\circ\blinn_1]_n.$$
Since the outer diagram is commutative, we have
\[equation]{\label{shower}\beta_{n-m_2}^n(\ACj,\ACg)\Big(\bppo[\GKA,\id,f]_{n+m_1}\Big)=[\GFA,\id,0]_{n-m_2}.}
Since the right square is Cartesian by the first statement of Proposition \ref{haolemma}, the equation \eqref{shower} implies that
$$\bppo[\GKA,\id,f]_{n+m_1}=[\GFA,\id,g]_n$$
	for some $g\in\Hom_{\OO F}(\OO F^{kh\vee},\GF(A))$. Therefore, by Corollary \ref{rmindu}, there exists an isogeny $\map\psi\GKA\GFA$ such that the following diagrams commute
$$
	\xymatrix{
		\GK\ar[rr]^{\pi^{m_1}{\biso_1}}\ar[d]_\id&&\GF\ar[d]^\id\\
		\GK\ar[rr]^{\psi}&&\GF\\
		}
\qquad\qquad
	\xymatrix{
		\mathcal O_K^{h\vee}\ar[d]_f&&\mathcal O_F^{kh\vee}\ar[ll]_{\blinn}\ar[d]^{g}\\
		\GKA\ar[rr]^{\psi}&&\GFA,\\
		}
$$
	which imply that $g=\psi\circ f\circ \blinn$ and $\psi=\pi^{m_1}\biso_1$.
Therefore we have $g=\pi^{m_1}\biso_1\circ f\circ \blinn$, which completes the proof.\end{proof}

\subsection{\Heegner cycles comparison}\label{capa}\comments{capa}
Now we define \Heegner cycles on $\GFK$ by similar ways as in Definition \ref{tcycle}.
\begin{defi}\lb{infinf}
	For any $(\biso,\blin)\in\equi KFh$, let $\cycy\infty\in \groo{\GFK}$ be the element obtained by  
	\begin{equation}\label{DEFICY}\comments{DEFICY}\cycy\infty:=\frac1{\dgpa m}\Big[{\gpp}_*\OO{\GKH}\Big].\end{equation}
\end{defi}
Clearly, this definition does not depend on $m$ by Proposition \ref{haha}.

Now we will compare cycles $\cycy n$ with $\cycy\infty$. Let $\ssss$ be the natural inclusion map 
\[equation]{\lb{fifinato}\map{\ssss}{\JJJ}{\GFK},}
and  for any $n\geq 1$ let $\sss$ be the map 
\[equation]{\lb{ififnato}\maps{\sss}{\JJJ}{\DEF {n+m}}{f}{[\GF,\id,f]_{n+m}}}
as defined in Definition \ref{Thes}.
The following proposition is the main result of this section.
\begin{prop}\label{juda}\comments{juda}
For any integer $m\geq 0$, $n>\m(\blin)$ and any $(\biso,\blin)\in\equi KFh$, we have
$$
\sss^*\cycy{n+m}=\ssss^*\cycy\infty\in\groo{\JJJ}.
$$\end{prop}
\begin{proof}
By definition of $\cycy\infty$ and $\cycy{n+m}$ in \eqref{DEFICY} and \eqref{zhaogebishengtianhaodelaopo}, it suffices to prove 
	\begin{equation}\[split]{
			&\frac1{\dpa{n+m}w}\Big[\sss^*\bppc{n+m}w_*\OO{\DEf{n+m+w}}\Big]\\=&\frac1{\dgpa w}\Big[\ssss^*\bppc\infty w_*\OO{\GKH}\Big]
		}
\end{equation}
for some $w\geq\m(\blin)$. Since $n>0$, we observe that the denominator of two sides are equal because
$$\dpa{n+m}w=q^{kh^2w}=\dgpa w.$$ It therefore suffices to show that
\begin{equation}\label{yuwang}\comments{yuwang}
\sss^*\bppc{n+m}w_*\OO{\DEf{n+m+w}}\cong\ssss^*\bppc\infty w_*\OO{\GKH}.
\end{equation}
Since $n>\m(\blin)$, by Proposition \ref{haolemma}, we obtain the following Cartesian diagram
$$
	\cate {\Jw}{{\JJJ}}{\DEf{n+m+w}}{\DEF{n+m},}{\bppc\infty w}{\ssc n{m+w}}{\sss}{\bppc{n+m} w}
$$
where all vertical maps,which are closed embeddings, are defined in the same way as in \eqref{FIBE:DEFIN} and labeled as shown in the diagram. We then obtain the following isomorphisms for the left hand side of \eqref{yuwang},
	\begin{equation}\lb{lafta}
\begin{split}
\sss^*\bppc{n+m}w_*\OO{\DEf{n+m+w}}\cong&\bppc\infty w_*\ssc n{m+w}^*\OO{\DEf{n+m+w}}\\\cong&
\bppc\infty w_*\OO{\Jw}.
\end{split}
\end{equation}
On the other hand, by using the following Cartesian diagram,
$$
\cate {\Jw}{\JJJ}{\GKH}{\GFK,}{\bppc\infty w}{\sgc\infty{m+w}}{\ssss}{\bppc\infty w}
$$
where all vertical maps are canonical closed embeddings and are labeled as shown, we obtain the following isomorphisms for the right hand side of \eqref{yuwang}
	\begin{equation}\lb{raita}
\begin{split}
	\ssss^*\bppc\infty w_*\OO{\GKH}\cong&\bppc\infty w_*\sgc\infty{m+w}^*\OO{\GKH}\\\cong&
\bppc\infty w_*\OO{\Jw}.
\end{split}
\end{equation}
Comparing \eqref{lafta} with \eqref{raita}, we complete the proof that two sides of \eqref{yuwang} are equal, and hence complete the proof of this proposition.\end{proof}

\section{Intersection comparison}\label{mi}\comments{mi}
From now, we take $k=2$, and then $K/F$ is a quadratic extension. In the previous section, we proved that there are two closed embeddings of formal schemes
\[equation]{\lb{location}
\xymatrix{\G&&\JJJJJ\ar[ll]_{\ssssss}\ar[rr]^{\sssss}&&\DEF n,}
}
and that our \Heegner cycles $\cycy\infty$ and $\cycy n$ satisfy
\[equation]{\label{hij}
\sssss^*\cycy{n} = \ssssss^*\cycy{\infty}.
}
Intuitively, the cycle $\cycy\infty$ is an approximation of $\cycy{n}$ for $n\rightarrow\infty$. Loosely speaking, when $n$ is sufficiently large, one can choose a large enough $M$ in the diagram \eqref{location} such that the intersection $\cycn1\infty\otimes_{\GF^{2h}}\cycn2\infty$ is contained completely inside the formal neighborhood $\GF^{2h}[\pi^M]$. Then from the diagram \eqref{location} and the equation \eqref{hij}, we identify the intersection numbers on sufficiently high levels with the one on the infinite level together.
The purpose of this section is to gather commutative algebra arguments to verify this intuition.

From this section onwards, we consider two quadratic extensions $K_1$, $K_2$ of $F$, which are not necessarily isomorphic. Our focus in this section is to prove the following key theorem:
\begin{thm}[Intersection Comparison]\label{THM:ICT}\comments{THM:ICT}
For any $(\biso_1,\blin_1)\in\equi{K_1}Fh$ and $(\biso_2,\blin_2)\in\equi {K_2}Fh$, if the length of $\cyio\otimes\cyit$ is finite, then there exists $N>0$(see \eqref{howbigN}) such that for all $n\geq N$, we have
\begin{equation}\chi(\cycn1n\litimes_{\DEF n}\cycn2n)=\chi(\cycn1\infty\litimes_{\GF^{2h}}\cycn2\infty).
\end{equation}
\end{thm}
\begin{rem}
The reader might be interested in the case where the length of $\cyio\otimes\cyit$ is infinite, in which case, it is easy to show that the length of $\cycn10\litimes_{\DEF 0}\cycn20$ is infinite as well. We omit the proof here since this case is not important for the application of the linear AFL. However, the author does not know if it is true that the length of $\cycn1n\litimes_{\DEF 0}\cycn2n$ is infinite as well for $n\geq 1$.
\end{rem}
Since we can regard $\G$ as a formal scheme over $\CO F$ by the morphisms $$\G\longrightarrow \Spf\CFF q\longrightarrow \Spf \CO F,$$
we only need to work in the category of formal schemes over $\CO F$, and hence we abbreviate $\length_{\CO F}(\bullet)$ to $\length(\bullet)$ for the rest of the article.

\subsection{Outline of proof}
We will prove Theorem \ref{THM:ICT} in 3 steps. 

\noindent\textbf{Step 1}: We reduce the calculation of intersection numbers to intersection multiplicities by proving the following identities:
\begin{equation}\label{fafa}\comments{fafa}
\chi(\cycn1n\litimes_{\DEF n}\cycn2n)=\length(\cycn1n\otimes_{\DEF n}\cycn2n),
\end{equation}
\begin{equation}\label{meme}\comments{meme}
\chi(\cycn1\infty\litimes_{\GF^{2h}}\cycn2\infty)=\length(\cycn1\infty\otimes_{\G}\cycn2\infty)
\end{equation}
when the right hand sides of them are finite. 

\noindent\textbf{Step 2}: By using the following closed embeddings, 
\begin{equation}\label{huijia}\comments{huijia}
\map{\sssss}{\JJJJJ}{\DEF n},
\end{equation} 
\begin{equation}\label{jiahui}\comments{jiahui}
\map{\ssssss}{\JJJJJ}{\G},
\end{equation} 
we prove that 
\begin{equation}\label{kunjuan}\comments{kunjuan}
\begin{split}
&\length\Big(\sssss_*\sssss^*\big(\cycn1n\otimes_{\DEF n}\cycn2n\big)\Big)\\=&\length\Big(\ssssss_*\ssssss^*\big(\cycn1\infty\otimes_{\G}\cycn2\infty\big)\Big)
\end{split}
\end{equation}
for any $n>\max\left(\m(\blin_1),\m(\blin_2)\right)+M$ by Proposition \ref{juda}.

\noindent\textbf{Step 3}: 
We prove that the approximation of \eqref{huijia} and \eqref{jiahui} will eventually catch the whole intersection as $M\rightarrow\infty$, in other words,
\begin{equation}\label{eshu}\comments{eshu}
\sssss_*\sssss^*\big(\cycn1n\otimes_{\DEF n}\cycn2n\big)=\cycn1n\otimes_{\DEF n}\cycn2n,
\end{equation} 
\begin{equation}\label{ishu}\comments{ishu}
\ssssss_*\ssssss^*\big(\cycn1\infty\otimes_{\G}\cycn2\infty\big)=\cycn1\infty\otimes_{\G}\cycn2\infty
\end{equation}
for any $n>M$, where $M$ is a large integer depending only on $(\biso_1,\blin_1)$ and $(\biso_2,\blin_2)$. With this choice of $M$, and letting $N=M+\max\left(\m(\blin_1),\m(\blin_2)\right)+1$, the finiteness assumption for the length of $\cycn1n\otimes_{\DEF n}\cycn2n$ in Step 1 follows from the finiteness assumption for the length of $\cycn1\infty\otimes_{\G}\cycn2\infty$, and therefore we will complete the proof of Theorem \ref{THM:ICT} once we finished these three steps.

\subsection{Step 1: Reduce to intersection multiplicities}\lb{euler}
Note that the following identity
\[equation]{\lb{intdefi}
\chi(\mathcal{F}\litimes_X \mathcal{G})=\sum_{i,j} (-1)^{i+j}\length_{\OO F}R^jp_*\left(\ttor_X^{i}(\mathcal{F},\mathcal{G})\right)
} holds for any coherent sheaves $\mathcal{F}$ and $\mathcal{G}$ of $\mathcal O_X$-modules, where $X$ is a formal scheme over $\CO F$ with structure map $\map pX{\Spf \CO F}$. Note that $\DEF n$ and $\GF^{2n}$ are affine formal schemes. To prove \eqref{fafa} and \eqref{meme}, it suffices to show that
\begin{equation}\label{chi}\comments{chi}\mathrm{Tor}_{\DEF n}^i(\cycn1n,\cycn2n)=0,\end{equation} 
\begin{equation}\label{shi}\comments{shi}\mathrm{Tor}_{\GF^{2h}}^i(\cycn1\infty,\cycn2\infty)=0\end{equation} 
	for any $i>0$, and that higher direct images vanish for $\DEF n\longrightarrow\Spf \CO F$ and $\G\longrightarrow\Spf\CO F$. However, since the vanishing of higher direct images follows directly from the fact that $\G$ and $\DEF n$ are affine formal schemes, the results \eqref{fafa} and \eqref{meme} follow directly from \eqref{chi} and \eqref{shi}. We prove \eqref{chi} and \eqref{shi} in the following lemmas.	
\begin{lem}[Acyclicity Lemma]\label{acyc}\cite[\href{http://stacks.math.columbia.edu/tag/00N0}{Tag 00N0}]{stacks-project}
\newcommand{\depth}{\mathrm{depth}}
Let $A$ be a Noetherian local ring, and let 
$$M^\bullet=0\rightarrow M^h\rightarrow\cdots \rightarrow M^0$$
be a complex of $A$-modules such that $\depth_A(M^i)\geq i$. If $\depth_A(H^{i}(M^\bullet))=0$ for any $i$, then $M^\bullet$ is exact.
\end{lem}
\[proof]{See Stacks Project \cite[\href{http://stacks.math.columbia.edu/tag/00N0}{Tag 00N0}]{stacks-project}.}
The following lemma gives some condition for the vanishing of higher torsion groups.
\begin{lem}\label{LEM2}\comments{LEM2}
\newcommand{\depth}{\mathrm{depth}}
\cite[\href{http://stacks.math.columbia.edu/tag/0B01}{Tag 0B01}]{stacks-project}
Suppose $A\rightarrow B_1$ and $A\rightarrow B_2$ are homomorphisms of regular local rings such that
\begin{enumerate}
\item $\dim(A)=2h$, $\dim(B_1)=\dim(B_2)=h$;
\item $\depth_A(B_1)=\depth_A(B_2)=h$;
\item $\length_A(B_1\otimes_AB_2)<\infty$.
\end{enumerate}
Then the higher torsion groups vanish, that is,
$$\mathrm{Tor}_A^i(B_1,B_2)=0$$
for any $i>0$.
\end{lem}
\begin{proof}\footnote{This proof is written according to \cite[\href{http://stacks.math.columbia.edu/tag/0B01}{Tag 0B01}]{stacks-project}.}
\newcommand{\depth}{\mathrm{depth}}
Since $A$ is a regular local ring, we have $\depth_A(A)=2h$.  
Therefore, there is a finite free $A$-module resolution $F_\bullet\rightarrow B_1$ of length
$$\text{depth}_A(A)-\text{depth}_A(B_1)=2h-h=h.$$
Note that $\mathrm{Tor}_A^i(B_1,B_2)$ is the i'th cohomology group of the complex $F_\bullet\otimes_A B_2\rightarrow B_1\otimes_A B_2$. Using Lemma \ref{acyc}, it suffices to prove the following claims:
\[itemize]{
\item $\depth_A(\mathrm{Tor}_A^i(B_1,B_2))=0$,
\item $\depth_A(F_i\otimes_A B_2)=\depth_A(B_2)=h\geq i$.
	}
Note that each cohomology group of the complex $F_\bullet\otimes_A B_2\rightarrow B_1\otimes_A B_2$ is a finite module over the Artinian ring $B_1\otimes_A B_2$. Hence $\depth_A(\mathrm{Tor}_A^i(B_1,B_2))=0$, which verifies the first claim. As each $F_i$ is a free $A$-module, we have $\depth_A(F_i\otimes_A B_2)=\depth_A(B_2)=h\geq i$, which verifies the second claim and therefore completes the proof.
\end{proof}
In fact, we are in a special situation where the depth of our ring is the same as its dimension, which we will describe in the following lemma.
\begin{lem}\label{LEM1}\comments{LEM1}
Suppose $A$,$B$ are regular local rings with residue field $\CFF q$. Let $\map{f}{\Spf B}{\Spf A}$ be a map such that $f^{-1}(\Spec\CFF q)$ is an Artinian scheme. We have
$$\mathrm{depth}_A(B)=\text{dim}(B).$$
\end{lem}
\begin{proof}
Let $J = \{f_1,f_2,\cdots f_n\}\subset\mathfrak m_A$ be a subset that generates the maximal ideal $\mathfrak m_A$ of $A$.
The scheme $f^{-1}(\Spec\CFF q)$ is the spectrum of the ring
$$
	B/(f_1,f_2,\cdots,f_n) = B/\mathfrak m_AB.
$$
Consider the subset $J':=\{f_i\}_{i\in \mathcal J}\subset J$, where
$$
	\mathcal J = \left\{i: \dim(B/(f_1,\cdots,f_{i-1}))>\dim(B/(f_1,\cdots,f_{i-1},f_i))\right\}.
$$
	Note that $\mathrm{depth}_A(B)$ is the number of elements in a regular sequence of $B$. It suffices to show that $J'$ is a regular sequence of length $\dim(B)$. Indeed,
since $f^{-1}(\Spec\CFF q)$ is Artinian, we have $\dim(B/(f_1,\cdots,f_{n}))=0$. Since $\dim(B/(f_1,\cdots,f_{i-1}))-\dim(B/(f_1,\cdots,f_{i-1},f_i))$ equals either $0$ or $1$, we have
$$
	\mathrm{depth}_A(B)=\#\mathcal J = \dim(B),
$$
which completes the proof.
\end{proof}
Now we are ready to complete the proof of Step 1.
\begin{proof}[Proof of Step 1: \eqref{chi} and \eqref{shi}]We prove these two identities together. For any 
	$$m\geq\max\left\{\m(\blin_1),\m(\blin_2)\right\},$$ 
we put
$$
A=\OO{\DEF n},\qquad B_1=\bpc m{\biso_1}{\blin_1}_*\OO{\DEf{n+m}},\qquad B_2=\bpc m{\biso_2}{\blin_2}_*\OO{\DEf{n+m}}
$$
	for the proof of \eqref{chi}, and put
$$A=\OO{\GF^{2h}},\qquad B_1=\gpc m{\biso_1}{\blin_1}_*\OO{\GK^{h}},\qquad B_2=\gpc m{\biso_2}{\blin_2}_*\OO{\GK^{h}}$$
for the proof of \eqref{shi}.
	Note that Proposition \ref{haolemma} implies that
	$$\eta^{\infty+m}_\infty(\biso_i,\blin_i)^{-1}\left(\Spec \CFF q\right)=\GK^h[\pi^m\biso_i\blin_i]\cong \eta^{n+m}_n(\biso_i,\blin_i)^{-1}\left(\Spec \CFF q\right)$$
for $i=1,2$.
Since $\J$ is an Artinian scheme,  we have  $\mathrm{depth}_A(B_1)=\dim(B_1)$ and $\mathrm{depth}_A(B_2)=\dim(B_2)$ by Lemma \ref{LEM1}. By Lemma \ref{LEM2}, it then suffices to show that
	\[itemize]{
\item $\dim(A)=2h$, $\dim(B_1)=\dim(B_2)=h$,
\item $\length_A(B_1\otimes_AB_2)<\infty$,
}
which are clearly true. This completes the proof of Step 1.
\end{proof}
\subsection{Step 2: Multiplicities inside the thickening}
\begin{proof}[Proof of Step 2: ]In this step, we prove \eqref{kunjuan}, which states that the following elements
$$
\sssss_*\sssss^*\big(\cycn1n\otimes_{\DEF n}\cycn2n\big),\qquad \ssssss_*\ssssss^*\big(\cycn1\infty\otimes_{\G}\cycn2\infty\big)
$$
have the same length.
Abbreviate $\delta_{i,n}:=\cycn in$ and $\delta_{i,\infty}:=\cycn i\infty$ for $i=1,2$. Note that
$$
\length\left(\sssss_*\sssss^*\left(\delta_{1,n}\otimes_{\DEF n}\delta_{2,n}\right)\right) = \length\left(\sssss^*\delta_{1,n}\otimes_{\G[\pi^M]}\sssss^*\delta_{2,n}\right)
$$
and
$$
\length\left(\ssssss_*\ssssss^*\left(\delta_{1,\infty}\otimes_{\G}\delta_{2,\infty}\right)\right) = \length\left(\ssssss^*\delta_{1,\infty}\otimes_{\G[\pi^M]}\ssssss^*\delta_{2,\infty}\right).
$$
It suffices to prove for $i=1,2$,
$$
\sssss^*\cycn i{n} = \ssssss^*\cycn i{\infty},
$$
which follows directly from Proposition \ref{juda} when $n>\max(\m(\blin_1),\m(\blin_2))+M$. This completes the proof of Step 2.
\end{proof}
\subsection{Step 3: Actual Multiplicity}In this step, our goal is to show that
\begin{equation*}
\sssss_*\sssss^*\big(\cycn1n\otimes_{\DEF n}\cycn2n\big)=\cycn1n\otimes_{\DEF n}\cycn2n,
\end{equation*} 
\begin{equation*}
\ssssss_*\ssssss^*\big(\cycn1\infty\otimes_{\G}\cycn2\infty\big)=\cycn1\infty\otimes_{\G}\cycn2\infty.
\end{equation*}
\begin{lem}\label{RESTR}\comments{RESTR}
Let $A$, $B$ be Noetherian local rings and let $\mathfrak m_A$, $\mathfrak m_B$ be their maximal ideals respectively. Let $\mathcal F$ be a coherent sheaf of $\mathcal O_{\Spf A}$-module. Suppose there is a closed embedding
$$\map{s}{\Spf B/\mathfrak m_B^M}{\Spf A}$$
	such that $\length_A(s_*s^*\mathcal F)< M$. If 
either of the following conditions hold
\[itemize]{
\item
	$\dim(A)=\dim(B)$ and both $A$ and $B$ are complete regular local rings with the same residue field,
\item  $A=B$ and $s$ is a morphism over $\Spf A$,
} 
	then we have
	$$s_*s^*\mathcal F\cong\mathcal F.$$
\end{lem}
\begin{proof}
	Let $E=\mathcal F(\Spf A)$ be the $A$-module of global sections of $\mathcal F$, and we denote the induced map of $s$ by 
	$$
	\map{\mathcal S}A{B/\mathfrak m_B^M}.
	$$
Since $s$ is a closed embedding, the homomorphism $\mathcal S$ is surjective. Since $\Spf A$ and $\Spf B$ are affine formal schemes, it suffices to prove that
	$$
	E\otimes_A B/\mathfrak m_B^M = E,
	$$
	which is equivalent to
	\[equation]{\lb{nenenene}
	\mathcal S^{-1}(\mathfrak m_B^M) E = 0.
}
	The key step of the proof is showing that $\mathcal S^{-1}(\mathfrak m_B^M) = \mathfrak m_A^M$, which is clearly true when $A=B$ with $s$ a morphism over $\Spf A$. Therefore, it remains to prove for the other case of complete regular local rings of the same dimension.
Since $A$ and $B$ are complete local rings,
we have $\mathcal S(\mathfrak m_A)\subset \mathfrak m_B$ and $\mathcal S^{-1}(\mathfrak m_B)\subset\mathfrak m_A$. Therefore the surjective ring homomorphism $\mathcal S$ gives rise to a surjective graded ring homomorphism
$$
	\map{\widetilde{\mathcal S}}{\bigoplus_{n=0}^{M-1}{\mathfrak m_A^{n}/\mathfrak m_A^{n+1}}}{\bigoplus_{n=0}^{M-1}{\mathfrak m_B^{n}/\mathfrak m_B^{n+1}}},
$$
	where we denote $\mathfrak m_A^{0}:=A$ and $\mathfrak m_B^{0}:=B$. 
	Letting $d:=\dim(A)=\dim(B)$ and $k:=A/\mathfrak m_A\cong B/\mathfrak m_B$, we have
	$$
	\dim_{k}\left(\mathfrak m_A^{n}/\mathfrak m_A^{n+1}\right) = \frac{(n+d-2)!}{(d-1)!(n-1)!} =  \dim_{k}\left(\mathfrak m_B^{n}/\mathfrak m_B^{n+1}\right),
	$$
	which implies that $\widetilde{\mathcal S}$ is an isomorphism
	. Therefore in both cases, we obtain the key identity
	\[equation]{\lb{ningmeng}
	\mathcal S^{-1}(\mathfrak m_B^M) = \mathfrak m_A^M.
	}
	Since $\length_{A}(E)<M$, the descending chain
$$
E\supset \mathfrak m_A E\supset\cdots\supset \mathfrak m_A^M E
$$
must have the property that $\mathfrak m_A^iE=\mathfrak m_A^{i+1}E$ for some $i\leq M-1$, which implies $\mathfrak m_A^{i+1}E=0$ by Nakayama lemma. Hence $\mathfrak m_A^ME=0$, which together with \eqref{ningmeng} prove \eqref{nenenene}, and therefore complete the proof of this proposition.\end{proof}
Next, we finish the proof of Step 3.
\begin{proof}[Proof of Step 3]
Finally, we prove \eqref{eshu} and \eqref{ishu}. Let $\OO{\GF^{2h}}$ be the ring of global sections of the structure sheaf of $\GF^{2h}$, and let $\mathfrak m_{\GF^{2h}}$ be the corresponding maximal ideal. Since 
$$\deg(\xymatrix{\GF\ar[rr]^{\pi^M}&&\GF})=q^{2hM},$$
one sees that $\Spec \OO{\GF^{2h}}/\mathfrak m_{\GK^{2h}}^{q^{2hM}}$ is a closed subscheme of
	$$\ker\left( \xymatrix{\GF^{2h}\ar[rr]^{\pi^M}&&\GF^{2h}}\right)=\GF^{2h}[\pi^M],$$
and we denote the corresponding closed embedding by
$$
\map{j_M}{\Spec \OO{\GF^{2h}}/\mathfrak m_{\GK^{2h}}^{q^{2hM}}}{\JJJJJ}.
$$
Consider the following closed embeddings 
\begin{equation}\label{tougaoa}\comments{tougaoa}
\xymatrix{\Spf \OO{\GF^{2h}}/\mathfrak{m}_{\GF^{2h}}^{q^{2hM}}\ar[rr]^{\qquad j_M}&&\JJJJJ\ar[r]^\ssssss&\G},
\end{equation}
\begin{equation}\label{tougaob}\comments{tougaob}
\xymatrix{\Spf \OO{\GF^{2h}}/\mathfrak{m}_{\GF^{2h}}^{q^{2hM}}\ar[rr]^{\qquad j_M}&&\JJJJJ\ar[r]^\sssss&\DEF n}.
\end{equation}
	Let $m$ be an integer such that $m\geq \max\{\m(\blin_1),\m(\blin_2)\}$, and denote by $\DEf{1,m+n}$, $\DEf{2,m+n}$ the corresponding $\pi^{m+n}$-level Lubin Tate spaces associated to the pair $(\biso_1,\blin_1)$ and $(\biso_2,\blin_2)$. Put
$$
\mathcal F_n := \bpc m{\biso_1}{\blin_1}_*\OO{\DEf {1,n+m}}\otimes_{\DEF n}\bpc m{\biso_2}{\blin_2}_*\OO{\DEf {2,n+m}}
$$
and
$$
	\mathcal F_\infty := \gpc m{\biso_1}{\blin_1}_*\OO{\GG_{K_1}^{h}}\otimes_{\GF^{2h}}\gpc m{\biso_2}{\blin_2}_*\OO{\GG_{K_2}^{h}}.
$$
	Note that $n\geq 1$, and therefore $\deg\left(\GG_{K_i}^h\xrightarrow{\pi^m}\GG_{K_i}^h\right) = q^{2hm}= \dpa {i,n}m$ for $i=1,2$, which implies that
$$
	[\mathcal F_\infty] = q^{4hm}\cycn1\infty\otimes_{\G}\cycn2\infty,\qquad
	[\mathcal F_n] = q^{4hm}\cycn1n\otimes_{\DEF n}\cycn2n.
$$
Therefore \eqref{kunjuan} implies $\length(\sssss^*\mathcal F_n)=\length(\ssssss^*\mathcal F_\infty)$.
It then suffices to prove $\sssss_*\sssss^*\mathcal F_n=\mathcal F_n$ and $\ssssss_*\ssssss^*\mathcal F_\infty=\mathcal F_\infty$. Let $M$ be any integer such that
\begin{equation}\label{howbigm}\comments{howbigm}
\begin{split}
M&>\frac1{2h}\log_q\length(\mathcal F_\infty).\\
\end{split}
\end{equation}
Note that we have a canonical isomorphism
$$
	\OO{\GF^{2h}}/\mathfrak{m}_{\GF^{2h}}^{q^{2hM}}\cong \OO{\GF^{2h}[\pi^M]}/\mathfrak{m}_{\GF^{2h}[\pi^M]}^{q^{2hM}},
$$
	and that $\length(\sssss^*\mathcal F_n) = \length(\ssssss^*\mathcal F_\infty)<q^{2hM}$.  Applying Lemma \ref{RESTR} to the map $j_M:\Spf \OO{\GF^{2h}[\pi^M]}/\mathfrak{m}_{\GF^{2h}[\pi^M]}^{q^{2hM}}\rightarrow \JJJJJ$, we have
$$
	j_{M*}j_M^*\left(\ssssss^*\mathcal F_\infty\right)=\ssssss^*\mathcal F_\infty,\qquad
	j_{M*}j_M^*\left(\sssss^*\mathcal F_n\right)=\sssss^*\mathcal F_n.
$$
Applying Lemma \ref{RESTR} again to the map $\ssssss\circ j_M:\Spf \OO{\GF^{2h}}/\mathfrak{m}_{\GF^{2h}}^{q^{2hM}}\rightarrow \G$, and to the map $\sssss\circ j_M:\Spf \OO{\GF^{2h}}/\mathfrak{m}_{\GF^{2h}}^{q^{2hM}}\rightarrow \DEF n$, we have
$$
(\ssssss\circ j_M)_*(\ssssss\circ j_M)^*\mathcal F_\infty = \mathcal F_\infty,\qquad (\sssss\circ j_M)_*(\sssss\circ j_M)^*\mathcal F_n = \mathcal F_n.
$$
Therefore $\ssssss_*\ssssss^*\mathcal F_\infty = \mathcal F_\infty$ and $\sssss_*\sssss^*\mathcal F_n = \mathcal F_n$  as desired.\end{proof}

\begin{proof}[Proof of Theorem \ref{THM:ICT}]
The theorem follows from combining all steps. Note that Step 2 works for $N\geq M+\max(\m(\blin_1),\m(\blin_2))+1$, and that Step 3 works for such $M$ and $m$ that satisfies \eqref{howbigm} and $m\geq\max(\m(\blin_1),\m(\blin_2))$. Therefore, any choice of $N$ such that 
\begin{equation}\label{howbigN}
N \geq \frac1{2h}\log_q\length(\cycn1\infty\otimes_{\GF^{2h}}\cycn2\infty) + 3\max(\m(\blin_1),\m(\blin_2)) + 1
\end{equation}
completes the proof of Theorem \ref{THM:ICT}.\end{proof}

\section{Computation of intersection numbers on the infinite level}\label{zhaojile}\comments{zhaojile}
In this section, we give an explicit formula for
\begin{equation}\label{mubiao}\comments{mubiao}
\length(\cycn1\infty\otimes_{\G}\cycn2\infty).
\end{equation} 
We keep the same notation as in Section \ref{mi}. We recall some elementary linear algebra facts for formal modules, which will be used to obtain our main results in Theorem \ref{yanchanghui}, which applies generally whether $K_1\cong K_2$ or $K_1\not\cong K_2$. Then we simplify the result to Theorem \ref{xiaomubiao} for the case of $K_1\cong K_2$. Let $K$ be either $K_1$ or $K_2$.

\subsection{Some linear algebra for formal modules}\label{notata}\comments{notata}
We briefly review some elementary linear algebra concepts that will be used in our calculation.
\subsubsection{Formal matrix multiplication}
For $R=F,K,\OO F$ and $\OO K$, we fix canonical isomorphisms 
$$R^n\cong \matt n1R\quad\text{and}\quad R^{n\vee}\cong \matt1nR.$$
The paring $R^{n\vee}\otimes_R R^n\longrightarrow R$ is given by the matrix multiplication in an obvious way. 
\[defi]{[Formal vector]\lb{fove}
When $R=F$ or $K$, for any formal $\OO R$-module $\GG$ over $\CFF q$, and for any $A\in\CCCC$, we represent each element $f\in\GG^n(A)$ as the $n\times 1$ matrix 
$$
f=\[pmatrix]{{f_1}\\{\vdots}\\{f_n}},
$$
such that $f_i=f(e_i)\in\GG(A)$ for any $1\leq i\leq n$, where $e_i\in\mathcal O_R^\vee$ is the $1\times n$ matrix $(0\;\cdots\;1\;\cdots\;0)$ with i'th entry 1 and 0 elsewhere.
}

\[defi]{[Formal matrix multiplication]\lb{fm}
When $R=F$ or $K$, for any formal $\OO R$-module $\GG$ over $\CFF q$, and any matrix $M=(m_{ij})_{1\leq i\leq a}^{1\leq j\leq b}\in\matt ab{\End(\GG)}$, we define the formal left multiplication map $\xymatrix{\GG^b\ar[r]^{M}&\GG^a}$, functorially in $A\in\CCCC$, by
$$
\maps{M}{\GG^b(A)}{\GG^a(A)}{\[pmatrix]{{f_1}\\{\vdots}\\{f_b}}}{\[pmatrix]{{[m_{11}]_\GG(f_1)[+]_\GG\cdots[+]_\GG [m_{1b}]_\GG(f_b)}\\{\vdots}\\{[m_{a1}]_\GG(f_1)[+]_\GG\cdots[+]_\GG [m_{ab}]_\GG(f_b)}}_{\cdot}}
$$
}

\[rem]{Under the identification $\GG^n(A)\cong\Hom_{\OO R}(\OO R^{n\vee},\GG(A))$ for any $A\in\CCCC$, one sees that for any $f\in\GG^n(A)$, after expressing $f$ as an $n\times 1$ matrix as in Definition \ref{fove}, the map $f:\OO R^{n\vee}\rightarrow \GG(A)$ is obtained from the formal matrix multiplication $v\mapsto v\cdot f$ for any $v\in\OO R^{n\vee}$.
	}
Let $\map{\biso'}{\GG_1}{\GG_2}$ be a quasi-isogeny of two formal groups. By abuse of notation, we also use the same symbol $\map{\biso'}{\GG_1^n}{\GG_2^n}$ for the map obtained by applying $\biso'$ componentwisely on $\GG_1^n$ for any $n\geq0$.

In this section, we work with two quadratic extensions $K_1/F$ and $K_2/F$. Let $\GG_{K_1}=\GG_F$ (resp. $\GG_{K_2}=\GF$) be the formal $\OO F$-module with an extra $\OO{K_1}$- (resp. $\OO{K_2}$) -action induced from a fixed embedding $K_1$ (resp. $K_2$) $\subset \End(\GF)\otimes_{\OO F}F$. We fix these embeddings throughout the section. Since $\GG_{K_1}=\GG_{K_2}=\GF$ as formal $\OO F$-modules, there are canonical isomorphisms  
$$\Hom_{\OO F}(\GG_{K_1},\GF)\cong\End_{\OO F}(\GF)\cong\Hom_{\OO F}(\GG_{K_2},\GF),$$
which identifies the quasi-isogenies $\map{\biso_1}{\GG_{K_1}^{2h}}{\GF^{2h}}$ and $\map{\biso_2}{\GG_{K_2}^{2h}}{\GF^{2h}}$ with scalar matrices $\biso_1I_{2h},\biso_2I_{2h}\in\matt{2h}{2h}{D_F}$ through $D_F=\End_{\OO F}(\GF)\otimes_{\OO F}F$.

Let $a,b$ be two non-negative integers. For any matrix $M\in\matt ab{K_i}$, we denote by $M_{\GG_{K_i}}$ its corresponding image in $\matt ab{D_F}$ through the embedding $K_i\subset D_F$ for any $i=1,2$. Note that even if $K:=K_1\cong K_2$, the embeddings $K_1\subset D_F$ and $K_2\subset D_F$ may have different images, where one might have $M_{\GG_{K_1}}\neq M_{\GG_{K_2}}$  for some $M\in\matt abK$. However, if $M\in\matt abF$, then we view $M$ as a matrix in $\matt ab{D_F}$ naturally through the embedding of the center $F\rightarrow D_F$. 

For ease of exposition, if a construction is made for both $\GG_{K_1}$ and $\GG_{K_2}$, then we will omit subscripts.
The map $\map\gpp{\GK^h}{\GF^{2h}}$ in Definition \ref{aho} can be represented as a matrix in $\matt{2h}h{D_F}$, which we describe now explicitly. 
\[defi]{\lb{matu}
Let $M_\blin\in \matt{2h}hK$ be such that for any $1\times 2h$ matrix $(a_1\;\cdots\;a_{2h})\in F^{2h\vee}$, we have
$$
\blinn (a_1\;\cdots\;a_{2h}) = (a_1\;\cdots\;a_{2h})M_\blin\in K^{h\vee}.
$$
}
\[prop]{For any $(\biso,\blin)\in\equi KFh$, and any integer $m\geq \m(\blin)$, the map $\map{\gpp}{\GK^{h}}{\GF^{2h}}$ can be decomposed as described by the following commutative diagram  
$$
\xymatrix{
	\GK^h\ar[rrd]_{\gpp}\ar[rrrr]^{(M_\blin)_{\GK}}&&&&\GK^{2h}\ar[lld]^{\pi^m\biso}\\
	&&\GF^{2h}.
	}
$$
}
\[proof]{For any $A\in\CCCC$, and any
$$
f=\[pmatrix]{f_1\\\vdots\\f_h}\in \GK^h(A),
$$
let $f'=\gpp f=\pi^m\biso\circ f\circ\blinn$. For any $a=(a_1\;\cdots\;a_{2h})\in\mathcal O_F^{2h\vee}$, we have
\[equation*]{\[split]{
	f'(a) &= \pi^m\biso\circ f\circ \blinn (a)\\
	&=\pi^m\biso\circ f\left((a_1\;\cdots\;a_{2h})M_\blin\right)\\
	&=\pi^m\biso\left((a_1\;\cdots\;a_{2h})_{\GK}\cdot(M_\blin)_{\GK}\cdot f\right).
}}
Since $\map{\pi^m\biso}{\GK}{\GF}$ is an isogeny of formal $\OO F$-modules, we have $\pi^m\biso\circ[x]_{\GK}=[x]_{\GF}\circ\pi^m\biso$ for any $x\in\OO F$. Therefore the preceding results imply
$$
f'(a)=(a_1\;\cdots\;a_{2h})\cdot\pi^m\biso\left((M_\blin)_{\GK}\cdot f\right).
$$
Hence
$$
f'=\pi^m\biso((M_\blin)_{\GK}\cdot f),
$$
which completes the proof.}
\[rem]{Using a similar method, for any $(\ACj,\ACg)\in\equi FF{2h}$, and any integer $m\geq \m(\ACg)$, the map $\map{\gppbj{}{m}}{\GF^{2h}}{\GF^{2h}}$ can be decomposed as described by the following commutative diagram
$$
\xymatrix{
	\GF^{2h}\ar[rrd]_{\gppbj{}m}\ar[rrrr]^{g}&&&&\GF^{2h}\ar[lld]^{\pi^m\ACj}\\
	&&\GF^{2h}.
	}
$$
}
Note that $g M_{\blin} = M_{g\blin}$ for any $g\in\GL_{2h}(F)$ and any $\blinn\in\Hom_F(F^{2h\vee},K^{h\vee})$.

\subsubsection{Reduced norm and degree}
Note that $\dim_F(D_F)=4h^2$. For any $n\geq 0$, the algebra $\matt nn{D_F}$ is a $(2h)^2\times n^2$-dimensional $F$-central simple algebra, for which $\matt nn{D_F}\otimes_F F^{\mathrm{alg}}$ is isomorphic to $\matt{2hn}{2hn}{F^\mathrm{alg}}$ over the algebraic closure $F^{\mathrm{alg}}$ of $F$. The reduced norm of an element $g\in\matt nn{D_F}$ is simply its determinant as an element of $\matt{2hn}{2hn}{F^\mathrm{alg}}$. In particular, in cases of $n=2h$, $n=h$, and $n=1$, we denote by 
$$\mathrm{NRD}(g),\qquad \mathrm{Nrd}(g'), \qquad\text{and}\qquad\mathrm{nrd}(g'')$$ 
the reduced norms of elements $g\in\matt {2h}{2h}{D_F}$, $g'\in\matt hh{D_F}$ and $g''\in D_F$ respectively. We will show that the absolute value of the reduced norm for an invertible element is related to the degree of the quasi-isogeny induced by it. From now on, we denote $\OO{D_F}:=\End_{\OO F}(\GF)$, which is the maximal order of $D_F$. We start with the following basic case.
\[lem]{\lb{shabusha}
	For any $\gamma\in \OO{D_F}\cap D_F^\times$, we have
	$$
	\deg\left(\xymatrix{\GF\ar[r]^{\gamma}&\GF}\right)=\left|\mathrm{nrd}(\gamma)\right|_F^{-1}.
	$$
	}
\[proof]{ Let $\GG_L$ be the formal $\mathcal O_L$-module over $\CFF q$ obtained from $\GF$ by choosing a maximal subfield $L\subset D_F$ such that $\gamma\in L$, where the tangent space $\mathrm{Lie}(\GG_L)=\mathrm{Lie}(\GG_F)$ is equipped with an $\mathcal O_L$-algebra structure via the action of $\mathcal O_L$ on $\mathrm{\GG_F}$. Note that the height of $\GG_L$ is $1$. Let $\breve L$ be the completion of the maximal unramified extension of $L$. By Lubin and Tate \cite{lubin1966formal}, there is a deformation $\widetilde{\GG_L}$ of $\GG_L$ over $\CO L$ with its Tate-module isomorphic to $\OO L$. It suffices to prove this lemma for $\widetilde{\GG_L}$. Note that the degree of the induced map 
$$\deg\left(\xymatrix{\widetilde{\GG_L}\ar[r]^{\gamma}&\widetilde{\GG_L}}\right)$$ 
equals the number of elements in the cokernel of the induced map for their Tate modules 
$$\#\coker\left(\xymatrix{\OO L\ar[r]^{\gamma}&\OO L}\right),$$ 
which equals $|\gamma|_L^{-1} = |\mathrm{nrd}(\gamma)|_F^{-1}$. This completes the proof of the Lemma.
	}
More generally, we have the following.
\begin{lem}\label{isonimei}\comments{isonimei}
	For any $g\in\matt hh{\OO{D_F}}$, we have
	$$\deg\left(\xymatrix{ {\GF^h}\ar[r]^{g}&{\GF^h}}\right)=\|\Nrd(g)\|_F^{-1}.$$
\end{lem}
\begin{proof}
	By Cartan decomposition of matrix algebras over division algebras, there exists $u_1,u_2\in\GL_h(\OO{D_F})$ and a diagonal matrix $t=\mathrm{diag}(t_1,t_2,\cdots,t_h)$ such that
	$$g=u_1tu_2,$$ 
	which decomposes the map $\map {g}{\GF^h}{\GF^h}$ into the following homomorphisms
$$
	\xymatrix{
		\GF^h\ar[rr]^{u_2}&&\GF^h\ar[rr]^{t}&&\GF^h\ar[rr]^{u_1}&&\GF^h.
		}
$$
This implies that the degree of $g:\GF^h\rightarrow\GF^h$ equals the degree of $t:\GF^h\rightarrow\GF^h$ since $u_1:\GF^h\rightarrow\GF^h$ and  $u_2:\GF^h\rightarrow\GF^h$ are isomorphisms. Finally, note that
	$$
	\deg\left(\xymatrix{\GF^h\ar[r]^{t}&\GF^h}\right)=\prod_{i=1}^h\deg\left(\xymatrix{\GF\ar[r]^{[t_i]_{\GF}}&\GF}\right)=\prod_{i=1}^h\|\nrd(t_i)\|_F^{-1}=\|\Nrd(t)\|_F^{-1}=\|\Nrd(g)\|_F^{-1},
	$$
which completes the proof.
\end{proof}

\subsubsection{Associated quasi-isogeny}
Let $\sigma$ be the non-trivial element in $\Gal(K/F)$, and let $M_{\blin}^{\sigma}$ be the matrix obtained by applying $\sigma$ entrywise to $M_{\blin}$. 
The following construction is the key to our calculation.
\[defi]{\lb{quaiso}
We associate the pair $(\biso,\blin)$ to the following matrix 
\[equation]{\label{lashi}\comments{lashi}\Diso:=\biso \cdot\[pmatrix]{{M_\blin}&{M_\blin^\sigma}}_{\GK}\in \matt{2h}{2h}{D_F},
}
which represents the composition of the following quasi-homomorphisms
$$
\xymatrix{
	\GK^{2h}
	\ar[rrrrr]^{\[pmatrix]{{M_\blin}&{M_\blin^\sigma}}_{\GK}}&&&&&\GK^{2h}\ar[rr]^{\biso}&&\GF^{2h}.
	}
$$
}
Let $|\mathrm{Disc}_{K/F}|_F$ be the norm of the relative discriminant of $K/F$. We show that $\Diso$ is in fact a quasi-isogeny.
\begin{lem}\lb{disco}
For any $(\biso,\blin)\in\equi KFh$, we have
$$
|\mathrm{NRD}(\Delta(\biso,\blin))|_F=|\mathrm{Disc}_{K/F}|_F^{h^2}
\neq 0.
$$
In particular, $\Diso\in\GL_{2h}(D_F)$.
\end{lem}
\begin{proof}
Let $\mu\in\OO K$ be a generator such that $\OO K=\OO F[\mu]$. Let
$$
M_0=\cc{I_h}{\mu I_h}.
$$
	Since right multiplication by $M_0$ induces an isomorphism of $\mathcal O_F$-modules $\OO F^{2h\vee}\rightarrow \OO K^{h\vee}$, there exsits $g\in \GL_{2h}(F)$ such that $\blan = gM_0$. Then we have
	$$\Vol(\blin^{\vee}(\OO F^{2h\vee}))=\Vol(\OO F^{2h\vee}\blan) = \Vol(\OO F^{2h\vee}gM_0)=\Vol(\OO F^{2h\vee}g).$$
Hence	
$$
|\det(g)|_F=\frac{\Vol\left(\OO F^{2h\vee}g\right)}{\Vol\left(\OO F^{2h\vee}\right)}=\frac{\Vol\left(\blinn(\mathcal O_F^{2h\vee})\right)}{\Vol\left(\OO F^{2h\vee}\right)}=q^{\mathrm{Height}(\blin)},
$$
which implies that $|\mathrm{NRD}(g)|_F=|\det(g)^{2h}|_F=q^{2h\mathrm{Height}(\blin)}$.

On the other hand, note that the quasi-isogeny $\biso:\GK\rightarrow\GF$ is of degree $|\nrd(\biso)|_F^{-1}$. By the definition of $\mathrm{Height}(\biso)$, we have
	$$
	|\nrd(\biso)|_F=q^{-\mathrm{Height}(\biso)},
	$$
	which implies $|\mathrm{NRD}(\biso I_{2h})|_F=|\nrd(\biso)^{2h}|_F=q^{-2h\mathrm{Height}(\biso)}$.
Finally,
$$
|\mathrm{NRD}\left((M_0\;M_0^\sigma)_{\GK}\right)|_F = \left|\mathrm{Nm}_{K/F}
	\left(
	\mathrm{det}_K
	\mm{I_h}{I_h}{\mu I_h}{\mu^\sigma I_h}
	\right)
	\right|_F^{h}
=|\mathrm{Disc}_{K/F}|_F^{h^2}.
$$
	Note that $\Delta(\biso,\blin)=\biso\cdot g\cdot (M_0\;M_0^{\sigma})_{\GK}$, and $\mathrm{Height}(\blin)=\mathrm{Height}(\biso)$. Combining the preceding equations, we complete the proof of this lemma.
\end{proof}

\subsection{Decomposition of $\cycy\infty$}
We fix the following notation throughout the proof of Lemma \ref{kuangxie} and Theorem \ref{yanchanghui}. Abbreviate $\left[\OO X\right]$ to $\left[X\right]$ for any formal scheme $X$. For any $\blinn\in\Hom_{F}(F^{2h\vee},K^{h\vee})$ such that $\blinn$ is an isomorphism, consider the Iwasawa decomposition of the matrix $\left(\blan\quad M_\blin^\sigma\right)\in\GL_{2h}(K)$ 	
$$
\left(M_\blin\quad M_\blin^\sigma\right) = \Gamma_\blin\cdot\mm{P_\blin}{*}{}{Q_\blin}\in\GL_{2h}(K),
$$
where $Q_\blin,P_\blin\in\GL_{2h}(K)$ and $\Gamma_\blin\in\GL_{2h}(\mathcal O_K)$. We also fix an integer $m$ such that $\pi^m\biso M_\blin\in\matt {2h}h{\mathcal O_{D_F}}$. Therefore $\pi^m\biso P_\blin, \pi^m\biso Q_\blin\in\matt hh{\mathcal O_{D_F}}\cap\GL_h(D_F)$, which implies that they are corresponding to actual isogenies. Furthermore, let $i_{\biso,\blin}$ and $p_{\biso,\blin}$ be the maps defined by
$$
\iii:\xymatrix{\GF^h\ar[rrrr]^{\biso\cdot\gaga\cdot\bisoo\cc{I_h}{0_h}}&&&&\GF^{2h}},
$$
and
$$
\ppp:\xymatrix{\GF^{2h}\ar[rrrr]^{\rr{0_h}{I_h}\biso\cdot\ggga\cdot\bisoo}&&&&\GF^{h}}.
$$
Note that $\iii$ is a closed embedding, and that these maps depend on the choice of $\Gamma_\blin$. The following lemma gives a decomposition of the cycle $\cycy\infty$.
\begin{lem}\label{kuangxie}\comments{kuangxie}
We have
	\[equation]{\lb{zia}
		\cycy\infty=\left|\Nrd\left(\biso\cdot(P_\blin)_{\GK}\right)\right|_F^{-1}\left[\mathrm{Im}\left(\xymatrix{\GF^h\ar[rr]^{\iii}&&\GF^{2h}}\right)\right],}
	\[equation]{\lb{zib}
	\cycy\infty=\left|\Nrd\left(\biso\cdot(P_\blin)_{\GK}\right)\right|_F^{-1}\left[\mathrm{Ker}\left(\xymatrix{\GF^{2h}\ar[rr]^{\ppp}&&\GF^{h}}\right)\right].
}
\end{lem}

\begin{proof}
Since
$$
\pi^m\biso M_\blin = \pi^m\biso\left(M_\blin\quad M_\blin^\sigma\right)\cc{I_h}{0} = \pi^{m}\biso\Gamma_\blin\mm{P_\blin}{*}{}{Q_\blin}\cc{I_h}{0} = \biso\Gamma_\blin\biso^{-1}\cdot\cc{I_h}{0} \cdot\pi^{m}\biso P_\blin,
$$
	the map $\gpp:\xymatrix{\GK^h\ar[rrr]^{\pi^m\cdot (M_\blin)_{\GK}}&&&\GF^{2h}}$ can be decomposed into 
	\[equation]{\lb{tingche}\xymatrix{\GKH\ar[rr]^{\pi^{m}\biso\cdot\pipa}&&\GF^h\ar[rr]^{\cc{I_h}0}&&\G\ar[rr]^{\biso\cdot{\gaga}\cdot\biso^{-1}}&&\G}. 
	}
Note that $\map{\pi^{m}\biso\cdot\pipa}{\GK^h}{\GF^h}$ is an isogeny, that $\map{\cc{I_h}0}{\GF^h}\G$ is a closed embedding, and that $\map{\biso\cdot{\gaga}\cdot\biso^{-1}}\G\G$ is an isomorphism. This implies that
	\[equation]{\lb{fabiao}
	\gpp_*[\GK^h] = \deg\left(\xymatrix{\GKH\ar[rr]^{\pi^m\cdot\biso\cdot( P_\blin)_{\GK}}&&\GKH}\right)\left[\mathrm{Im}\left(\xymatrix{\GF^h\ar[r]^{\iii}&\GF^{2h}}\right)\right].
	}
By Definition \ref{infinf}, we have
	\[equation*]{\lb{fabiaobiao}
	\gpp_*[\GK^h] = \deg\left(\xymatrix{\GK^h\ar[r]^{\pi^m}&\GK^h}\right)\cycy\infty.
	}
Comparing this identity with \eqref{fabiao}, we complete the proof of \eqref{zia}. Moreover, considering the following exact sequence 
\begin{equation}
\shexkk{\GF^h}{\GF^{2h}}{\GF^h}{\cc{I_h}0}{\rr0{I_h}},
\end{equation}
	and twisting it by the isomorphism $\map{\biso^{-1}\cdot\gaga\cdot\biso}{\GF^{2h}}{\GF^{2h}}$, we obtain an exact sequence
\begin{equation}
	\shexkk{\GF^h}{\GF^{2h}}{\GF^h}{\biso^{-1}\cdot(\Gamma_{\blin})_{\GK}\cdot\biso\cdot\cc{I_h}0}{\rr{0}{I_h}\cdot\biso^{-1}\cdot(\Gamma_{\blin}^{-1})_{\GK}\cdot\biso},
\end{equation}
which implies that
\begin{equation}\label{heheda}\mathrm{Im}\left(\xymatrix{\GF^h\ar[r]^{\iii}&\GF^{2h}}\right)=\mathrm{Ker}\left(\xymatrix{\GF^{2h}\ar[r]^{\ppp}&\GF^{h}}\right).\end{equation} 
	Comparing this identity with \eqref{zia}, we finish the proof of \eqref{zib}, and hence complete the proof of the lemma. 
\end{proof}

\subsection{Computation of the intersection number}
In this part, we give a formula for the intersection number $\chi(\cycn1\infty\litimes_{\GF^{2h}}\cycn2\infty)$. We give the most general version of the intersection formula by the following theorem.
\begin{thm}\label{yanchanghui}\comments{yanchanghui}
Suppose $(\biso_i,\blin_i)\in\equi{K_i}Fh$ for $i=1,2$. If the following expression 
	\[equation]{\lb{fofofo}
|\mathrm{Disc}_{K_1/F}|_F^{-h^2}\cdot\left|\mathrm{Nrd}\left(\[pmatrix]{{0_h}&{I_h}}\cdot\Delta(\biso_1,\blin_1)^{-1}\cdot\Delta(\biso_2,\blin_2)\cdot\cc{I_h}{0_h}\right)\right|_F^{-1}
	}
is finite, then it equals to the intersection number $\chi(\cycn1\infty\litimes_{\GF^{2h}}\cycn2\infty)$.
\end{thm}
\begin{proof}
By (\ref{meme}), we have 
$$
\chi(\cycn1\infty\litimes_{{\G}}\cycn2\infty) = \length(\cycn1\infty\otimes_{{\G}}\cycn2\infty)
$$
when the right hand side is finite, which, by Proposition \ref{kuangxie}, equals
	\begin{equation}\lb{bashibashi}
\|\Nrd(\biso_1\cdot\pipao\cdot\biso_2\cdot\pipat)\|_F^{-1}
	\length\left(\left[\mathrm{Im}
	\left(\xymatrix{\GF^{h}\ar[r]^{\iiit}&\GF^{2h}}\right)
	\times_{\G}
	\Ker
	\left(\xymatrix{\GF^{2h}\ar[r]^{\pppo}&\GF^{h}}\right)\right]\right).
\end{equation}
The map $\iiit:\GF^h\rightarrow\GF^{2h}$ induces a natrual map
	$$\Ker\left(\xymatrix{\GF^{h}\ar[rr]^{\pppo\circ\iiit}&&\GF^{h}}\right)\longrightarrow
\mathrm{Im}
	\left(\xymatrix{\GF^{h}\ar[r]^{\iiit}&\GF^{2h}}\right)
	\times_{\G}
	\Ker
	\left(\xymatrix{\GF^{2h}\ar[r]^{\pppo}&\GF^{h}}\right),
$$
which is an isomorphism since $\iiit$ is a closed embedding. Therefore 
	\[equation]{\lb{bishibishi}
\length\left(\left[\mathrm{Im}
	\left(\xymatrix{\GF^{h}\ar[r]^{\iiit}&\GF^{2h}}\right)
	\times_{\G}
	\Ker
	\left(\xymatrix{\GF^{2h}\ar[r]^{\pppo}&\GF^{h}}\right)\right]\right) = \left|\Nrd(\pppo\cdot\iiit)\right|_F^{-1}.
	}
Note that
$$
\pppo\cdot\iiit = 
	\rr{0_h}{I_h}\cdot\biso_1\cdot\gggao\cdot\bisoo_1\cdot \biso_2\cdot\gagat\cdot\bisoo_2 \cdot\cc{I_h}{0_h},
$$
$$
\biso_2\cdot\gagat\cdot\bisoo_2\cdot\cc{I_h}{0_h}=\Disot\cdot\cc{I_h}{0_h}\cdot\pppat\cdot\biso_2^{-1},$$
and that
$$\[pmatrix]{0_h&I_h} \cdot\biso_1\cdot\gggao\cdot\bisoo=\biso_1 \cdot \qaqao \cdot\[pmatrix]{0_h&I_h}\cdot\Disoo^{-1}.
$$
Therefore,
	\[equation]{\lb{bachai}
	\pppo\cdot\iiit = \biso_1\cdot\qaqao\cdot\[pmatrix]{0_h&I_h}\cdot\Disoo^{-1}
	\cdot\Disot\cdot\[pmatrix]{I_h\\0}\cdot\pppat\cdot\biso_2^{-1}.
	}
	Combining the preceding equations \eqref{bashibashi},\eqref{bishibishi} and \eqref{bachai}, the intersection number 
	equals
	\[equation]{\lb{chaiji}
\left|\Nrd\left(\biso_1^2\pipao \qaqao\right)\right|_F^{-1}\Big|\Nrd\left(\rr0{I_h}\cdot\Disoo^{-1}
\cdot\Disot\cdot\cc{I_h}0\right)\Big|_F^{-1}.
	}
Note that
$$
	\left|\Nrd\left(\biso_1^2\pipao \qaqao\right)\right|_F=\left|\mathrm{NRD}\left(\mm{\biso_1 I_h}{}{}{\biso_1I_h}\mm{P_{\blin_1}}{*}{0}{Q_{\blin_1}}_{\GG_{K_1}}\right)\right|_F.
$$
Since $|\det\Gamma_{\blin_1}|_F=1$, the above equation can be simplifeid to
	\[equation]{\lb{shayi}|\mathrm{NRD}(\Disoo)|_F=|\mathrm{Disc}_{K/F}|_F^{h^2}
	}
	by Lemma \ref{disco}. Combining \eqref{chaiji} and \eqref{shayi}, we complete the proof of this theorem. 
\end{proof}

\subsection{The invariant polynomial and the resultant formula}\label{xia}\comments{xia}
In this section, we restrict ourselves to the case where the two fixed embeddings $K_1\subset D_F$ and $K_2\subset D_F$ are the same. In particular, we have $K_1\cong K_2$ and $\GG_{K_1}=\GG_{K_2}$. We abbreviate them to $K\subset D_F$ and $\GK$ respectively. We introduce the invariant polynomials to simplify the formula \eqref{fofofo} further. We abbreviate $\biso:=\biso_1$,$\blin:=\blin_1$, and note that there exists elements $g\in\GL_{2h}(F)$ and $\ACj\in D_F^\times$ describing the relative position of $\biso_1,\biso_2,\blin_1$ and $\blin_2$ by
\[equation]{\lb{tratra}
\biso_2=\ACj\circ\biso_1,\qquad M_{\blin_2}=gM_{\blin_1},
}
where in this case, we have
\[equation]{\lb{endless}
\Delta(\biso_1,\blin_1)^{-1}\cdot\Delta(\biso_2,\blin_2)=\Delta(\biso,\blin)^{-1}\cdot\ACj\cdot\ACg\cdot\Delta(\biso,\blin).
}
Therefore, to calculate $\chi(\ca{\biso}{\blin}\infty\litimes_{\GF^{2h}}\ca{\ACj\biso}{\ACg\blin}\infty)$, it suffices to simplify the following expression
\[equation]{\lb{samesame}\left|\Nrd\left(\rr0{I_h}\cdot\Diso^{-1}\cdot\ACj\cdot\ACg\cdot\Diso\cdot\cc{I_h}0\right)\right|_F^{-1}.}
Letting $D_K\subset D_F$ be the centralizer of $K$, we note that $\mathcal O_{D_K}:=D_K\cap \mathcal O_{D_F}$ is the endomorphism ring of $\GK$ as a formal $\OO K$-module. Note that the maps $\biso:\GK\rightarrow\GF$ and $\blinn:F^{2h\vee}\rightarrow K^{h\vee}$ induce $F$-algebra embeddings 
$$\matt hhK\rightarrow \matt{2h}{2h}F\qquad\text{and}\qquad\OD Kh\rightarrow \OD F{2h}$$
by $x\mapsto g$ and $y\mapsto \ACj$ such that $\blan x = g\blan$ and $\biso\circ y=\ACj\circ\biso$ respectively. Denote their images by $\matt{2h}{2h}F_+$ and $D_{F_+}$. Note that
$$
\matt{2h}{2h}F_+\cong \matt hhK,\qquad D_{F+}=\biso\cdot D_K\cdot \bisoo.
$$
Let $K'\subset\matt{2h}{2h}F_+$ be the center of $\matt{2h}{2h}F_+$. Obviously $K'\cong K$. Define
$$\matt{2h}{2h}F_-:=\{\ACg\in\matt{2h}{2h}F: \ACg k=k^\sigma \ACg \quad \text{ for any }\; k\in K'\},$$
and
$$D_{F-}:=\{\ACj\in D_F: \ACj k=k^\sigma \ACj \quad\text{ for any }\; k\in \biso K\biso^{-1}\}.$$
Then we have
\[equation]{\lb{decompa}
\matt{2h}{2h}F=\matt{2h}{2h}F_+\oplus\matt{2h}{2h}F_-,\qquad D_{F}=D_{F+}\oplus D_{F-}
}
because they are the eigenspace decompositions of $D_F$ and $\matt{2h}{2h}F$ as left $K$-linear spaces under the right $K$-multiplications. With respect to this decomposition, for any $\ACj\in\OD Fh$ and any $\ACg\in\matt{2h}{2h}F$, we write
$$\ACj=\ACj_++\ACj_-,\qquad \ACg=\ACg_++\ACg_-.$$
When $\ACj$ (resp. $\ACg$) is invertible, conjugating it by a trace 0 element $\mu\in \bisoo\cdot K^\times\cdot\biso\subset D_F^\times$ (resp. $\mu\in K'^\times\subset\GL_{2h}(F)$), we have
\[equation]{\lb{eeu}
\mu(\ACj_++\ACj_-)\mu^{-1}=(\ACj_+-\ACj_-),\qquad\text{resp.}\quad
\mu(\ACg_++\ACg_-)\mu^{-1}=(\ACg_+-\ACg_-),
}
which implies that $\ACj_+-\ACj_-$(resp. $\ACg_+-\ACg_-$) is also invertible. 
\begin{defi}
For any $\ACj\in\OD F{2h}^\times$ and $\ACg\in\GL_{2h}(F)$, we define
	\begin{equation}\label{tiana}\comments{tiana}\[split]{
			\ACj^\prime_\biso:=(\ACj_+-\ACj_-)^{-1}\ACj_+(\ACj_++\ACj_-)^{-1}\ACj_+&\in D_F,\\
		\ACg^\prime_\blin:=(\ACg_+-\ACg_-)^{-1}\ACg_+(\ACg_++\ACg_-)^{-1}\ACg_+&\in\matt{2h}{2h}F.
}
\end{equation}
\end{defi}
\begin{rem}\lb{rerere}
	When $\ACj_+$ (resp. $\ACg_+$) is invertible, we can write
$$\ACj^\prime_\biso=(\ACj_+-\ACj_-)^{-1}\ACj_+(\ACj_++\ACj_-)^{-1}\ACj_+=(1-\ACj_+^{-1}\ACj_-\ACj_+^{-1}\ACj_-)^{-1},$$
$$\ACg^\prime_\blin=(\ACg_+-\ACg_-)^{-1}\ACg_+(\ACg_++\ACg_-)^{-1}\ACg_+=(1-\ACg_+^{-1}\ACg_-\ACg_+^{-1}\ACg_-)^{-1}.$$
\end{rem}
\begin{lem}We have $\ACj_\biso^\prime\in D_{F+}=\biso\cdot\OD Kh\cdot\bisoo$ and $\ACg_\blin^\prime\in {\matt{2h}{2h}F_+}$.
\end{lem}
\begin{proof}Note that the assignment $\ACg\mapsto\ACg_\biso'$ induces a morphism $\GL_{2h}(F)\rightarrow\matt{2h}{2h}F$ of algebraic varieties over $F$. Since $\matt{2h}{2h}F_+$ is a Zariski closed subset, it suffices to prove this lemma for the following Zariski-dense subset
	$$
	\{\ACg\in\GL_{2h}(F):\ACg_+\in\GL_h(K)\}.
	$$
By Remark \ref{rerere}, we can write $\ACg'_\blin=(1-(\ACg_+^{-1}\ACg_-)^2)^{-1}$. Note that for any $k\in K'$, we have
$$
\ACg_+^{-1}\ACg_-\ACg_+^{-1}\ACg_- k = \ACg_+^{-1}\ACg_- k^\sigma \ACg_+^{-1}\ACg_- = k\ACg_+^{-1}\ACg_-\ACg_+^{-1}\ACg_-,
$$
	which implies $(\ACg_+^{-1}\ACg_-)^2\in\matt{2h}{2h}F_+$, and therefore $\ACg'_\blin\in\matt{2h}{2h}F_+$. A similar argument shows that $\ACj_\biso^\prime\in D_{F+}$, which completes the proof.
\end{proof}
\begin{defi}\lb{invap}
	We define the invariant polynomial of $\ACj\in\OD F{2h}^\times$ (resp. $\ACg\in\GL_{2h}(F)$) with respect to $\biso$ (resp. $\blin$) as the $K$-coefficient reduced characteristic polynomial of $\ACj_\biso^\prime$ as an element in the $K$-central simple algebra $D_{F+}$ (resp. $\matt{2h}{2h}F_+$). We denote this polynomial by $P_\ACj^{\biso}$ (resp. $P_\ACg^\blin$). We also abbreviate them to $P_\ACj$ and $P_\ACg$ if the choice of $\biso$ and $\blin$ is clear from the context.
\end{defi}
From the definition, we know that these polynomials are of degree $h$. 
\begin{cor}\lb{jiangdiao}
	For any $x\in F^\times$, we have $P_{x\gamma}=P_{\gamma}$ and $P_{xg}=P_g$ for any $\ACj\in D_F^\times$ and $g\in\GL_{2h}(F)$.
\end{cor}
\[proof]{This is obvious since $(x\gamma)_+=x\gamma_+$, $(x\gamma)_-=x\gamma_-$, $(xg)_+=xg_+$ and $(xg)_-=xg_-$. }
\begin{prop}\lb{coeffa}
	The coefficients of the polynomials $P_\ACj^{\biso}$ and $P_\ACg^\blin$ are in $F$.
\end{prop}
\begin{proof}
	It suffices to prove this proposition for a Zariski-dense subset consisting of $\ACj\in D_F^\times$ such that $\ACj_+,\ACj_-\in D_F^\times$. Let $p(x)=x^n+\sum_{i=0}^{n-1} a_ix^i$ be the minimal polynomial of $\ACj'_{\biso}$ with coefficients $a_i\in \biso K\biso^{-1}\subset D_F$. Since $\ACj_+^{-1}\ACj_- a_i = a_i^\sigma \ACj_+^{-1}\ACj_-$, and $\ACj_+^{-1}\ACj_-$ commutes with $\ACj'_{\biso}$, we have $p^\sigma(\ACj'_{\biso})=\ACj_+^{-1}\ACj_- p(\ACj'_{\biso})\ACj_-^{-1}\ACj_+=0$, which implies that $p^\sigma(x)$ is also the minimal polynomial of $\ACj'_{\biso}$. Therefore $p(x)$ has coefficients in $F$, which implies $P_\ACj^{\biso}(x)\in F[x]$. A similar method will show that $P_\ACg^{\blin}(x)\in F[x]$. We complete the proof.
\end{proof}
\begin{defi}\label{RESA}\comments{RESA}
For any $\ACj\in\OD Fh$ and $\ACg\in\gth$, we define the relative resultant
\begin{equation}\label{shashi}\comments{shashi}
	\mathrm{Res}(P_\ACj^\biso,P_\ACg^\blin):=\prod_{1\leq i,j\leq h}(x_i-y_j)
\end{equation}
as the resultant of polynomials $P_\ACj^\biso$ and $P_\ACg^\blin$. Here $x_i$ and $y_i$ are roots of $P_\ACj^\biso$ and $P_\ACg^\blin$ for $1\leq i\leq h$, with multiplicities counted.
\end{defi}
Using \eqref{samesame}, we will prove the following theorem. 
\begin{thm}\label{xiaomubiao}\comments{xiaomubiao}
	Suppose $(\biso,\blin)\in\equi KFh$ and $(\ACj,\ACg)\in\equi FF{2h}$.
	We have
	$$\chi(\ca{\biso}{\blin}\infty\litimes_{\GF^{2h}}\ca{\ACj\biso}{\ACg\blin}\infty)=\yabi\cdot\|\mathrm{Res}(P_\ACj^\biso,P_\ACg^\blin)\|_F^{-1}
$$when the right hand side is finite. 
\end{thm}
This theorem is a direct consequence of the following proposition.
\begin{prop}\lb{duijie}
	For any $\ACj\in D_F^\times, \ACg\in\GL_{2h}(F)$ and $(\biso,\blin)\in\equi KFh$, we have
$$
	\left|\Nrd\left(\rr0{I_h}\cdot\Diso^{-1}\cdot\ACj\cdot\ACg\cdot\biso\cdot\blan\right)\right|_F^{-1}=\|\mathrm{Res}(P_\ACj^\biso,P_\ACg^\blin)\|_F^{-1}|\det\ACg\cdot \nrd(\ACj)|_F^{-h}.
$$
\end{prop}
\begin{proof}
It suffices to prove this proposition for those $\ACj$ such that $\ACj_+$ is invertible. Since
$$
	\[pmatrix]{I_h&0_h}\cdot\Diso^{-1}\cdot\Diso\cdot\[pmatrix]{I_h\\0_h}=
	\[pmatrix]{0_h&I_h}\cdot\Diso^{-1}\cdot\Diso\cdot\[pmatrix]{0_h\\I_h}=
	I_h,
$$
$$
	\[pmatrix]{0_h&I_h}\cdot\Diso^{-1}\cdot\Diso\cdot\[pmatrix]{I_h\\0_h}=
	\[pmatrix]{I_h&0_h}\cdot\Diso^{-1}\cdot\Diso\cdot\[pmatrix]{0_h\\I_h}=
	0_h,
$$
	and $\Diso=\biso\cdot(M_\blin\; M_\blin^\sigma)$, we have
	\[equation]{\lb{yaya}
\[pmatrix]{I_h&0_h}\cdot\Diso^{-1}\cdot\biso\cdot\blan=
	\[pmatrix]{0_h&I_h}\cdot\Diso^{-1}\cdot\biso\cdot\blan^\sigma=
	I_h,
	}
	\[equation]{\lb{yoyo}
\[pmatrix]{I_h&0_h}\cdot\Diso^{-1}\cdot\biso\cdot\blan^\sigma=
	\[pmatrix]{0_h&I_h}\cdot\Diso^{-1}\cdot\biso\cdot\blan=
	0_h.
	}
	Let $$X:=\rr0{I_h}\cdot\Diso^{-1}\cdot\ACj\cdot\ACg\cdot\biso\cdot\blan.$$
	It suffices to prove
$
	\Nrd\left(X\right)^4=\mathrm{nrd}(\ACj)^{4h}\det(\ACg)^{4h}\Res_{\biso,\blin}(\ACj,\ACg)^4
$ for those $\ACj,\ACg$ such that $\ACj_+$, $\ACg_+$ are invertible.
By constructions of $g_+$ and $g_-$, there exists $h\times h$ matrices $x_+,x_-\in\matt hhK$ such that
	\[equation]{\lb{mumu}
\ACg_+\blan = \blan x_+,\qquad \ACg_-\blan = \blan^\sigma x_-.
	}
Similarly for $\ACj_+$ and $\ACj_-$, we have
	\[equation]{\lb{miemie}
	(\biso^{-1}\cdot \ACj_+\cdot\biso) \cdot \blan = \blan\cdot (\biso^{-1}\cdot \ACj_+\cdot\biso),\quad
	(\biso^{-1}\cdot \ACj_-\cdot\biso) \cdot \blan = \blan^\sigma\cdot( \biso^{-1}\cdot \ACj_-\cdot\biso).
	}
Using preceding equations \eqref{yoyo},\eqref{mumu},\eqref{miemie}, we have
	\[equation]{\lb{pp}
\[pmatrix]{0_h&I_h}\cdot\Diso^{-1}\cdot\ACj_+\cdot\biso\cdot\ACg_+\cdot\blan=
	0_h,}
	\[equation]{\lb{mm}
\[pmatrix]{0_h&I_h}\cdot\Diso^{-1}\cdot\ACj_-\cdot\biso\cdot\ACg_-\cdot\blan=
	0_h,}
	\[equation]{\lb{pm}
\[pmatrix]{I_h&0_h}\cdot\Diso^{-1}\cdot\ACj_+\cdot\biso\cdot\ACg_-\cdot\blan=
	0_h,}
	\[equation]{\lb{mp}
\[pmatrix]{I_h&0_h}\cdot\Diso^{-1}\cdot\ACj_-\cdot\biso\cdot\ACg_+\cdot\blan=
	0_h.}
	 Note that $\ACj=\ACj_++\ACj_-$, that $\ACg=\ACg_++\ACg_-$, and that $g_+,g_-$ commute with $\ACj_+,\ACj_-$ and $\biso$. Using equations \eqref{pp} and \eqref{mm}, we have
	\[equation]{\lb{xialiu}
	X = \rr{0_h}{I_h}\cdot\Diso^{-1}\cdot(\ACj_+\cdot\ACg_-+\ACj_-\cdot\ACg_+)\cdot\biso\cdot\blan.
}
	By \eqref{miemie}, we have $(\bisoo\cdot\ACj_-\cdot\ACj_+^{-1}\biso)\cdot \blan =\blan^\sigma (\bisoo\cdot\ACj_-\cdot\ACj_+^{-1}\biso)$, which implies
$$
	\ACj_-\cdot\ACj_+^{-1}\cdot\Diso = \Diso\[pmatrix]{&I_h\\I_h}\cdot(\bisoo\cdot\ACj_-\cdot\ACj_+^{-1}\cdot\biso).
$$
Letting
$$
X':=(\bisoo\cdot\ACj_-\cdot\ACj_+^{-1}\cdot\biso)^{-1}\cdot X\cdot (\bisoo\cdot\ACj_+^{-1}\cdot\ACj_-\cdot\biso),
$$
we have
\[equation]{\lb{taixialiu}
X'=\rr{I_h}{0_h}\cdot\Diso^{-1}\cdot(\ACj_+\cdot\ACg_-+\ACj_-\cdot\ACg_+)\cdot\biso\cdot\blan^\sigma.
}
Combining \eqref{xialiu} and \eqref{taixialiu}, we have
$$
\Diso^{-1}\cdot(\ACj_+\cdot\ACg_-+\ACj_-\cdot\ACg_+)\cdot\Diso = \[pmatrix]{0_h&X'\\X&0_h}.
$$
Since $\mathrm{Nrd}(X)=\mathrm{Nrd}(X')$, we have
$$
(-1)^h\mathrm{Nrd}(X)^2=\mathrm{NRD}\[pmatrix]{0_h&X'\\X&0_h}=\mathrm{NRD}(\ACj_+\ACg_-+\ACj_-\ACg_+)
=\mathrm{NRD}(\ACj_+(\ACj_+^{-1}\ACj_-+\ACg_-\ACg_+^{-1})\ACg_+).
$$
Letting $\mu\in \biso\cdot K^\times\cdot\bisoo$ be such that $\mu^\sigma=-\mu$, we have
$$
\mu(\ACj_+^{-1}\ACj_-+\ACg_-\ACg_+^{-1})\mu^{-1}=-\ACj_+^{-1}\ACj_-+\ACg_-\ACg_+^{-1},
$$
which implies $\mathrm{NRD}(\ACj_+^{-1}\ACj_-+\ACg_-\ACg_+^{-1})=\mathrm{NRD}(-\ACj_+^{-1}\ACj_-+\ACg_-\ACg_+^{-1})$. Therefore,
\[equation]{\lb{xx}
\mathrm{Nrd}(X)^4=\mathrm{NRD}(\ACj_+^2\ACg_+^2)\mathrm{NRD}((\ACg_-\ACg_+^{-1})^2-(\ACj_+^{-1}\ACj_-)^2)=\mathrm{NRD}(\ACj_+^2\ACg_+^2)\mathrm{NRD}((\ACj'_\biso)^{-1}-(\ACg'_\blin)^{-1}).
}
By \eqref{eeu}, we have $\mathrm{NRD}(g_++g_-)=\mathrm{NRD}(g_+-g_-)$ and  $\mathrm{NRD}(\ACj_++\ACj_-)=\mathrm{NRD}(\ACj_+-\ACj_-)$, which implies 
\[equation]{\lb{xxx}
\mathrm{NRD}(\ACj'_\biso) = \mathrm{NRD}(\ACj_+^2)\mathrm{NRD}(\ACj^{-2}),\qquad
\mathrm{NRD}(\ACg'_\blin) = \mathrm{NRD}(\ACg_+^2)\mathrm{NRD}(\ACg^{-2}).
}
Note that, as elements of $K$-central simple algebras, the characteristic polynomials of $\ACj'_\biso\in D_{F+}$ and $\ACg'_\blin\in\matt{2h}{2h}F_+$ are given by $P_\ACj^\biso$ and $P_\ACg^\blin$. Therefore, as elements in $D_F$ and $\matt{2h}{2h}F$, their characteristic polynomials are $(P_\ACj^\biso)^2$ and $(P_\ACg^\blin)^2$. Let $L\subset D_{F+}$ be a maximal subfield containing $\ACj'_\biso$ and $\biso\cdot K\cdot\bisoo$. Since $\ACj'_\biso$ commutes with $\ACg'_\blin$, we have 
$$\ACj'_\biso-\ACg'_\blin\in\matt{2h}{2h}L\subset\matt{2h}{2h}{D_F}.$$
Therefore $\det_L(\ACj'_\biso-\ACg'_\blin)=(P_g^\blin(\ACj'))^2$, and we have
\[equation]{\lb{xxxx}
\mathrm{NRD}(\ACj'_\biso-\ACg'_\blin)=\mathrm{Nm}_{L/F}(P_g^\blin(\ACj'))^2=\mathrm{Res}((P_\ACg^\blin)^2,(P_\ACj^\biso)^2) = \mathrm{Res}(P_\ACg^\blin,P_\ACj^\biso)^4.
}
Combining preceding equations \eqref{xx},\eqref{xxx},\eqref{xxxx}, we have
$$
\mathrm{Nrd}(X)^4=\mathrm{NRD}(\ACj_+^2\ACg_+^2)\mathrm{NRD}((\ACj'_\biso)^{-1}(\ACg'_\blin)^{-1})\mathrm{NRD}(\ACj'_\biso-\ACg'_\blin)=\mathrm{NRD}(\ACj^2\cdot\ACg^2)\mathrm{Res}(P_\ACg^\blin,P_\ACj^\biso)^4.
$$
Noting that $\mathrm{NRD}(\ACj)=\mathrm{nrd}(\ACj)^{2h}$ and $\mathrm{NRD}(\ACg)=\mathrm{det}(\ACg)^{2h}$, we complete the proof.
\end{proof}

\begin{proof}[Proof of the Theorem \ref{xiaomubiao}]
	By \eqref{endless} and Theorem \ref{yanchanghui}, we know that 
$$
\chi(\ca{\biso}{\blin}\infty\litimes_{\GF^{2h}}\ca{\ACj\cdot\biso}{\ACg\cdot\blin}\infty)=\yabi\left|\Nrd\left(\rr0{I_h}\cdot\Diso^{-1}\cdot\ACj\cdot\ACg\cdot\biso\cdot\blan\right)\right|_F^{-1}.
$$
	Since $(\ACj,\ACg)$ is an equi-height pair, we have $|\det(\ACg)\nrd(\ACj)|_F=1$. the theorem follows from the Proposition \ref{duijie}.
\end{proof}\lb{nonzero}
Finally we prove the following result, which will be used in the rest of the article.
\begin{prop}\lb{nonde}
	For any $\ACj\in D_F^\times$ and $\ACg\in\GL_{2h}(F)$, if $P^\biso_\ACj$ is an irreducible polynomial and $P^\biso_\ACj(0)P^\biso_\ACj(1)\neq 0$, then $\mathrm{Res}_{\biso,\blin}(\ACj,\ACg)\neq 0$. 
\end{prop}
\begin{proof}
	Since $P_\ACj^\biso$ and $P_\ACg^\blin$ are monic polynomials of degree $h$, it suffices to show that $P_\ACj^\biso\neq P_\ACg^\blin$. Suppose for the sake of contradiction that $P_\ACj^\biso= P_\ACg^\blin$.  Since these polynomials are irreducible and $P^\biso_\ACj(0)\neq 0$, we have $\mathrm{nrd}(\ACj'_\biso)\neq 0$ and $\det(\ACg'_\blin)\neq 0$, which implies that $\ACj_+$ and $\ACg_+$ are invertible. In this case, by Remark \ref{rerere}, the characteristic polynomial of $(\ACj_+^{-1}\ACj_-)^2$ and $(\ACg_+^{-1}\ACg_-)^2$ will be the same. Since $P^\biso_\ACj$ is irreducible and $P^\biso_\ACj(1)\neq 0$, we have $\ACj_+^{-1}\ACj_-\neq 0$. Note that there exists matrices $x_+,x_-\in\matt hhK$ such that $\ACg_+M_\blin=M_\blin x_+$ and $\ACg_-M_\blin=M_\blin^\sigma x_-$. This implies that $(\ACg_+^{-1}\ACg_-)^2M_\blin=M_\blin xx^\sigma$ where $x=x_+^{-1}x_-$. So $P_{g}^\blin$ is the characteristic polynomial of $xx^\sigma$. Let $L=F[(\ACj_+^{-1}\ACj_-)^2]$, and let $D_L\subset D_F$ be the centralizer of $L$. Since $[L:F]=h$, one sees that $D_L$ is a quaternion algebra over $L$. For any $k\in K$, since both $\ACj_+^{-1}\ACj_-$ and $k$ commute with $(\ACj_+^{-1}\ACj_-)^2$, we have $\ACj_+^{-1}\ACj_-\in D_L$ and $K\subset D_L$. Let $\lambda = (\ACj_+^{-1}\ACj_-)^2$. Since $\lambda$ can be identified with an eigenvalue of $xx^\sigma$, there exists $y\in LK$ such that $\lambda =yy^{\sigma_L}$ where $\sigma_L$ is the non-trivial element in $\Gal(LK/L)$. Since $\ACj_+^{-1}\ACj_-y=y^{\sigma_L}\ACj_+^{-1}\ACj_-$, we have
	$$
	0 = (\ACj_+^{-1}\ACj_-)^2-yy^{\sigma_L} = (\ACj_+^{-1}\ACj_- + y)(\ACj_+^{-1}\ACj_- - y^{\sigma_L}),
	$$
	which implies either $\ACj_+^{-1}\ACj_-=-y$ or $\ACj_+^{-1}\ACj_- = y^{\sigma_L}$. This contradicts to the assumption that $\ACj_+^{-1}\ACj_-$ does not commute with some elements in $KL$.
\end{proof}

\section{Proof of main theorems}\label{FORMULAA}\comments{FORMULAA}
Combining Theorems \ref{THM:ICT} and \ref{yanchanghui} in previous sections, we have essentially calculated the intersection number $\chi(\cycn1n\litimes_{\DEF n}\cycn2n)$ for sufficiently large $n$, where the minimal choice of $n$ depends on the pairs $\biso_1,\biso_2$ and $\blin_1,\blin_2$. This illustrates the fact that the intersection number induced by two fixed pairs will eventually be stable as the level increases. 
In this section, we eliminate the restriction on $n$ by proving a formula relating intersection numbers on higher and lower levels. Then we combine all formulas together to give a proof of the main theorems.
\subsection{Intersection number formula for low levels}
Let 
$$
R_n:=\ker\left(\xymatrix{
	\GL_{2h}(\OO F)\ar[rr]^{g\mapsto \overline g}&&\GL_{2h}(\OO F/\pi^n)
}\right).
$$
That is, $R_0=\GL_{2h}(\OO F)$, and $R_n=I_{2h}+\pi^n\matt{2h}{2h}{\OO F}$ for any $n\geq 1$. By associating each $g\in \GL_{2h}(\mathcal O_F)$ to the morphism
$$
\map{\beta^n_n(\id,g)}{\DEF n}{\DEF n}
$$
as defined in \eqref{haww}, the group $\GL_{2h}(\mathcal O_F)$ naturally acts on $\DEF n$ for any $n\geq 0$. This action is the natural action of $\GL_{2h}(\mathcal O_F)$ on $\pi^n$-level structures of $\DEF n$.
Note that $R_n\subset \GL_{2h}(\mathcal O_F)$ is the maximal subgroup that acts trivially on $\DEF n$, and that the quotient group $R_n/R_{n+m}$ is the automorphism group of $\DEF{n+m}$ as a formal scheme over $\DEF n$ via the transition map $\map{\pb nm}{\DEF{n+m}}{\DEF n}$ for any $n,m\geq 0$. We put additional subscripts and write $\DEf{1,n}$, $\DEf{2,n}$ for $\pi^n$-level Lubin--Tate spaces associated to $\GG_{K_1}$ and $\GG_{K_2}$ respectively. Moreover, for each $i=1,2$, $n\geq 0$, and $m\geq 0$, we simply denote by $\DEf{i,n+m}\rightarrow\DEf{i,n}$ the corresponding transition map.

We state the key theorem in this subsection as follows. The reader willing to accept this theorem may skip this subsection.
\begin{thm}\lb{project}Suppose $(\biso_i,\blin_i)\in\equi{K_i}Fh$ for $i=1,2$. For any $n,m\geq 0$, we have
\[equation]{\lb{lala}
	\chi(\cycn1{n}\litimes_{\DEF{n}}\ca{\biso_2}{\blin_2}{n}) = \frac{\sum_{k\in R_n/R_{n+m}}\chi(\cycn1{n+m}\litimes_{\DEF{n+m}}\ca{\biso_2}{k\blin_2}{n+m})}{
		\mathrm{deg}(\DEf{1,m+n}\rightarrow \DEf{1,n})
		\mathrm{deg}(\DEf{2,m+n}\rightarrow \DEf{2,n})
		}.
}
\end{thm}
Intuitively, up to some multiplicities, the cycles $\ca{\biso_2}{k\blin_2}{n+m}$ run over all possible cycles that map to $\ca{\biso_2}{\blin_2}{n}$ under the transition map $\DEF{n+m}\rightarrow\DEF n$ when $k$ runs through $R_n/R_{n+m}$.
Before proceeding to the proof, we fix some notations. Let $w_1$,$w_2$ be integers such that $w_1\geq \m(\blin_1)$ and $w_2\geq \m(\blin_2)$. Let
\[equation]{\lb{shava}
\mathcal F_1:=\bco {m+n}{w_1}_*\OO{\DEf{1,m+n+w_1}},\qquad
\mathcal F_2:=\bct {m+n}{w_2}_*\OO{\DEf{2,m+n+w_2}}
}
be coherent sheaves of $\OO{\DEF{m+n}}$-modules, and let $\map{\pb nm}{\DEF{m+n}}{\DEF n}$ be the transition map. Note that by Definition \ref{tcycle}, we have
$$
[\mathcal F_i]=\deg\left(\DEf{i,m+n+w_i}\rightarrow \DEf{i,m+n}\right)Z_{m+n}(\biso_i,\blin_i)
$$for $i=1,2$. We use the following two lemmas to write both sides of the equation \eqref{lala} in terms of $\mathcal F_1$ and $\mathcal F_2$.
\begin{lem}\lb{laba}
	For the left hand side of \eqref{lala}, we have
	\[equation]{\lb{leftleft}
\chi(\cycn1{n}\litimes_{\DEF{n}}\ca{\biso_2}{\blin_2}{n})=\frac{
	\chi(\mathcal F_1\otimes^{\mathbb L}_{\DEF{m+n}}{\pb nm}^*{\pb nm}_*\mathcal F_2)
	}
	{\deg(\DEf{1,m+n+w_1}\rightarrow\DEf{1,n})\deg(\DEf{2,m+n+w_2}\rightarrow\DEf{2,n})}.
	}
\end{lem}
\begin{proof}
	Using the projection formula, we have
\[equation*]{
	\chi(\mathcal F_1\otimes^{\mathbb L}_{\DEF {m+n}}{\pb nm}^*{\pb nm}_*\mathcal F_2)=
\chi({\pb nm}_*\mathcal F_1\otimes^{\mathbb L}_{\DEF n}{\pb nm}_*\mathcal F_2).
	}
	Note that ${\pb nm}_* Z_{m+n}(\biso_i,\blin_i) = \deg\left(\DEf{i,m+n}\rightarrow\DEf{i,n}\right)Z_n(\biso_i,\blin_i)$ for $i=1,2$, which implies that
$$
	\left[{\pb nm}_*\mathcal F_i\right] = \deg(\DEf{i,m+n+w_i}\rightarrow\DEf{i,n})\ca{\biso_i}{\blin_i}{n}\quad \text{for }i=1,2.
$$
Combining all preceding equations together, we complete the proof of the lemma.
\end{proof}
\begin{lem}\lb{raba}
	For each summand of the right hand side of \eqref{lala}, we have
	\[equation]{\lb{rightright}
	\chi(\cycn1{n+m}\litimes_{\DEF{n+m}}\ca{\biso_2}{k\blin_2}{n+m})=\frac{\chi(\mathcal F_1\otimes_{\DEF{m+n}} k^{-1*}\mathcal F_2)}{\deg(\DEf{1,m+n+w_1}\rightarrow\DEf{1,m+n})\deg(\DEf{2,m+n+w_2}\rightarrow\DEf{2,m+n})}.
}
\end{lem}
\begin{proof}
This lemma follows directly from the observation that 
$$
	[k^{-1*}\mathcal F_i] = [k_{*}\mathcal F_i] =\deg(\DEf{i,m+n+w_i}\rightarrow\DEf{i,m+n}) \ca{\biso_i}{k\blin_i}{n+m}
$$
	for $i=1,2$ and any $k\in\GL_{2h}(\OO F)$.
\end{proof}

Combining \eqref{leftleft} and \eqref{rightright} together, it suffices to prove that
\[equation]{\lb{topprove}
\chi(\mathcal F_1\otimes^{\mathbb L}_{\DEF{n+m}}{\pb nm}^*{\pb nm}_*\mathcal F_2)=\chi\left(\mathcal F_1\litimes_{\DEF{m+n}} \bigoplus_{k\in R_n/R_{n+m}} k^*\mathcal F_2\right).
}
Therefore, it is important to compare ${\pb nm}^*{\pb nm}_*\mathcal F_2$ with $\oplus_{k\in R_n/R_{n+m}} k^*\mathcal F_2$. The following lemma is useful to construct a comparison morphism between them.
\begin{lem}\lb{inja}There is an injective canonical morphism of sheaves of $\mathcal O_{\DEF{m+n}}$-modules
	\[equation]{\lb{injinj}
	\xymatrix{\pbpb^*\pbpb_*\OO{\DEF{m+n}}\ar[r]& \bigoplus_{k\in R_n/R_{n+m}}\haofan^*\OO{\DEF{m+n}}}
	}
such that its induced map 
	\[equation]{\lb{isoiso}
	\xymatrix{\pbpb^*\pbpb_*\OO{\DEF{m+n}}\left[\frac1\pi\right]\ar[r]^\cong& \bigoplus_{k\in R_n/R_{n+m}}\haofan^*\OO{\DEF{m+n}}\left[\frac1\pi\right]}
	}
is an isomorphism.
\end{lem}
\begin{proof}
	Since $\DEF{m+n}$ is affine, it suffices to describe a coherent sheaf $\mathcal F$ by taking its global sections. In this proof, by abuse of notation, we will use the same symbol $\mathcal F$ to denote the set of global sections $\mathcal F(\DEF{m+n})$. 
For any $s\in \OO{\DEF{m+n}}$, we denote $k^*s$ by $s^k$. Note that by definition, we have the following isomorphisms of $\OO{\DEF{m+n}}$-modules
$$
{\pb nm}^*{\pb nm}_*\mathcal F_2\cong \OO{\DEF{m+n}}\otimes_{\DEF n}\mathcal F_2,\qquad
k^*\mathcal F_2\cong \OO{\DEF{m+n}}\otimes_{\DEF{m+n},k}\mathcal F_2.
$$
Here the twisted tensor product $\OO{\DEF{m+n}}\otimes_{\DEF{m+n},k}\mathcal F_2$ is defined so that $a^k\otimes v = 1\otimes av$ for all $a\in \OO{\DEF{m+n}}$ and $v\in \mathcal F_2$. The assignment $a\otimes b\mapsto a\otimes b$ defines a canonical morphism $$\OO{\DEF{m+n}}\otimes_{\DEF n} \OO{\DEF{m+n}}\longrightarrow \OO{\DEF{m+n}}\otimes_{\DEF{m+n},k} \OO{\DEF{m+n}}$$ of sheaves of $\OO{\DEF{m+n}}$-modules for any $k\in R_n/R_{n+m}$, which gives rise to a canonical morphism 
$$
	\xymatrix{\pbpb^*\pbpb_*\OO{\DEF{m+n}}\ar[r]& \bigoplus_{k\in R_n/R_{n+m}}\haofan^*\OO{\DEF{m+n}}}
$$
	of sheaves of $\OO{\DEF{m+n}}$-modules as desired. By inverting $\pi$, we obtain the rigid generic fiber $\DEF n^\circ$ of $\DEF n$ for any $n\geq 0$. It is well known that the transition maps $\DEF{n+m}^\circ\longrightarrow \DEF n^\circ$ of the rigid generic fibers of Lubin--Tate spaces are \'etale for any $n,m\geq 0$. Since $R_n/R_{n+m}$ acts on the \'etale cover $\DEF{n+m}^\circ\longrightarrow \DEF n^\circ$ faithfully, and 
	$$
	\#\left(R_n/R_{n+m}\right)=\deg(\DEF{n+m}^\circ\longrightarrow \DEF n^\circ),
	$$
	the induced map \eqref{isoiso} is an isomorphism. Since $\DEF {m+n}$ is flat over $\Spf \CO F$, the natural map
	$$
	\pbpb^*\pbpb_*\OO{\DEF{m+n}}\longrightarrow \pbpb^*\pbpb_*\OO{\DEF{m+n}}\left[\frac1\pi\right]
	$$
	is injective, which implies that the map \eqref{injinj} is also injective.
\end{proof}
In the following two lemmas, we review some sufficient conditions for the vanishing of Euler characteristics.

\begin{lem} [Serre's multiplicity vanishing theorem]
Let $R$ be a regular local ring, and let $\mathfrak{p}$, $\mathfrak{q}$ be primes of $R$ such that $R/\mathfrak{p}\otimes_R R/\mathfrak{q}$ is an Artinian ring. If $\dim(R/\mathfrak{p})+\dim(R/\mathfrak{q})<\dim(R)$, then 
$$\chi(R/\mathfrak{p}\litimes_R R/\mathfrak{q})=0.$$
\end{lem}
\begin{proof}This was proven in 1985 by Paul C. Roberts\cite{roberts1985vanishing}.\end{proof}
\begin{lem}\label{fuzhu}\comments{fuzhu}
Let $M,N$ be finite modules over a regular Noetherian local ring $A$ such that $M\otimes_A N$ is of finite length. If $\dim(\Supp(M))+\dim(\Supp(N))<\dim(A)$, 
then
$$\chi(M\litimes_A N)=0.$$
\end{lem}
\begin{proof}
Consider a filtration 
$$0=M_n\subset \cdots\subset M_0=M$$ 
such that $M_{i}/M_{i+1}\cong A/\mathfrak{p}_i$ for some associated primes $\mathfrak{p}_i\in\mathrm{Ass}(M)$ for any $0\leq i<n$. Similarly for $N$, let 
$$0=N_r\subset \cdots\subset N_0=N$$
be  a filtration such that $N_{j}/N_{j+1}\cong A/\mathfrak{q}_j$ for some $\mathfrak{q}_j\in\mathrm{Ass}(N)$ for any $0\leq j<r$.	
Then,	
$$\chi(M\litimes_{A}N)=\sum_{\substack{1\leq i< n\\1\leq j< r}}\chi(A/\mathfrak{p}_i\litimes_A A/\mathfrak{q}_j).$$
However, the fact that $$\dim(A/\mathfrak{p}_i)+\dim(A/\mathfrak{q}_j)\leq \dim(\Supp(M))+\dim(\Supp(N))<\dim(A)$$ implies that each summand $\chi(A/\mathfrak{p}_i\litimes_A A/\mathfrak{q}_j)=0$ for any $1\leq i< n,1\leq j< r$ by Serre's multiplicity vanishing theorem. Hence
$\chi(M\litimes_{A}N)=0$ as desired.\end{proof}
Combining all previous lemmas, we prove the following proposition, which implies \eqref{toprove} and therefore proves Theorem \ref{project}. Actually, Proposition \ref{nababa} below proves something stronger than \eqref{toprove}; this is necessary for our application in the next subsection.
\begin{prop}\lb{nababa}
Suppose $(\biso_2,\blin_2)\in\equi{K_2}Fh$. For any integer $m,n\geq 0$, and any coherent sheaf $\mathcal F$ of $\mathcal O_{\DEF{m+n}}$-modules such that $\dim(\Supp(\mathcal F))=h$, we have
\[equation]{\lb{toprove}
\chi(\mathcal F\otimes^{\mathbb L}_{\DEF{n+m}}{\pb nm}^*{\pb nm}_*\ca{\biso_2}{\blin_2}{n+m})=\sum_{k\in R_n/R_{n+m}}\chi\left(\mathcal F\litimes_{\DEF{m+n}}  \ca{\biso_2}{k\blin_2}{n+m}\right).
}
\end{prop}
\begin{proof}
\newcommand{\gas}{\mathcal F}
\newcommand{\grass}{\mathcal J}
Recall $\mathcal F_2$ as defined in \eqref{shava}. It suffices to prove that$$
\chi(\mathcal F\otimes^{\mathbb L}_{\DEF{n+m}}{\pb nm}^*{\pb nm}_*\mathcal F_2)=\chi\left(\mathcal F\litimes_{\DEF{m+n}}  \bigoplus_{k\in R_n/R_{n+m}}k^*\mathcal F_2\right).
$$
Let $\mathcal J$ be the cokernel of the map in \eqref{injinj}. By Lemma \ref{inja}, there is an exact sequence 
\begin{equation}\label{saqua}\comments{saqua}\xymatrix{
0\ar[r]&\pbpb^*\pbpb_*\OO{\DEF{m+n}}\ar[r]&\bigoplus_{k\in R_n/R_{n+m}}\haofan^*\OO{\DEF{m+n}}\ar[r]&\mathcal{J}\ar[r]&0
}.\end{equation}
	Since $k^*\OO{\DEF{m+n}}$ is free, we have $\mathrm{Tor}_{{\DEF{m+n}}}^1(k^*\OO{\DEF{m+n}},\mathcal F_2)=0$. Tensoring the sequence \eqref{saqua} by $\mathcal F_2$, we have an exact sequence
\begin{equation*}\label{lalala}\comments{lalala}
\xymatrix{
0\ar[r]&\ttor_{{\DEF{m+n}}}^1(\mathcal J,\mathcal F_2)\ar[r]&\pbpb^*\pbpb_*\mathcal F_2\ar[r]&\bigoplus_{k\in R_n/R_{n+m}}\haofan^*\mathcal F_2\ar[r]&\mathcal{J}\otimes_{\DEF{m+n}}\mathcal F_2\ar[r]&0,
}
\end{equation*}
	where $\ttor_{{\DEF{m+n}}}^1(\mathcal J,\mathcal F_2)$ is the sheaf of $\OO{\DEF{m+n}}$-modules associated to $\mathrm{Tor}_{\DEF{m+n}}^1(\mathcal J,\mathcal F_2)$. 
Then it suffices to prove that
\begin{equation*}\label{kuangzao}\comments{kuangzao}
\chi\Big(\mathcal F\litimes_{\DEF{m+n}}(\grass\otimes_{\DEF{m+n}}\mathcal F_2)\Big)=0;\qquad \chi\Big(\mathcal F\litimes_{\DEF{m+n}}
	{\ttor_{{\DEF{m+n}}}^1(\mathcal J,\mathcal F_2)
}\Big)=0.
\end{equation*}
	By Lemma \ref{inja}, tensoring the sequence \eqref{saqua} by $\OO{\DEF{m+n}}\left[\frac1\pi\right]$, we obtain an isomorphism 
$${\pbpb^*\pbpb_*\OO{\DEF{m+n}}\Big[\frac1\pi\Big]}\longrightarrow{\bigoplus_{k\in R_n/R_{n+m}}\haofan^*\OO{\DEF{m+n}}\Big[\frac1\pi\Big]},$$ 
which implies that
	\begin{equation}\label{chaosile}\comments{chaosile}\mathcal{J}\otimes_{{\DEF{m+n}}}\OO{\DEF{m+n}}\Big[\frac1\pi\Big]=0.\end{equation} 
		Let $V(\pi)\subset\DEF {m+n}$ be the vanishing locus of $\pi$.		
Therefore, we have $\Supp(\mathcal{J})\subset V(\pi)$, which implies that
		$$\Supp(\mathcal{J}\otimes_{\DEF{m+n}} \mathcal F_2)\subset V(\pi)\cap \Supp(\mathcal F_2),$$
		and that 
		$$\Supp(\ttor_{{\DEF{m+n}}}^1(\mathcal{J},\mathcal F_2))\subset V(\pi)\cap \Supp(\mathcal F_2).$$ 
It then suffices to show that for any coherent sheaf $\mathscr M$ of $\OO{\DEF{m+n}}$-modules with $\Supp(\mathscr M)\subset V(\pi)\cap \Supp(\mathcal F_2)$, one has
		$$
		\chi(\mathcal F\litimes_{\DEF{m+n}}\mathscr M)=0.
		$$
Indeed, since $\pi$ is not a zero-divisor for $\mathcal F_2$, we have
$$
	\dim(V(\pi)\cap \Supp(\mathcal F_2))\leq \dim(\Supp(\mathcal F_2))-1=h-1.
$$
However, by what given, we have
$$
	\dim(\Supp(\mathcal F))=h,
$$
which implies that
$$
	\dim(V(\pi)\cap \Supp(\mathcal F_2))+	\dim(\Supp(\mathcal F))\leq 2h-1 <\dim(\DEF{m+n}).
$$
	Hence $\chi(\mathcal F\litimes_{\DEF{m+n}}\mathscr M)=0$ by Lemma \ref{fuzhu}, which completes the proof.\end{proof}
\[proof]{[Proof of Theorem \ref{project}]By Lemma \ref{laba} and Lemma \ref{raba}, it suffices to prove \eqref{topprove}, which is a special case when $\mathcal F=\mathcal F_1$ in Proposition \ref{nababa}. }
\subsection{A generalization of the intersection number formula for correspondences}\lb{heck}
Note that the group $D_F^\times\times\GL_{2h}(\mathcal O_F/\pi^n)$ acts on $\DEF n^\sim$, and therefore also acts on $\groo{\DEF n^\sim}$. Here the group $\groo{\DEF n^\sim}$ is defined in \S\ref{kgroup}.  In fact, we can consider a more general action on $\groo{\DEF n^\sim}$. Let $\mathscr{H}(R_n\backslash \GL_{2h}(F)/R_n)$ be the set of compactly supported double $R_n$-invariant $\QQ$-valued functions. For some particular elements
$$(\ACj,f)\in D_F^\times\times \mathscr{H}(R_n\backslash \GL_{2h}(F)/R_n),$$ 
we define their actions on $\groo{\DEF n^\sim}$ via certain correspondences. Our main result in this subsection is Theorem \ref{FORMULA}, which is a key theorem to generalize our intersection formula for those correspondences.

Let $\ACg_0\in\GL_{2h}(F)$, and let $n,m$ be arbitrary non-negative integers such that $m\geq\m(g_0)$. Consider a function 
\begin{equation}\label{shenmeya}\comments{shenmeya}
	f_{R_n\ACg_0 R_n,m}(x):=\[cases]{C&x\in R_n\ACg_0 R_n;\\
	0&x\notin R_n\ACg_0 R_0,
	}
\end{equation}
where $C=\deg(\DEF{m+n}\longrightarrow\DEF n)/\Vol(R_n\cdot \ACg_0\cdot R_n)$ is the constant such that 
$$\int_{\GL_{2h}(F)}f_{R_n\ACg_0 R_n,m}(x)\dd x=\deg(\DEF{m+n}\longrightarrow\DEF n).$$
Let $\ACj\in D_F^\times$ be an element such that $\mathrm{Height}(\ACj)=\mathrm{Height}(\ACg_0)$. To the pair $(\ACj,f_{R_n\ACg_0 R_n,m})$ we attach the following correspondence
\begin{equation*}\label{corres}\comments{corres}
\xymatrix{
	\beta_n^n(\ACj,f_{R_n\ACg_0 R_n,m}):\DEF n&&\DEF{n+m}\ar[ll]_(.3){\pbpb}\ar[rr]^{\bpb m\ACj{\ACg_0}}&&\DEF n\\
}.
\end{equation*} 
Recall that two correspondences $X\leftarrow Z\rightarrow Y$ and $X\leftarrow Z'\rightarrow Y$ are the same if and only if there exists an isomorphism $Z\rightarrow Z'$ such that the following diagram commutes
$$
\xymatrix{Z\ar[r]\ar[rd]\ar[d]&X\\
	Y&Z'.\ar[l]\ar[u]}
$$
Therefore replacing $g_0$ in the preceding definition by any element in $R_n\ACg_0R_n$ gives the same correspondence, which implies that our notation $\beta_n^n(\ACj,f_{R_n\ACg_0 R_n,m})$ is well-defined. 

For any  coherent sheaf of $\mathcal O_{\DEF n}$-modules $\mathcal F$ on $\DEF n$, the definition
\[equation]{\lb{dafi}\beta_n^n(\ACj,f_{R_n\ACg_0 R_n,m})^*\mathcal F := \pbpb_*\bpb m\ACj{\ACg_0}^*\mathcal F}
also defines an action of the pair $(\ACj,f_{R_n\ACg_0 R_n,m})$ on $\groo{\DEF n}$.
Note that when $(\ACj,\ACg_0)$ is not an equi-height pair, the same construction gives the action of $(\ACj,f_{R_n\ACg_0 R_n,m})$ on $\groo{\DEF n^\sim}$, but in that case the subset $\groo{\DEF n}\subset \groo{\DEF n^\sim}$ is not invariant. Our main result of this subsection is the following theorem. The reader willing to accept this theorem may skip this subsection.
\begin{thm}\label{FORMULA}\comments{FORMULA}
	Suppose $(\biso_i,\blin_i)\in\equi{K_i}Fh$ for $i=1,2$ and $(\ACj,\ACg_0)\in\equi FF{2h}$.	
Let $f=f_{R_ng_0R_n,m}$. We have
	\begin{equation*}\chi\left(\cycn1n\litimes_{\DEF n}\beta_n^{n}(\ACj,f)^*\cycn2n\right)
		=
		{\sum_{x\in R_n/R_{n+m}}\chi\left(Z_{n}(\ACj\biso_1,g_0x\blin_1)\litimes_{\DEF{n}}{\cycn2{n}}\right)}
		.
\end{equation*}
\end{thm}
\begin{proof}By definition in \eqref{dafi}, we have
$$
\chi\left(\cycn1n\litimes_{\DEF n}\beta_n^{n}(\ACj,f)^*\cycn2n\right)=\chi\left(\cycn1n\litimes_{\DEF n}\pbpb_*\bpb m\ACj{\ACg_0}^*\cycn2n\right),
$$
which, by using the projection formula, equals
$$
	\chi\left({\pp{m}^*\cycn1n} \litimes_{\DEF{n+m}} \pro^{*}{\cycn2n}\right).
$$
Note that $  \dpapapa\cycn1n=\pp{m}_{*}\cycn{1}{n+m}$. We have
	$${
\begin{split}
	&\dpapapa\chi\left({\pp{m}^*\cycn1n} \litimes_{\DEF{n+m}} \pro^{*}{\cycn2n}\right)\\
	=\;&\chi\left({\pp{m}^*\pp{m}_{*}\cycn{1}{n+m}}\litimes_{\DEF{n+m}}{\pro^{*}\cycn2n}\right),
\end{split}
}
$$
which, by Proposition \ref{nababa}, equals 
	\[equation]{\lb{zhaolaopo}
\sum_{x\in R_n/R_{n+m}}\chi(\ca{\biso_1}{x\blin_1}{n+m}\litimes_{\DEF{m+n}}\pro^*\cycn2n).
	}
Using the projection formula again, the preceding equation can be simplified to
\begin{equation}\label{zhongji}\comments{zhongji}
	\sum_{x\in R_n/R_{n+m}}\chi\left({\pro_*\ca{\biso_1}{x\blin_1}{n+m}}\litimes_{\DEF n}{\cycn2n}\right).
\end{equation}
Note that, by Proposition \ref{wanquanbuzhi}, we have
$$\pro_*\ca{\biso_1}{x\blin_1}{n+m}=\dpapapa\ca{\ACj\biso_1}{\ACg_0 x\blin_1}{n}.$$ 
	Plugging this equation back to \eqref{zhongji}, we complete the proof of this theorem. 
\end{proof}
\subsection{Proof of the main theorem}
In this subsection, we combine all results of previous sections together to prove the main theorem. Fix three equi-height pairs $(\biso_1,\blin_1)\in\equi{K_1}Fh$, $(\biso_2,\blin_2)\in\equi{K_2}Fh$ and $(\ACg_0,\ACj)\in\equi{F}F{2h}$. Recall the definition of $f_{R_ng_0R_n,m}$ in \eqref{shenmeya}, and the definition of $\Delta(\biso_i,\blin_i)\in\matt{2h}{2h}{D_F}$ in \eqref{lashi} for $i=1,2$.  We denote
$$
\mathbf R(g):=\mathrm{Nrd}\left(\[pmatrix]{{0_h}&{I_h}}\cdot\Delta(\biso_1,\blin_1)^{-1}\cdot\ACj^{-1}\cdot g^{-1}\cdot\Delta(\biso_2,\blin_2)\cdot\cc{I_h}{0_h}\right)\in F.
$$
Recall that $\yabii|\mathbf R(g)|_F^{-1}$ is the intersection number $\chi(Z_\infty(\ACj\cdot\biso_1,\ACg\cdot\blin_1)\litimes_{\GF^{2h}}Z_\infty(\biso_2,\blin_2))$ by Theorem \ref{yanchanghui}.
Note that our intersection formula in Theorem \ref{yanchanghui} only holds for 
\[equation*]{
\mathbf R(g)
\neq 0. 
}
Intuitively, this condition is describing that the cycle $Z_\infty(\ACj\cdot\biso_1,\ACg\cdot\blin_1)$ intersects properly with the cycle $Z_\infty(\biso_2,\blin_2)$.
When $K_1\cong K_2$, the non-degeneracy condition holds automatically for any $g\in\GL_{2h}(F)$ when $P_\ACj^\biso$ is irreducible and $P_\ACj^\biso(0)P_\ACj^\biso(1)\neq 0$ by Proposition \ref{nonde}. When $K_1\not\cong K_2$, there is a similar criterion based on constructions of invariant polynomials for bi-quadratic cases, which was discussed in detail in \cite{howard2020bi}. When $L=F$, $K_1$, or $K_2$, we define their local zeta functions by $\zeta_L(x):=(1-q_L^{-x})^{-1}$, where $q_L$ is the cardinality of the residue field of $\mathcal O_L$. We state our main theorem as follows.
\begin{thm}\lb{chaochao}
	Assume the non-degeneracy condition
\[equation]{\lb{codi}
\mathbf R(g)
\neq 0
}
for any $g\in R_ng_0R_n$.
Suppose $(\biso_i,\blin_i)\in\equi{K_i}Fh$ for $i=1,2$, $(\ACj,\ACg_0)\in\equi FF{2h}$, $n\geq 0$ and $m\geq \m(g_0)$. Let $f:=f_{R_ng_0R_n,m}$. Then the intersection number 
	$$\chi\left(\cycn1n\litimes_{\DEF n}\beta_n^{n}(\ACj,f)^*\cycn2n\right)$$
	is finite and we have the following formula
$$
\chi\left(\cycn1n\litimes_{\DEF n}\beta_n^{n}(\ACj,f)^*\cycn2n\right)
		=
	C_n\cdot\yabii\int_{\GL_{2h}(F)}f(x)|\mathbf R(x)|_F^{-1}\dd x,	
$$
where $C_n=1$ for any $n>0$, and
	\[equation]{\lb{bagayalo}
C_0=\frac
	{\prod_{i=1}^{h}\zeta_{K_1}(i)\zeta_{K_2}(i)}
	{\prod_{i=1}^{2h}\zeta_F(i)}.
	} 
\end{thm}
Before proceeding to the proof, we state a lemma for simplifying integrals.
\begin{lem}\lb{counting}
	For any real-valued function $X(x)$ on $\GL_{2h}(F)$ and any $g\in\GL_{2h}(F)$, we have
	$$
	\int_{R_n}\int_{R_n}X(ygx)\dd y\dd x = \frac{\Vol(R_n)^2}{\Vol(R_n\cdot g\cdot R_n)}\int_{R_n\cdot g\cdot R_n}X(t)\dd t.
	$$
\end{lem}
\begin{proof}
Applying a substitution $y\rightarrow yx^{-1}\ACg^{-1}$, we have
	\begin{equation}\label{seminar}\comments{seminar}\int_{R_n}\int_{R_n}X(yg x)\dd y\;\dd x = \int_{R_n}\int_{R_n\cdot\ACg x}X(y)\dd y\;\dd x.
\end{equation}
Note that
$$
	x\in R_n\;\text{ and } y\in R_n\cdot gx \iff y\in R_n\cdot g\cdot R_n\;\text{ and }x\in R_n\cap g^{-1}\cdot R_n\cdot y.
$$
By changing the order of the integration, we have
$$
\int_{R_n}\int_{R_n\cdot\ACg x}X(y)\dd y\;\dd x=\int_{R_n\cdot g\cdot R_n}\int_{R_n\cap g^{-1}\cdot R_n\cdot y}X(y)\dd x\;\dd y,
$$
which equals
	\[equation]{\lb{quanshiwode}
	\int_{R_n\cdot g\cdot R_n}\Vol\left(R_n\cap g^{-1}\cdot R_n\cdot y\right)X(y)\dd y.
	}
Note that for any $t\in R_n$, we have 
$$\left(R_n\cap g^{-1}\cdot R_n\cdot y\right)\cdot t=R_n\cap g^{-1}\cdot R_n\cdot yt,$$
and 
$$R_n\cap g^{-1}\cdot R_n\cdot ty=R_n\cap g^{-1}\cdot R_n\cdot y,$$
which imply that the value $\Vol\left(R_n\cap g^{-1}\cdot R_n\cdot y\right)$ does not depend on $y$ for all $y\in R_n\cdot g\cdot R_n$.
Letting $X$ be a function that constantly equals $1$, we obtain
$$
\int_{R_n}\int_{R_n}\dd y\;\dd x = \int_{R_n\cdot g\cdot R_n}\Vol\left(R_n\cap g^{-1}\cdot R_n\cdot y\right)\dd y,
$$
which implies that
$$
	\Vol\left(R_n\cap g^{-1}\cdot R_n\cdot y\right)=\frac{\Vol(R_n)^2}{\Vol(R_n\cdot g\cdot R_n)}.
$$
Combining this equation with \eqref{quanshiwode}, we complete the proof of this lemma.\end{proof}
\begin{proof}[Proof of Theorem \ref{chaochao}]
By Theorem \ref{FORMULA}, we have
$$
\chi\left(\cycn1n\litimes_{\DEF n}\beta_n^{n}(\ACj,f)^*\cycn2n\right)=
{\sum_{x\in R_n/R_{n+m}}\chi\left(Z_{n}(\ACj\biso_1,g_0x\blin_1)\litimes_{\DEF{n}}{\cycn2{n}}\right)},
$$
	which, by Theorem \ref{project}, can be written further into
\[equation]{\lb{intform}
	\frac{\sum_{x\in R_n/R_{n+m}}\sum_{y\in R_{n}/R_{n'}}\chi\left(Z_{n'}(\ACj\biso_1,yg_0x\blin_1)\litimes_{\DEF{n'}}{\cycn2{n'}}\right)}
	{\deg\left(\DEf{1,n'}\rightarrow\DEf{1,n}\right)\deg\left(\DEf{2,n'}\rightarrow\DEf{2,n}\right)}
}
for any $n'\geq n$. 
	Let $U=R_n\ACg_0 R_n$ be the support of $f$. Note that the non-degeneracy condition \eqref{codi} implies that $|\mathbf R(g)|_F^{-1}$ is a continuous function on $U$. Therefore there is a uniform bound $|\mathbf R(g)|_F^{-1}<M$ for any $g\in U$ since $U$ is compact. Using \eqref{howbigN}, we can choose a large enough bound $N$ such that the Theorem \ref{THM:ICT} holds uniformly for $(\ACj\biso_1,k\blin_1)$ and $(\biso_2,\blin_2)$ whenever $k\in U$. Letting $n'$ be a value such that $n'>N$, by Theorem \ref{THM:ICT}, each summand of \eqref{intform} simplifies to
	\[equation]{\lb{infinitel}
\chi\left(Z_{n'}(\ACj\biso_1,yg_0x\blin_1)\litimes_{\DEF{n'}}{\cycn2{n'}}\right)
=
	\chi\left(Z_{\infty}(\ACj\biso_1,yg_0x\blin_1)\litimes_{\GF^{2h}}{\cycn2{\infty}}\right),
	}
	which, by Theorem \ref{yanchanghui}, equals
$$
	\yabii|\mathbf R(yg_0x)|_F^{-1}.
$$
	Since the action of $\beta^{n'}_{n'}(\id,t)$ on $\DEF{n'}$ is trivial for any $t\in R_{n'}$, we have 
$$Z_{n'}(\ACj\biso_1,tyg_0x\blin_1)=\beta^{n'}_{n'}(\id,t)_*Z_{n'}(\ACj\biso_1,yg_0x\blin_1)=Z_{n'}(\ACj\biso_1,yg_0x\blin_1).$$
	Therefore $|\mathbf R(tyg_0x)|_F=|\mathbf R(yg_0x)|_F$ for any $t\in R_{n'}$, which implies that
	\[equation*]{
	|\mathbf R(yg_0x)|_F^{-1}=\frac{\int_{R_{n'}\cdot y}|\mathbf R(tg_0x)|_F^{-1}\dd t}{\int_{R_{n'}}\dd t}.
	}
	Then, summing up the simplified summands over $R_n/R_{n'}$, we have
	\[equation]{\lb{jifenjifen}
\sum_{y\in R_{n}/R_{n'}}\chi\left(Z_{n'}(\ACj\biso_1,yg_0x\blin_1)\litimes_{\DEF{n'}}{\cycn2{n'}}\right)=
\frac{\yabii\int_{R_{n}}|\mathbf R(tg_0x)|_F^{-1}\dd t}{\int_{R_{n'}}\dd t}.
		}
	Combining preceding equations \eqref{intform},\eqref{infinitel}, and \eqref{jifenjifen}, we obtain
\[equation]{\lb{midamida}
\[split]{
	&\chi\left(\cycn1n\litimes_{\DEF n}\beta_n^{n}(\ACj,f)^*\cycn2n\right)\\
	=\;&
	\frac{\yabii\sum_{x\in R_n/R_{m+n}}\int_{R_{n}}|\mathbf R(tg_0x)|_F^{-1}\dd t}
	{
		\deg\left(\DEf{1,n'}\rightarrow\DEf{1,n}\right)\deg\left(\DEf{2,n'}\rightarrow\DEf{2,n}\right)
	\int_{R_{n'}}\dd t
	}.
	}
}
Since $\m(g_0)\leq m$, we have $g_0^\vee(\mathcal O_F^{2h\vee})\supset\pi^m\mathcal O_F^{2h\vee}$, which implies that 
	$$
	g_0\cdot R_{n+m}\cdot g_0^{-1}\subset R_{n}.
	$$
Therefore for any $t\in R_{n+m}$, we have $yg_0tx\in y\cdot R_{n}\cdot g_0x = R_{n}\cdot yg_0x$, which implies that the integral $\int_{R_n}|\mathbf R(yg_0x)|_F^{-1}\dd y$ is invariant when replacing $x\mapsto tx$.  Then we have
	\[equation]{\lb{haipia}
	\sum_{x\in R_n/R_{m+n}}\int_{R_{n}}|\mathbf R(yg_0x)|_F^{-1}\dd y=\frac{\int_{R_n}\int_{R_n}|\mathbf R(yg_0x)|_F^{-1}\dd y\dd x
	}
	{
		\int_{R_{m+n}}\dd x	
		},
	}
	which, by Lemma \ref{counting}, equals
$$
\frac	{\Vol(R_n)^2}{\Vol(R_n\cdot g_0\cdot R_n)\Vol(R_{m+n})}
\;	{\int_{R_n\cdot g_0\cdot R_n}|\mathbf R(g)|_F^{-1}\dd g
	}
	.
$$
Moreover, note that by definition of $f$ in \eqref{shenmeya}, we have
	\[equation]{\lb{hehekaka}
\frac{\deg(\DEF{m+n}\rightarrow \DEF n)}{\Vol(R_n\cdot g_0\cdot R_n)}\int_{R_n\cdot g_0\cdot R_n}|\mathbf R(g)|_F^{-1}\dd g=\int_{\GL_{2h}(F)}f(g)|\mathbf R(g)|_F^{-1}\dd g.
}
Note that $\deg(\DEF{m+n}\rightarrow \DEF n)=\Vol(R_{n})/\Vol(R_{m+n})$.
Combining preceding equations \eqref{midamida},\eqref{haipia}, and \eqref{hehekaka}, the intersection number $\chi\left(\cycn1n\litimes_{\DEF n}\beta_n^{n}(\ACj,f)^*\cycn2n\right)$ equals
\[equation]{\lb{factors}
\frac{\yabii\Vol(R_n)}
	{
		\deg\left(\DEf{1,n'}\rightarrow\DEf{1,n}\right)\deg\left(\DEf{2,n'}\rightarrow\DEf{2,n}\right)
	\Vol(R_{n'})
	}
\int_{\GL_{2h}(F)}f(g)|\mathbf R(g)|_F^{-1}\dd g		
	.
}
To simplify each leading factors, note that for any integer $r$ and $i=1,2$, we have
\[equation]{\lb{zetasmall}
\deg\left(\DEf{i,r}\rightarrow\DEf{i,0}\right)=\#\GL_h(\mathcal O_{K_i}/\pi^{r})=q^{2{r}h^2}\prod_{j=1}^h\left(\zeta_{K_i}(j)\right)^{-1},
}
\[equation]{\lb{defsmall}
\deg\left(\DEf{i,n'}\rightarrow\DEf{i,n}\right)=\frac{\deg\left(\DEf{i,n'}\rightarrow\DEf{i,0}\right)}{\deg\left(\DEf{i,n}\rightarrow\DEf{i,0}\right)},
}
\[equation]{\lb{zetabig}
\Vol(R_0/R_{r})=\#\GL_{2h}(\mathcal O_F/\pi^{r})=q^{4{r}h^2}\prod_{j=1}^{2h}\left(\zeta_F(j)\right)^{-1},
}
and
\[equation]{\lb{defbig}
\Vol(R_n)/\Vol(R_n')=\frac{\Vol(R_0/R_{n'})}{\Vol(R_0/R_{n})}.
}
If $n=0$, then putting $r=n'$ in \eqref{zetasmall} and \eqref{zetabig}, we have
$$
\frac{\Vol(R_0)/\Vol(R_{n'})}
	{
		\deg\left(\DEf{1,n'}\rightarrow\DEf{1,0}\right)\deg\left(\DEf{2,n'}\rightarrow\DEf{2,0}\right)	
	}
	=\frac{q^{4{n'}h^2}\prod_{j=1}^{2h}\zeta_F(j)^{-1}}{\left(\prod_{j=1}^h\zeta_{K_1}(j)^{-1}\right) q^{2{n'}h^2}\left(\prod_{j=1}^h\zeta_{K_2}(j)^{-1}\right) q^{2{n'}h^2}},
$$
where one sees that the term $q^{4n'h^2}$ cancels out. So we complete the proof for $n=0$. When $n>0$, using \eqref{defsmall},\eqref{defbig} together with \eqref{zetasmall}, \eqref{zetabig}, we have $\Vol(R_n)/\Vol(R_{n'})=q^{4(n'-n)h^2}$ and $\deg\left(\DEf{1,n'}\rightarrow\DEf{1,n}\right)=\deg\left(\DEf{2,n'}\rightarrow\DEf{2,n}\right)=q^{2(n'-n)h^2}$, which completes the proof of the theorem for all $n$.
\end{proof}
\begin{proof}[Proof of main Theorems \ref{BIQ},\ref{SUPER},\ref{girls} and \ref{babys}]
	The Theorem \ref{BIQ} is identical with Theorem \ref{chaochao}. Let $\Delta_i=\Delta(\biso_i,\blin_i)$ for $i=1,2$. The Theorem \ref{SUPER} is a special case of Theorem \ref{chaochao} when $K_1\cong K_2$, $\biso_2=\ACj\biso_1$, and $\blin_2=\ACg\blin_1$. In this case we have $\Delta_2=\ACj\ACg\Delta_1$. Since $(\ACj,\ACg)$ is an equi-height pair, we have, by Proposition \ref{duijie},
	$$
	\left|\mathrm{Nrd}\left(\[pmatrix]{{0_h}&{I_h}}\cdot\Delta_1^{-1}\cdot \ACj_0^{-1}\cdot x^{-1}\cdot\ACg\cdot\ACj\cdot\Delta_1\cdot\cc{I_h}{0_h}\right)\right|_F=\left|\mathrm{Res}_{\biso_1,\blin_1}\left(P^{\biso_1}_{\ACj_0^{-1}\ACj},P^{\blin_1}_{x^{-1}\ACg}\right)\right|_F,
	$$
	which is non-zero everywhere for all $g\in\GL_{2h}(F)$ when $P_{\ACj_0^{-1}\ACj}^{\biso_1}$ is irreducible and $P_{\ACj_0^{-1}\ACj}^{\biso_1}(0)P_{\ACj_0^{-1}\ACj}^{\biso_1}(1)\neq 0$ thanks to Proposition \ref{nonde}. Therefore Theorem \ref{chaochao} implies Theorem \ref{SUPER}. The Theorems \ref{girls} and \ref{babys} are special cases of Theorem \ref{SUPER} when $m=0$, $\ACj=\id$, and $g_0=\id$.
\end{proof}

\subsection{Generalization of the theorem to non-equi-height pairs}\lb{removal}
When defining \Heegner cycles $Z_\bullet(\biso,\blin)$ on $\DEF\bullet^\sim$ in Definition \ref{tcycle}, we do not put any restrictions for the height of $(\biso,\blin)$. In this subsection, we emphasize the fact that the intersection number involving non-equi-height pairs is essentially the same as the one corresponding to equi-height pairs.

Let $\varpi_K$ be a uniformizer of $D_K$, and let $\vv_D(\bullet)$ be the normalized discrete valuation of $D_F$. Note that $\vv_D(\varpi_K)=2$ if $K/F$ is unramified and $\vv_D(\varpi_K)=1$ if $K/F$ is ramified.

For any pairs $(\biso,\blin)$ as in \eqref{pairsa}, it is clear from the definition that $Z_n^\sim(\biso,\blin) = Z_n^\sim(\biso\cdot \varpi_K^m,\blin)$ for any $m$. Therefore, if $m=\frac{\height(\blin)-\height(\biso)}{\vv_D(\varpi_K)}$ is an integer, then $(\biso\cdot\varpi_K^m,\blin)$ is an equi-height pair, and $Z_n^\sim(\biso,\blin) = Z_n^\sim(\biso\cdot \varpi_K^m,\blin)$. Therefore, after such a replacement, we could apply the main theorem to these cycles.

It remains to consider the pairs $(\biso,\blin)$ such that $K/F$ is unramified and $\height(\blin)-\height(\biso)\in2\ZZ+1$, which we call odd pairs. We call it an even pair if $\height(\blin)-\height(\biso)\in2\ZZ$. If $(\biso,\blin)$ is an odd pair, we have $Z_n^{(2j)}(\biso,\blin)=0$ and $Z_n^{(2j+1)}(\biso,\blin)\neq0$ for any $j\in\ZZ$. Therefore, the cycle $Z_n(\biso,\blin)^\sim$ can not intersect with other cycles on $\DEF n^{(0)}$. Instead, one needs to consider $\DEF n^{(1)}$. Let $\varpi_F$ be a uniformizer of $D_F$. Using an isomorphism $\map{\beta_n^n(\varpi_F,\id)}{\DEF n^{(1)}}{\DEF n^{(0)}}$, we have
$$
\beta_n^n(\varpi_F,\id)_*\left(Z_n^{(1)}(\biso_1,\blin_1)\litimes_{\DEF n^{(1)}}Z_n^{(1)}(\biso_2,\blin_2)\right)=
Z_n^{(0)}(\varpi_F\cdot\biso_1,\blin_1)\litimes_{\DEF n^{(0)}}Z_n^{(0)}(\varpi_F\cdot\biso_2,\blin_2),
$$
which essentially translates the intersection problems from $\DEF n^{(1)}$ to $\DEF n^{(0)}$. If neither $(\varpi_F\cdot\biso_1,\blin_1)$ nor $(\varpi_F\cdot\biso_2,\blin_2)$ is an odd pair, we can reduce this to the equi-height pairs case and our main theorem applies. The only exception happens if $(\biso_1,\blin_1)$ is an odd pair but $(\biso_2,\blin_2)$ is an even pair when $K_1\cong K_2$ are unramified extensions over $F$. However, these two cycles do not intersect since they are located on different components. Therefore, we can reduce all non-trivial intersection problems of \Heegner cycles to the equi-height pairs case.

\section{Proof of the linear AFL for {$\mathrm{GL}_2(F)$}}\label{EVI}\comments{EVI}
In this section, we prove the linear AFL for the case of $h=1$. Assume $K/F$ is an unramified quadratic extension of non-archimedean local fields of odd residue characteristic.
\subsection{Notation}
Let $H':=\GL_1(F)\times\GL_1(F)=F^\times\times F^\times$ and $G':=\GL_{2}(F)$. We identify $H'$ with a subgroup of $G'$ by the blockwise diagonal embedding. Let $\eta_{K/F}$ be the quadratic character associated to $K/F$. We lift the characters $|\bullet|_F$ and $\eta_{K/F}$ to characters on $H'$ by defining $\eta_{K/F}(x)=\eta_{K/F}(x_1^{-1}x_2)$ and $|x|_F=|x_1^{-1}x_2|_F$ for any $x=(x_1,x_2)\in F^\times\times F^\times$. For any $g\in G'$ and any real-valued test function $f$ on $G'$, we define the orbital integral
\begin{equation}\lb{yazi}
	\OOF\left(g,f,s\right):=\int_{\frac{H^\prime\times H^\prime}{I(g)}}f\left(h_1^{-1}g h_2\right)\eta_{K/F}(h_2)|h_1h_2|_F^s\dd h_1\dd h_2
\end{equation}
as a complex-valued function with $s\in\CC$, where $I(g):=\{(h_1,h_2)\in H'\times H':h_1g=gh_2\}$. 

Fix an embedding $K^\times\subset D_F^\times$ and an element $\ACj\in D_F^\times$.  Denote the invariant polynomial of $\gamma\in D_F^\times$ by $P_\ACj$. Let $g(\ACj)\in G'$ be an element having the same invariant polynomial with $\ACj$. We assume that $P_\ACj$ is irreducible and $P_\ACj(0)P_\ACj(1)\neq 0$. 

Let $H:=\GL_{1}(K)=K^\times$ and $G:=\GL_{2}(F)$. Let $H\subset G$ be an embedding which restricts to $\OO K^\times\subset\GL_{2}(\mathcal O_F)$. Denote by $P_\ACg$ the corresponding invariant polynomial for any $\ACg\in G$. We fix the non-trivial element $\sigma\in\Gal(K/F)$, and denote by $k^\sigma$ the Galois conjugate of $k$.
For any $m\geq n\geq 0$, we denote
$$
g_{m,n} := \[pmatrix]{\pi^m&\\&\pi^n}.
$$
Let $M$ be an integer such that $M>\m(g_{m,n})$.
\subsection{The linear AFL}
Let $Z_0^{\sim}=Z_0^{\sim}(\id,\id)$ be CM cycles on $\DEF 0^\sim$ as in Definition \ref{tcycle}. Recall that $Z_0^{\sim}$ is a collection of cycles on each $\DEF 0^{(n)}$ for any integer $n$.  We are interested in its intersection with cycles 
$${Z'_0}^\sim:=\beta_{0*}^M\beta^{M}_0(\ACj,g_{m,n})^*Z_0^{\sim}(\id,\id)$$
as constructed in Section \ref{heck}. In this expression, the symbol ${Z'_0}^\sim$ is an abbreviation since it depends on the choice of $\ACj,g_{m,n}$ and $M$. Recall that we have defined the function
\[equation]{\lb{yabiui}
f_{R_0g_{m,n}R_0,M}(g)=\[cases]{\frac{\deg(\DEF M\rightarrow\DEF 0)}{\mathrm{Vol}(R_0g_{m,n}R_0)}&\text{ if }g\in R_0g_{m,n}R_0\\0&\text{otherwise}}
}
in \eqref{shenmeya}. For $h=1$, the linear AFL is the following theorem, which we will prove in this section.
\[thm]{\lb{rara}Let $f=f_{R_0g_{m,n}R_0,M}$ be the function as described in \eqref{yabiui}, we have
$$
\chi\left(Z_0^{(0)}\otimes_{\DEF 0} {{Z'_0}^{(0)}}\right) = \left.\pm(2\ln q)^{-1}\frac{\dd}{\dd s}\right|_{s=0}\OOF(g(\ACj),f,s).
$$
Here the sign is chosen so that the right hand side is non-negative.
}
Note that the pair $(\ACj,g_{m,n})$ is not necessarily an equi-height pair. We will show that the conjecture essentially reduces to the case of equi-height pairs with $n=0$. Let
$$\hbar:=\mathrm{Height}(\ACj)-\mathrm{Height}(g_{m,n})$$
be the difference of the height of the two elements. 
 Recall the map $\beta^{M}_0(\ACj,g_{m,n}):\DEF M^\sim\rightarrow\DEF0^\sim$ restricts to 
$$
\map{\beta^{M}_0(\ACj,g_{m,n})}{\DEF M^{(0)}}{\DEF 0^{(-\hbar)}}
$$
as defined in \eqref{haww}, and recall that by using the superscript $(-\hbar)$ on $Z_0^{(-\hbar)}$ we mean the part of $Z_0^\sim$ that is located on $\DEF 0^{(-\hbar)}$ as described in \eqref{zhaogebishengtianhaodelaopo}.
\[lem]{\lb{shability}If $\hbar$ is an odd integer, then
$$
{Z'_0}^{(0)}=0.
$$
If $\hbar$ is an even integer, then $(\pi^{-\frac \hbar2}\ACj,g_{m,n})$ is an equi-height pair, and
$$
{Z'_0}^{(0)} =\beta_{0*}^M\beta^{M}_0(\pi^{-\frac \hbar2}\ACj,g_{m,n})^*Z_0^{(0)}.
$$
}
\[proof]{	Note that
$$
{Z'_0}^{(0)} = \beta_{0*}^M\beta^{M}_0(\ACj,g_{m,n})^*Z_0^{(-\hbar)}.
$$
By definition, $Z_0^{(2i)}\neq 0$ and $Z_0^{(2i+1)}=0$ for all $i\in\ZZ$. If $\hbar$ is an odd integer, then we have
$$
Z_0^{(-\hbar)}=0,
$$
which implies that ${Z'_0}^{(0)}=0$. 

Suppose $\hbar$ is an even number. Since
$$
\mathrm{Height}(g_{m,n}) = -\hbar+\mathrm{Height}(\ACj)=\mathrm{Height}(\pi^{-\frac \hbar2})+\mathrm{Height}(\ACj)=\mathrm{Height}(\pi^{-\frac \hbar2}\ACj),
$$
we know that $(\pi^{-\frac \hbar2}\ACj,g_{m,n})$ is an equi-height pair. Furthermore, since
$$
\beta^{M}_0(\ACj,g_{m,n}) = \beta_0^0(\pi^{\frac \hbar2},\id)\circ \beta^{M}_0(\pi^{-\frac \hbar2}\ACj,g_{m,n})
$$
and $\beta_0^0(\pi^{\frac \hbar2},\id)^*Z_0^{(-\hbar)} = Z_0^{(0)}$, the cycle ${Z'_0}^{(0)}$ can be written as
$$
\beta_{0*}^M\beta^{M}_0(\ACj,g_{m,n})^*Z_0^{(-\hbar)} = \beta_{0*}^M\beta^{M}_0(\pi^{-\frac \hbar2}\ACj,g_{m,n})^*\beta_0^{0*}(\pi^{\frac \hbar2},\id)Z_0^{(-\hbar)}=\beta_{0*}^M\beta^{M}_0(\pi^{-\frac \hbar2}\ACj,g_{m,n})^*Z_0^{(0)}
$$
as desired.}
Since we are in the equi-height situation now, we are able to apply our intersection formula to $\chi\left(Z_0^{(0)}\otimes_{\DEF 0} {{Z'_0}^{(0)}}\right)$.
\[lem]{\lb{aribi}
If $\mathrm{Height}(\ACj)-m\in 2\ZZ$, we have
$$
\chi\left(Z_0^{(0)}\otimes_{\DEF 0} {{Z'_0}^{(0)}}\right) = \frac{\zeta_K(1)^2}{\zeta_F(1)\zeta_F(2)}\int_{\GL_{2}(F)}f_{R_0g_{m,n}R_0,M}(g)\Big|\Res(P_\ACj,P_\ACg)\Big|_F^{-1}\dd g.
$$ 
If $\mathrm{Height}(\ACj)-m\in 2\ZZ+1$, we have
$$
\chi\left(Z_0^{(0)}\otimes_{\DEF 0} {{Z'_0}^{(0)}}\right)=0.
$$
}
\[proof]{If $\mathrm{Height}(\ACj)-\mathrm{Height}(g_{m,n})$ is odd then ${Z'_0}^{(0)}=0$. Therefore $\chi\left(Z_0^{(0)}\otimes_{\DEF 0} {{Z'_0}^{(0)}}\right)=0$.

If $\hbar = \mathrm{Height}(\ACj)-\mathrm{Height}(g_{m,n})$ is even, then by Lemma \ref{shability} we have
$$
\chi\left(Z_0^{(0)}\otimes_{\DEF 0} {{Z'_0}^{(0)}}\right) = \chi\left(Z_0^{(0)}\otimes_{\DEF 0} {\beta_{0*}^M\beta^{M}_0(\pi^{-\frac \hbar2}\ACj,g_{m,n})^*Z_0^{(0)}}\right).
$$
By Theorem \ref{chaochao}, we have
$$
\chi\left(Z_0^{(0)}\otimes_{\DEF 0} {{Z'_0}^{(0)}}\right) = \frac{\zeta_K(1)^2}{\zeta_F(1)\zeta_F(2)}\int_{\GL_{2}(F)}f_{R_0g_{m,n}R_0,M}(g)\Big|\Res(P_{\pi^{-\frac \hbar2}\ACj},P_\ACg)\Big|_F^{-1}\dd g.
$$
Therefore the lemma follows from the fact that $P_{\pi^{-\frac \hbar2}\ACj}=P_{\ACj}$ by Corollary \ref{jiangdiao}.}
Furthermore, it is easy to see that
$$
f_{R_0g_{m,n}R_0,M}(g)=f_{R_0g_{m-n,0}R_0,M}(\pi^{-n}g)
$$
by \eqref{shenmeya}. Since $P_{\pi^{-n}g}=P_g$ by Corolary \ref{jiangdiao},  we have
$$
\int_{\GL_{2}(F)}f_{R_0g_{m,n}R_0,M}(g)\Big|\Res(P_{\ACj},P_\ACg)\Big|_F^{-1}\dd g =\int_{\GL_{2}(F)}f_{R_0g_{m-n,0}R_0,M}(\pi^{-n}g)\Big|\Res(P_{\ACj},P_{\pi^{-n}\ACg})\Big|_F^{-1}\dd g.
$$
By replacing $g\mapsto \pi^{n}g$, the integral can be simplified to
$$
\int_{\GL_{2}(F)}f_{R_0g_{m-n,0}R_0,M}(g)\Big|\Res(P_{\ACj},P_{\ACg})\Big|_F^{-1}\dd g,
$$
which implies that the intersection number only depends on $m-n$.

The similar phenomenon also happens on the orbital integral side. Since $|\pi^n h'|_F=|(\pi^nx_1)^{-1}(\pi^nx_2)|_F=|x_1^{-1}x_2|=|h'|_F$ for any $h'=(x_1,x_2)\in H'$, we have

$${\[split]{
	&\int_{\frac{H^\prime\times H^\prime}{I(g)}}f_{R_0g_{m,n}R_0,M}\left(h_1^{-1}g h_2\right)\eta_{K/F}(h_2)|h_1h_2|_F^s\dd h_1\dd h_2
\\=&\;\int_{\frac{H^\prime\times H^\prime}{I(g)}}f_{R_0g_{m-n,0}R_0,M}\left((\pi^{m}h_1)^{-1}g h_2\right)\eta_{K/F}(h_2)|(\pi^mh_1)h_2|_F^s\dd h_1\dd h_2
\\=&\;\int_{\frac{H^\prime\times H^\prime}{I(g)}}f_{R_0g_{m-n,0}R_0,M}\left(h_1^{-1}g h_2\right)\eta_{K/F}(h_2)|h_1h_2|_F^s\dd h_1\dd h_2.
}}$$
This implies that we only need to prove the linear AFL for $f(x)=f_{R_0g_{m,0}R_0,M}$. Since $\mathrm{Height}(g_{m,0})=m$, we have
$$
\hbar = \mathrm{Height}(\ACj) - \mathrm{Height}(g_{m,0})  =  \mathrm{Height}(\ACj) - m.
$$ 
The intersection number is given by Lemma \ref{aribi} when $\hbar$ is even. The intersection number is $0$ when $\hbar$ is odd by Lemma \ref{shability}. To simplify our notations, let
$$
T_{m,0}(g):=\[cases]{1&\text{ if }g\in R_0g_{m,0}R_0;\\0&\text{ otherwise. }}
$$
Since $f_{R_0g_{m,0}R_0,M}$ is a scalar multiple of $T_{m,0}$, the theorem reduces to the following theorem.
\begin{thm}\lb{buxiangbuxiang}
	Suppose $P_\gamma(x)$ is the invariant polynomial of $g(\ACj)$ and $P_\gamma(0)P_\gamma(1)\neq 0$. For any $m$, if $\mathrm{Height}(\ACj)-m\in 2\ZZ$, then we have 
	\begin{equation}\lb{oi}
		\left.\pm(2\ln q)^{-1}\frac{\dd}{\dd s}\right|_{s=0}\OOF(g(\ACj),T_{m,0},s)=\frac{\zeta_K(1)^2}{\zeta_F(1)\zeta_F(2)}\int_{\GL_{2}(F)}T_{m,0}(g)\Big|\Res(P_\ACj,P_\ACg)\Big|_F^{-1}\dd g.
\end{equation}
	Otherwise, if $\mathrm{Height}(\ACj)-m\in 2\ZZ+1$, then the left hand side of the above equation is $0$. Here the sign is chosen so that the left hand side is non-negative.
\end{thm}
Call the left hand side of this identity the analytic side, and the right hand side the geometric side.
For any two sets $A$ and $B$, by $A-B$ we mean the subset of $A$ consisting of all elements that is not in $B$. Note that 
\begin{equation}\lb{miaoya}
T_{m,0}(g)=1 \iff \det g\in\pi^m\OO F^\times\text{ and }g\in\gl_2(\OO F)-\pi\gl_2(\OO F).
\end{equation}
  
\subsection{Matching orbits}
Fix the decomposition $\ACj = \ACj_++\ACj_-$ according to the embedding $K\subset D_F$ such that $\ACj_+ k=k\ACj_+$ and $\ACj_- k=k^\sigma\ACj_-$ for all $k\in K$. Note that when $h=1$, $D_F$ is a quaternion algebra over $F$. Let $\vv_D(\bullet)$, $\vv_K(\bullet)$ and $\vv_F(\bullet)$ be normalized valuations on $D_F$, $K$ and $F$. It is well-known that $\vv_D(x)=2\vv_F(x)=2\vv_K(x)$ for any $x\in F\subset K\subset D_F$. Since $P_\ACj$ is irreducible and $P_\ACj(1)P_\ACj(0)\neq 0$, there exists an element $c\in F$ such that $c\neq 0,1$ and $P_\ACj(X)=X-(1-c)^{-1}$.
Note that $P_\ACj$ is the characteristic polynomial of $\ACj_+(\ACj_++\ACj_-)^{-1}\ACj_+(\ACj_+-\ACj_-)^{-1}$ as an element in $K$. Therefore,
$$
\ACj_+(\ACj_++\ACj_-)^{-1}\ACj_+(\ACj_+-\ACj_-)^{-1}=(1-c)^{-1}\in F^\times,
$$
in particular, we have $\ACj_+\neq 0$. Note that this element can also be written as
$$
\ACj_+(\ACj_++\ACj_-)^{-1}\ACj_+(\ACj_+-\ACj_-)^{-1}=(1-\ACj_+^{-1}\ACj_-\ACj_+^{-1}\ACj_-)^{-1},
$$
which implies that $\ACj_+^{-1}\ACj_-\ACj_+^{-1}\ACj_-=c$.

\[lem]{\lb{commute}If $\gamma\in D_F^\times$ such that $\gamma x=x^\sigma\gamma$ for any $x\in K$, then $\vv_D(\gamma)$ is an odd integer. }

\[proof]{Since $\gamma x=x^\sigma \gamma$ for any $x\in K$, one sees that $\gamma^2$ commutes with all elements in $K$, which implies $\gamma^2\in K$. Since $\gamma x=x^\sigma \gamma$ for any $x\in K$, one sees that $\gamma(\gamma^2)=(\gamma^2)^\sigma\gamma$, which implies $(\gamma^2)=(\gamma^2)^\sigma$ and therefore $\gamma^2\in F$. Assume, for the sake of contradiction, that $\vv_D(\gamma)$ is even. Then $\vv_F(\gamma^2)=\frac12\vv_D(\gamma^2)=\vv_D(\gamma)$ is also an even integer. Since the norm subgroup of the unramified extension $K/F$ is $\pi^{2\ZZ}\cdot\OO K^\times$, there exists some element $y\in K$ such that $\gamma^2=yy^\sigma$. Since $\gamma y=y^\sigma\gamma$, we have
$$
(\gamma+y^\sigma)(\gamma-y) = \gamma^2 - yy^\sigma = 0.
$$
Since $D_F$ is a division algebra, we have either $\gamma=y$ or $\gamma=-y^\sigma$, which contradicts the assumption that $\gamma x=x^\sigma\gamma$ for any $x\in K$.  
}

Note that $c=(\ACj_+^{-1}\ACj_-)^2$, and that $\ACj_+^{-1}\ACj_-k=k^\sigma\ACj_+^{-1}\ACj_-$ for any $k\in K$. We have $\vv_F(c)\in 2\ZZ+1$. Similarly, we have  $\vv_D(\ACj_-)\in 2\ZZ+1$. Hence $\vv_D(\ACj_+)\in 2\ZZ$. Therefore $\vv_D(\ACj)=\vv_D(\ACj_++\ACj_-)=\min\{\vv_D(\ACj_+),\vv_D(\ACj_-)\}$. Moreover,
\[equation]{\lb{yijisuan}\[split]{\vv_D(\ACj)\in 2\ZZ\iff \vv_D(\ACj_+)>\vv_D(\ACj_-) \iff \vv_F(c)>0,\\ 
\vv_D(\ACj)\in 2\ZZ+1\iff \vv_D(\ACj_+)<\vv_D(\ACj_-) \iff \vv_F(c)<0.}}

\subsection{The calculation for the analytic side} 
Recall that $H'\subset G'\cong \GL_2(F)$ is the subgroup of diagonal matrices. For any $g\in G'$, the decomposition $g=g_++g_-$ gives a diagonal matrix $g_+$ and a matrix $g_-$ with zero entries on the diagonal.

Note that 
$g(\ACj)$ is an element with invariant polynomial $P_\ACj(X)=X-(1-c)^{-1}$. Since $c\neq 1$ or $0$, we know $(g(\ACj)_+^{-1}g(\ACj)_-)^2=c$. Since the diagonal entries of $g(\ACj)_+^{-1}g(\ACj)_-$ are $0$, there exists a diagonal matrix $\widetilde{h}$ such that
$$
\widetilde{h}^{-1}g(\ACj)_+^{-1}g(\ACj)_-\widetilde{h}=\mm0c10.
$$
Therefore
$$
\left(\widetilde{h}^{-1}g(\ACj)_+^{-1}\right)\cdot g(\ACj)\cdot \widetilde{h}=\mm1c11.
$$
 Since $g(\ACj)_+\in H'$ and $\widetilde{h}\in H'$, we have $\widetilde{h}^{-1}g(\ACj)_+^{-1}\in H'$ and $\widetilde{h}\in H'$. Since replacing $g(\ACj)$ by $h'g(\ACj)h''$ does not change the absolute value of the orbital integral \eqref{yazi} for any $h',h''\in H'$, we can assume 
$$
g(\ACj):=\mm1c11\in G'.
$$
Let $k$ be the integer such that $\vv_F(c)=2k+1$. Indeed, the invariant polynomial of $g(\ACj)$ is $P_\ACj(X)=X-(1-c)^{-1}$.  

\begin{prop}\lb{ayoha}
	Let $M=m-\min\{0,\vv_F(c)\}$. We have 
\begin{equation*}(1+q^{2s})\OOF(g(\ACj),T_{m,0},s) =
	\[cases]{0&\text{ if }M\in 2\ZZ+1
							;\\
	{-(1-q^{-2s})\left(q^{-Ms}+q^{(M+4k-2)s}\right)}
	&\text{ if }M\in 2\ZZ\smallsetminus\{0\}
	;\\
	{-q^{2\vv_F(c)s}+q^{-2s}}
	&\text{ if }M=0
	.}
\end{equation*}
\end{prop}
\begin{proof}
	Recall the definition of $\OOF(g(\gamma),T_{m,0},s)$ in \eqref{yazi}. Note that $I(g)\cong F^\times$ is the set of scalar matrices, which embeds diagonally in $H'\times H'\cong (F^\times\times F^\times)\times (F^\times\times F^\times)$. Therefore $H'\times H'/I(g)\cong F^\times\times F^\times\times F^\times$, with its Haar measure normalized by $\OO F^\times\times \OO F^\times\times \OO F^\times$.
Then, the orbital integral $\OOF(g(\gamma),T_{m,0},s)$ equals
\begin{equation}
	\int_{F^\times\times F^\times\times F^\times}T_{m,0}\left(\mm a{}{}b\mm1c11\mm{d}{}{}1\right)\eta_{K/F}(d)|a^{-1}bd|_F^{-s}\dd a^\times\dd b^\times\dd d^\times.
\end{equation}
Letting $x:=\vv_F(a)$, $y:=\vv_F(b)$ and $z:=\vv_F(d)$, we rewrite the above integral into
$$
\OOF(g(\gamma),T_{m,0},s)=\int_{F^\times\times F^\times\times F^\times}T_{m,0}\mm{ad}{ac}{bd}b(-1)^{z}q^{(y+z-x)s}\dd a^\times\dd b^\times\dd d^\times,
$$
which equals
$$
	\sum_{x,y,z\in\ZZ}\int_{\pi^x\OO F^\times\times \pi^y\OO F^\times\times \pi^z\OO F^\times}T_{m,0}\mm{ad}{ac}{bd}b(-1)^{z}q^{(y+z-x)s}\dd a^\times\dd b^\times\dd d^\times.
$$
	Furthermore let $w:=\vv_F(\det(g(\ACj)))=\vv_F(1-c)$. Since $\vv_F(c)=2k+1\neq 0$, we have $w=\min\{0,2k+1\}$. Then $M=m-w$. 
Note that each summand is either $(-1)^{z}q^{(y+z-x)s}$ or $0$, where it is non-zero if and only if 
$$
T_{m,0}\mm{ad}{ac}{bd}b=1,
$$
which, by \eqref{miaoya}, is equivalent to
$$
x+y+z=M\qquad\text{ and }\qquad\min\{x+z,x+(2k+1),y+z,y\}=0.
$$
Replacing $z$ by $M-x-y$, this restriction can be simplified to
	\[equation]{\lb{symm}
	z=M-x-y\qquad \text{ and }\qquad \min\{M-y,x+(2k+1),M-x,y\}=0.
	}
	Let $S_{M,k}=\{(x,y): \min\{M-y,x+(2k+1),M-x,y\}=0\}$. Then we have
	\[equation]{\lb{obt}
	\OOF(g(\gamma),T_{m,0},s)=\sum_{(x,y)\in S_{M,k}}(-1)^{M-x-y}q^{(M-2x)s}.
	}

From now on, we will discuss the orbital integral in two cases depending on the parity of $M$. Let us start with the case of $M\in 2\ZZ+1$. Since the set $S_{M,k}$ is invariant under the replacement $y\mapsto M-y$, we have
$$
	\sum_{(x,y)\in S_{M,k}}(-1)^{M-x-y}q^{(M-2x)s} = \sum_{(x,y)\in S_{M,k}}(-1)^{y-x}q^{(M-2x)s} = -\sum_{(x,y)\in S_{M,k}}(-1)^{M-x-y}q^{(M-2x)s},
$$
which implies $\OOF(g(\gamma),T_{m,0},s)=0$. 

In the case of $M\in 2\ZZ$, we decompose $S_{M,k}$ into the following subsets
	$${\[split]{
		S_{M,k}=\;&\{(-(2k+1),y):0\leq y< M\}\cup\{(M,y):0< y\leq M\}\\
		&\cup\{(x,0):-(2k+1)<x\leq M\}\cup\{(x,M):-(2k+1)\leq x<M\}.
		}}
$$
It suffices to calculate \eqref{obt} by each component. For the first two components $\{(-(2k+1),y):0\leq y< M\}$ and $\{(M,y):0< y\leq M\}$, since $M$ is even, we have
$$
	\sum_{y=1}^{M}(-1)^{M-x-y}q^{(M-2x)s} = 0 = \sum_{y=0}^{M-1}(-1)^{M-x-y}q^{(M-2x)s}.
$$
Therefore,
$$
	\OOF(g(\gamma),T_{m,0},s) = \sum_{x=-2k}^M (-1)^{M-x} q^{(M-2x)s} + \sum_{x=-2k-1}^{M-1} (-1)^{-x} q^{(M-2x)s},
$$
which equals
$$
	\frac{q^{(M+4k)s}+q^{(-M-2)s}}{1+q^{-2s}}+\frac{-q^{(M+4k+2)s}-q^{-Ms}}{1+q^{-2s}}=-\frac{(1-q^{-2s})\left(q^{-Ms}+q^{(M+4k+2)s}\right)}{1+q^{-2s}}.
$$

Suppose 
$M=0$. We have
$$
	S_{0,k}=\{(x,0):-(2k+1)\leq x\leq 0\}.
$$
Therefore
$$
	\OOF(g(\gamma),T_{m,0},s) = \sum_{x=-(2k+1)}^0 (-1)^{-x} q^{-2xs}=\frac{-q^{(4k+2)s}+q^{-2s}}{1+q^{2s}}.
$$
Combining all previous cases, we complete the proof of the proposition.
\end{proof}
\begin{proof}[Proof of Theorem \ref{buxiangbuxiang} in the case that $m-\mathrm{Height}(\gamma)\in 2\ZZ+1$]
	Note that $\mathrm{Height}(\gamma)=\vv_D(\gamma)$. By \eqref{yijisuan}, $\vv_D(\gamma)-\min\{\vv_F(c),0\}\in 2\ZZ$. This implies that 
	$$m-\min\{0,\vv_F(c)\}\in 2\ZZ+1\iff m-\mathrm{Height}(\gamma)\in 2\ZZ+1.$$
	In this case, by Proposition \ref{ayoha}, the orbital integral is constantly zero as desired.
\end{proof}
We consider the remaining cases. 
Note that $m-\mathrm{Height}(\gamma)\in 2\ZZ$ implies that 
$$m-\min\{0,\vv_F(c)\}\in 2\ZZ\iff m\in 2\ZZ\text{ and }\vv_F(c)>0, \text{ or }m\in 2\ZZ+1\text{ and }\vv_F(c)<0.$$
Therefore, when $m\in 2\ZZ\smallsetminus\{0\}$ and $\vv_F(c)>0$, or when $m\in 2\ZZ+1$ and $\vv_F(c)<0$, we have
\[equation]{\lb{caseone}
-\frac1{2\ln q}\left.\frac{\dd}{\dd s}\right|_{s=0}\OOF(g(\gamma),T_{m,0},s) = 1;
}
when $m=0$ and $\vv_F(c)>0$, we have 
\[equation]{\lb{casetwo}
-\frac1{2\ln q}\left.\frac{\dd}{\dd s}\right|_{s=0}\OOF(g(\gamma),T_{m,0},s) = \frac{\vv_F(c)+1}2.
}

\subsection{The calculation for the geometric side}
In this subsection, we calculate the geometric side of \eqref{yazi} in cases of \eqref{caseone} and \eqref{casetwo} respectively. 

Recall that $H:=\GL_1(K)=K^\times$ and $G:=\GL_{2}(F)$ are algebraic groups over $F$. Fix an embedding $H\subset G$ that restricts to $\OO K^\times\subset\GL_2(\OO F)$. Let $0_2$ be the zero matrix in $\matt22 F$. Note that 
$$\left(H\cap\matt22{\OO F}\right)\cup\{0_{2}\}\cong\OO K.$$
By abuse of notation, we also denote $\left(H\cap\matt22{\OO F}\right)\cup\{0_{2}\}$ by $\OO K$ if no confusion can arise. Fix the decomposition $\ACg = \ACg_++\ACg_-$ according to the embedding $H\subset G$ such that $\ACg_+ k=k\ACg_+$ and $\ACg_- k=k^\sigma\ACg_-$ for all $k\in H\subset G$. Let 
$$P_g(X)=X-b$$
be the invariant polynomial of $g$. Note that $b\in F$ and $b=\ACg_+(\ACg_+-\ACg_-)^{-1}\ACg_+(\ACg_++\ACg_-)^{-1}$. Therefore as an element in $\GL_{2}(F)$ we have $\det_F(b)=b^2$. On the other hand, by \eqref{eeu}, we have $\det(b)=\det(\ACg_+)^2\det(g)^{-2}$, which implies that
\[equation]{\lb{tiaoban}
b = \ACg_+(\ACg_+-\ACg_-)^{-1}\ACg_+(\ACg_++\ACg_-)^{-1} = \pm (\ACg_+\ACg_+^\sigma)\det(g)^{-1}\in F.
}
Note that $\zeta_K(s)=(1-q^{-2s})^{-1}$ and $\zeta_F(s)=(1-q^{-s})^{-1}$. Recall that $P_\ACj(X)=X-(1-c)^{-1}$. We have $\Res(P_\ACj,P_\ACg) = (1-c)^{-1}-b$. Therefore the geometric side of \eqref{oi} equals
\[equation]{\lb{shazicailianai}
\frac{1}{1+q^{-1}}\int_{\GL_{2}(F)}T_{m,0}(g)\left|(1-c)^{-1}-b\right|_F^{-1}\dd g.
}

\begin{lem}\lb{intega}
	If $T_{m,0}(g)=1$, then $g_+,g_-\in\matt22{\OO F}$. Furthermore, if $m>0$, then $g_+\in\OO K^\times$.
\end{lem}
\begin{proof}
	By definition \eqref{miaoya}, $T_{m,0}(g)=1$ implies that $g\in\matt 22F$ and $g\not\equiv 0$ modulo $\pi$. Let $\mu\in \OO K^\times$ be an element such that $\mu^\sigma=-\mu$. We have $2\mu\cdot g_+=\mu g+g\mu$ and $2\mu\cdot g_-=\mu g-g\mu$, which implies $g_+,g_-\in\matt 22{\OO F}$. Suppose $m>0$, and suppose for the sake of contradiction that $g_+\in\pi\OO K$. Then $g_+\equiv 0$ modulo $\pi$, which implies that $g_-=g-g_+\not\equiv 0$ modulo $\pi$. Note that there is an element $\sigma\in\GL_{2}(\OO F)$ such that $\sigma k\sigma^{-1}=k^\sigma$. Then $\sigma \cdot g_- \in K$, which implies $\sigma\cdot g_-\in\OO K^\times\subset\GL_2(\OO F)$. Therefore $g=g_++g_-\in\GL_2(\OO F)$, which contradicts to $T_{m,0}(g)=1$ since $m>0$.
\end{proof}

\subsubsection{Case 1: $m>0$}
By Lemma \ref{intega} and \eqref{tiaoban}, one sees that $T_{m,0}(g)=1$ implies $|b|_F=|\deg(g)^{-1}|_F=q^{m}>1$. Recall that $\vv_F(c)$ is an odd integer, and, in particular, $\vv_F(c)\neq 0$, which implies that $(1-c)^{-1}\in\OO F$. Therefore $|(1-c)^{-1}-b|_F=q^m$, which implies that the equation \eqref{shazicailianai} equals
$$
\frac1{1+q^{-1}}\int_{\GL_{2}(F)}T_{m,0}(g)q^{-m}\dd g.
$$
Note that under the Haar measure $\dd g$ normalized by $\GL_2(\OO F)$, we have
$$
\Vol\left(\GL_2(\OO F)\[pmatrix]{1&\\&\pi^m}\GL_2(\OO F)\right) = (1+q^{-1})q^{m},
$$
which implies that the equation \eqref{shazicailianai} equals $1$ as desired. Comparing this result with \eqref{caseone}, we complete the proof of the linear AFL for the case of $m>0$.

\subsubsection{Case 2: $m=0$} In this case, we only need to prove the formula for $\vv_F(c)>0$. Since the volume of the subset $\{g\in G:g_+=0\}$ is zero, we may assume that $\ACg_+$ is invertible for integration purposes.
Recall $P_\ACj(X)=X-(1-c)^{-1}$, and $P_\ACg(X)=X-(1-\ACg_+^{-1}\ACg_-\ACg_+^{-1}\ACg_-)^{-1}$. We have
\[equation*]{
\Res(P_\ACj,P_\ACg) = \frac{c - \ACg_+^{-1}\ACg_-\ACg_+^{-1}\ACg_-}{(1-c)(1-\ACg_+^{-1}\ACg_-\ACg_+^{-1}\ACg_-)}.
}
If $g\in\GL_2(\OO F)$, by \eqref{tiaoban}, we have $|(1-\ACg_+^{-1}\ACg_-\ACg_+^{-1}\ACg_-)^{-1}|_F=|g_+g_+^\sigma|_F$. Since $\ACg_-\ACg_+^{-1}=(\ACg_+^{-1})^\sigma\ACg_-$, we have $\ACg_+^{-1}\ACg_-\ACg_+^{-1}\ACg_-=\ACg_+^{-1}(\ACg_+^{-1})^\sigma\ACg_-^2$. Note that $\vv_F(c)>0$ implies $|1-c|_F=1$. Therefore,
\[equation]{\lb{resres}
|\Res(P_\ACj,P_\ACg)|_F = |\ACg_+\ACg_+^\sigma c - \ACg_-^2|_F.
}
Note that the embedding $\OO K^\times\subset\GL_2(\OO F)$ gives rise to an $\OO K$-action on $\OO F^2$, which induces an isomorphism $\OO F^2\rightarrow \OO K$. The non-trivial Galois conjugation of $\OO K$ corresponds to an element $\sigma\in\GL_2(\OO F)$ such that $\sigma k=k^\sigma\sigma$ for any $k\in H$. Therefore, for any $g\in\GL_2(\OO F)$, there exists $k_1,k_2\in H$ such that 
$$
g_+=k_1,\qquad g_-=k_2\sigma.
$$
For any $n\geq 1$, let 
\[equation]{\lb{omama}
\Omega_n:=\{g=k_1+k_2\sigma\in\GL_2(\OO F): k_2\in\pi^n\OO K\}.
}
Recall the decomposition \eqref{decompa}. In particular, we have $\pi^n\matt 22{\OO F}=\pi^n\OO K\oplus \pi^n\sigma\OO K$ for any $n\geq 1$, which implies that 
$$
\Omega_n = \OO K^\times + \pi^n\matt 22{\OO F}.
$$
Under the normalized Haar measure of $\GL_2(\OO F)$, the volume of $\Omega_n$ is 
\[equation]{\lb{onimei}
\Vol(\Omega_n)=\frac{\#(\OO K/\pi^n)^\times}{\#(\GL_2(\OO F/\pi^n))} = \frac{q^{-2n}}{1-q^{-1}}.
}
Therefore, note that $\ACg_+\ACg_+^\sigma c - \ACg_-^2\in\OO F^\times$ for any $g\notin\Omega_1$, we have
$$
\int_{\GL_2(\OO F) - \Omega_1}|\Res(P_\ACj,P_\ACg)|_F^{-1}\dd g = \int_{\GL_2(\OO F) - \Omega_1}\dd g = \Vol(\GL_2(\OO F))-\Vol(\Omega_1) = 1-\frac{q^{-2}}{1-q^{-1}}.
$$
For any $1\leq 2n<\vv_F(c)$, we have  $|\ACg_+\ACg_+^\sigma c - \ACg_-^2|_F=|\ACg_-^2|_F=q^{-2n}$ for any $g\in\Omega_n-\Omega_{n+1}$, which implies that
$$
\int_{\Omega_n - \Omega_{n+1}}|\Res(P_\ACj,P_\ACg)|_F^{-1}\dd g=\int_{\Omega_n - \Omega_{n+1}}|g_-^2|_F^{-1}\dd g=q^{2n}\left(\Vol(\Omega_n)-\Vol(\Omega_{n+1})\right)=1+q^{-1}.
$$
For $2n=\vv_F(c)+1$, we have  $|\ACg_+\ACg_+^\sigma c - \ACg_-^2|_F=|c|_F=q^{-\vv_F(c)}$ for any $g\in\Omega_n$, which implies that
$$\int_{\Omega_n}|\Res(P_\ACj,P_\ACg)|_F^{-1}\dd g = \int_{\Omega_n}|c|_F^{-1}=q^{\vv_F(c)}\Vol(\Omega_n)=\frac{q^{\vv_F(c)-2n}}{1-q^{-1}}=\frac{q^{-1}}{1-q^{-1}}.$$
Combining preceding equations, we have
$$
\int_{\GL_2(\OO F)}|\Res(P_\ACj,P_\ACg)|_F^{-1}\dd g = 1-\frac{q^{-2}}{1-q^{-1}}+\sum_{n=1}^{\frac{\vv_F-1}2}(1+q^{-1}) +\frac{q^{-1}}{1-q^{-1}} = \frac{(1+q^{-1})(\vv_F(c)+1)}2,
$$
which proves that
$$
\frac1{1+q^{-1}}\int_{\GL_2(\OO F)}|\Res(P_\ACj,P_\ACg)|_F^{-1}\dd g= \frac{\vv_F(c)+1}2
$$
as desired. Therefore we proved the linear AFL for the case of $m=0$. Combining all cases, we complete the proof of the linear AFL for $h=1$.

\bibliography{bib}
\bibliographystyle{amsalpha}
\end{document}